\documentclass[12pt]{article}
\usepackage{authblk}
\usepackage{graphicx}
\usepackage{subfig}
\usepackage{multirow}
\usepackage{amsfonts}
\usepackage[a4paper,left=2.4cm,right=2.3cm,top=2cm,bottom=1.8cm]{geometry}%
\usepackage{amsmath,amssymb,mathrsfs}
\usepackage{amsthm}
\usepackage{textcomp}
\usepackage{color}

\usepackage{hyperref}
\usepackage{float}
\usepackage{epsfig}
\usepackage{epstopdf}
\usepackage{array,booktabs}
\usepackage{tabularx}
\usepackage{cleveref}

\usepackage{algorithm}
\usepackage[noend]{algcompatible}

\usepackage{longtable}
\usepackage{orcidlink}

\usepackage{comment}
\usepackage{bm}

\newtheorem{proposition}{Proposition}
\newtheorem{lemma}{Lemma}
\newtheorem{theorem}{Theorem}
\newtheorem{remark}{Remark}

\newtheorem{assumption}{Assumption}
\newtheorem{corollary}{Corollary}

\def\keywords{\vspace{.5em}
{\textit{Keywords}:\,\relax%
}}

\begin{document}

\title{Power Homotopy for Zeroth-Order Non-Convex Optimizations} 

\author[1]{Chen Xu}
\affil[1]{Department of Engineering, Shenzhen MSU-BIT University, China. 
\authorcr Email: chen\_xu\_research@outlook.com
}





\maketitle

\begin{abstract}
The existing method of GS-PowerOpt solves the non-convex optimization problem of the form $\max_{\boldsymbol{x} \in \mathbb{R}^d} f(\boldsymbol{x})$ through maximizing a Gaussian-smoothed surrogate $F_{N,\sigma}(\boldsymbol{\mu}) = \mathbb{E}_{\boldsymbol{x}\sim\mathcal{N}(\boldsymbol{\mu},\sigma^2 I_d)}[e^{N f(\boldsymbol{x})}]$. We analyze the role of the smoothing radius $\sigma>0$ and identify a limitation of the fixed-$\sigma$ design used in GS-PowerOpt. Specifically, $\sigma$ induces an inherent exploration--refinement tradeoff: a larger $\sigma$ improves global exploration and finite-time surrogate optimization, but may distort the location of the surrogate maximizer; in contrast, a smaller $\sigma$ better preserves local structure but can weaken gradient signals away from high-value regions.

To address this limitation, we propose GS-PowerHP, a power-smoothed homotopy method with an incrementally decaying $\sigma$ schedule. The proposed mechanism uses larger smoothing radii in early iterations to maintain informative gradient signals when the iterate is far from high-value regions, and gradually decreases $\sigma$ to improve local refinement near the maximizer. We provide theoretical results showing that this decaying schedule improves the exploration--refinement tradeoff of fixed-$\sigma$ power smoothing. Empirically, GS-PowerHP consistently outperforms the fixed-$\sigma$ baseline and exhibits robust performance across different optimization tasks, including adversarial attacks on ImageNet ($d=150{,}528$), where it substantially improves over other smoothing-based zeroth-order methods.

\end{abstract}

\section{Introduction}
Zeroth-order (ZO) optimization methods optimize an objective function $f:\mathbb{R}^d\rightarrow \mathbb{R}$ without requiring its gradient. This makes them particularly useful for non-differentiable or non-convex problems, which are prevalent in machine learning and computer vision. 

A powerful category of zeroth-order optimization methods, namely homotopy for optimization \cite{MobahiFisher2015, Hazan2016, GaoSener2022, Lin2023},  are characterized by a surrogate objective, which is constructed as the Gaussian smoothed transform of the original:  $F_{\sigma}(\bm{\mu})=\mathbb{E}_{\bm{x}\sim\mathcal{N}(\bm{\mu},\sigma^2I_d)}[f(\bm{x})]$, where $\mathcal{N}$ denotes a Gaussian distribution, $I_d$ denotes an identity matrix of shape $d\times d$, and $\sigma>0$ is called the smoothing radius (or scaling parameter). This transformation smooths the landscape of 
$F$ compared to the original $f$, which can remove sharp local minima and facilitate locating the global optimum. The gradient of $F_{\sigma}(\bm{\mu})$ can be efficiently estimated (e.g., using the method in \cite{Nesterov2017}) with $f$'s value on points randomly sampled near $\bm{\mu}$, which enables stochastic gradient methods for optimizing $F_{\sigma}$. However, the optimum point $\bm{\mu}^*$ of $F_\sigma$ is in general away from the global optimum $\bm{x}^*$ of $f$, unless $\sigma$ is close to 0. Therefore, the standard homotopy applies a double-loop mechanism, where the outer loop incrementally decreases the smoothing radius $\sigma$ and the inner loop performs a stochastic gradient algorithm to find the maximizer of $F_{\sigma}$ under the current $\sigma$ value \cite{Hazan2016}.

The double-loop mechanism is time-consuming. The iteration complexity for the standard ZO homotopy is $O(d^{2}\varepsilon^{-4})$, as shown by \cite[Theorem 5.1]{Hazan2016}. \cite{Iwakiri2022} proposes a more efficient single-loop method (ZOSLGH) that updates the solution candidate and $\sigma$ at the same iteration, with its iteration complexity reduced to $O(d^2\varepsilon^{-2})$. However, ZOSLGH is only theoretically guaranteed to converge a local optimum of $f$. 

In our previous work \cite{GS-PowerOpt}, we proposed GS-PowerOpt, a smoothing-based method that power-transforms the objective before Gaussian smoothing, rather than adopting a $\sigma$-decreasing schedule. Specifically, it solves $\max_{\bm{x}\in\mathbb{R}^d}f(\bm{x})$ through targetting at the surrogate problem of \begin{equation}
\label{gs-poweropt-obj}
\max_{\bm{\mu}}F_{N,\sigma}(\bm{\mu}):=\mathbb{E}_{\bm{x}\sim\mathcal{N}(\bm{\mu},\sigma^2 I_d)}[e^{Nf(\bm{x})}]. 
\end{equation}
As illustrated in their Section 2.1, increasing the power will drag the surrogate's global optimum to $\bm{x}^*:=\arg\max_{\bm{x}}f(\bm{x})$. In fact, we proved that with a sufficiently large power $N$, all the stationary points of the surrogate objective fall in a small neighborhood of the global optimum point $\bm{x}^*$ of $f$, provided that $f$ has a unique global optimum.

While GS-PowerOpt demonstrates strong empirical performance, we show in our Section \ref{sec-motivation} that its reliance on a fixed smoothing parameter $\sigma$ introduces an inherent exploration–accuracy trade-off: a large $\sigma$ facilitates global exploration but can distort the location of the maximizer of the surrogate objective, whereas a small $\sigma$ preserves local geometry but yields weak gradient signals when far from optimal regions. We formalize this limitation by analyzing the dependence of iteration complexity and surrogate bias on $\sigma$. Motivated by this structural tension, we propose Power-Transformed Gaussian Homotopy (GS-PowerHP), which replaces the fixed-$\sigma$ mechanism with a principled decaying-$\sigma$ schedule. This strategy dynamically transitions from global exploration in early iterations to local refinement near the maximizer in later stages, and is directly derived from the geometry of the power-transformed Gaussian surrogate, supporting both theoretical convergence guarantees and strong empirical performance. 


\textbf{Contributions.} Our contributions are fourfold. 
\begin{enumerate}
\item We provide new theoretical insights into the behavior of the existing powerful method of GS-PowerOpt, identifying a limitation of the fixed-$\sigma$ mechanism due to the exploration-refinement trade-off.
\item Building on this understanding, we propose GS-PowerHP, a zeroth-order algorithm that integrates power transformation with a $\sigma$-decaying mechanism. To the best of our knowledge, such a combination has not been explicitly explored in prior zeroth-order methods. Empirically, this design substantially enhances performance on non-convex optimization tasks (Subsection \ref{sec-ablation-study}). 
\item We establish theoretical convergence guarantees for GS-PowerHP; see Corollaries~\ref{convergence-rate} and~\ref{summary-cor}. In particular, our finite-time convergence rate to first-order stationarity in Corollary \ref{convergence-rate} validates how the decaying-$\sigma$ schedule improves the fixed-$\sigma$ tradeoff: larger smoothing radii can be exploited during finite-time optimization for faster progress, while the smoothing radius is gradually reduced to improve local refinement.

\item We demonstrate consistent empirical superiority over existing smoothing-based zeroth-order optimization methods, achieving the strongest overall performance on ImageNet under challenging least-likely targeted black-box attack settings (Section \ref{sec-experiments}). 


\end{enumerate}

\textbf{Novelty}. The main novelty of our work lies in the following three aspects:
\begin{enumerate}
\item a theoretical identification of a fixed-$\sigma$ limitation specific to power-transformed Gaussian smoothing (Section \ref{sec-motivation}); 
\item the combination of a decaying $\sigma$ and the power transformed smoothing (Algorithm \ref{alg:gspowerhp});
\item and the theoretical and empirical evidence that the combination solves the limitation while retaining the global-optimum-neighborhood behavior (Section \ref{sec-convergence-analysis} and \ref{sec-experiments}).
\end{enumerate}

\section{Related Works} 
\label{sec-related}
The key distinction between our method, GS-PowerHP, and GS-PowerOpt \cite{GS-PowerOpt} is that GS-PowerHP adopts a decaying smoothing radius $\sigma$. This design is motivated by the trade-off between iteration complexity and solution accuracy, as discussed in the analytical analysis in Section \ref{sec-motivation}. Consistent with this motivation, the experiments in Subsection \ref{sample-efficiency-exp} show that, under identical per-iteration sampling budgets, GS-PowerHP requires fewer iterations and attains better solutions than GS-PowerOpt. Taken together, the theoretical analysis and empirical results suggest that, in power-based global optimization, the smoothing radius can be more effectively treated as an adaptive computational resource rather than as a static hyper-parameter.



There are a few other studies that also incorporates the power-transformed objective with Gaussian smoothing \cite{Dvijotham2014, Roulet2020, Chen2024, LiBeirami2023, aminian2025}. However, the first three studies do not establish a formal relationship between the power size and the geometric shift--the distance between the global optimum of the original objective function ($f$) and that of the smoothed surrogate $F_{N,\sigma}(\bm{\mu})$, while the objectives of the last two works are not to optimize the original objective. Furthermore, their methodologies do not contain an incremental $\sigma$-decreasing mechanism.

\section{Analysis on the Role of $\sigma$ in GS-PowerOpt}
\label{sec-motivation}
\subsection{The trade-off between Faster Exploration and Local Refinement}
GS-PowerOpt optimizes the power-transformed Gaussian-smoothed surrogate
\[
F_{N,\sigma}(\boldsymbol{\mu})
=
\mathbb{E}_{\boldsymbol{x}\sim\mathcal{N}(\boldsymbol{\mu},\sigma^2I_d)}
\left[e^{Nf(\boldsymbol{x})}\right].
\]
We show that, for this surrogate objective, the smoothing radius $\sigma$ plays two competing roles. A larger $\sigma$ improves exploration and can accelerate surrogate optimization, while a smaller $\sigma$ better preserves the local structure of the original objective. This tradeoff motivates our decaying-$\sigma$ mechanism.
\subsubsection{Benefits of a Larger $\sigma$: Faster Surrogate Optimization}
Corollary 3.9 in \cite{GS-PowerOpt} implies, as formalized in our Theorem \ref{sigma-iter-complexity}, that when the smoothing radius $\sigma$ specified in (\ref{gs-poweropt-obj}) lies in a neighborhood of $0$, the iteration complexity of GS-PowerOpt scales inversely with $\sigma^2$. This suggests that choosing a larger $\sigma$ can reduce the number of required iterations. 

\begin{theorem}
\label{sigma-iter-complexity}
Under the conditions specified in \cite[Corollary 3.9]{GS-PowerOpt}, for some constant $U>0$ that depends on $e^{Nf(\bm{x}^*)}$ and the learning rate parameter $\gamma$, whenever $\sigma\in (0, U)$, the iteration complexity of GS-PowerOpt is $O(\sigma^{-2}d^2\varepsilon^{-1})^{2/(1-2\gamma)}$. 
\end{theorem}

\subsubsection{Benefits of a Smaller $\sigma$: Better Local Refinement}
Since GS-PowerOpt maximizes the surrogate objective $F_{N,\sigma}$ for maximizing the original objective $f$, aligning the maxima of $F_{N,\sigma}$ and $f$ are important. We prove under mild assumptions on $f$ that these maxima are better aligned with smaller $\sigma$. Specifically, near any sharp local maximizer $\bm{y}^*$ of $f$, for sufficiently small $\sigma$, $F_{N,\sigma}$ has a unique local maximizer $\bm{z}^*$ satisfying $\|\bm{z}^*-\bm{y}^*\|=\mathcal{O}(\sigma^{1/p})$ for some $p>0$.

\begin{proposition}[Local maximizer preservation under small $\sigma$; informal]
\label{unique-max}
Assume that $f$ is twice continuously differentiable, Lipschitz continuous, upper bounded, and has a bounded Hessian. Suppose $\boldsymbol{y}^*$ is a locally strictly concave local maximizer of $f$. Then, for any fixed $N>0$ and any sufficiently small neighborhood of $\boldsymbol{y}^*$, there exists a sufficiently small $\sigma>0$ such that $F_{N,\sigma}$ has a unique local maximizer $\boldsymbol{z}^*$ in that neighborhood.
\end{proposition}

The next proposition further quantifies how close this surrogate maximizer is to the original local maximizer when $\boldsymbol{y}^*$ is sharp.

\begin{proposition}[Displacement bound under local sharpness; informal]
\label{mu-dist}
Suppose $\boldsymbol{y}^*$ is a locally sharp maximizer of $f$, i.e., for some $\kappa,p>0$,
\[
f(\boldsymbol{y}^*)-f(\boldsymbol{x})
\geq
\kappa\|\boldsymbol{x}-\boldsymbol{y}^*\|^p
\]
in a neighborhood of $\boldsymbol{y}^*$. If $F_{N,\sigma}$ has a local maximizer $\boldsymbol{\mu}^*_{N,\sigma}$ in this neighborhood, then
\[
\|\boldsymbol{\mu}^*_{N,\sigma}-\boldsymbol{y}^*\|
=
\mathcal{O}(\sigma^{1/p}).
\]
Therefore, smaller $\sigma$ leads to better local alignment between the surrogate maximizer and the corresponding maximizer of the original objective.
\end{proposition}
\begin{remark}
The formal version of the two propositions can be found in our appendix.
\end{remark}
As illustrated in Figure \ref{sigma-effect}(a), for a fixed power $N$, decreasing the smoothing parameter $\sigma$ reduces the discrepancy between the global maximizer of the original objective and that of its Gaussian-smoothed counterpart. 

Together, Theorem~\ref{sigma-iter-complexity} and Propositions~\ref{unique-max}--\ref{mu-dist} reveal an inherent limitation of the fixed-$\sigma$ mechanism in GS-PowerOpt. A large fixed $\sigma$ is beneficial for exploration and finite-time surrogate optimization, but may reduce local accuracy. A small fixed $\sigma$ improves local refinement, but can weaken exploration and increase the iteration complexity.

\subsection{Another $\sigma$ Trade-off from Empirical Findings}
\begin{figure}[htb!]
	\centering
      \begin{tabular}{cc}
        \subfloat[]{\includegraphics[scale=0.3]{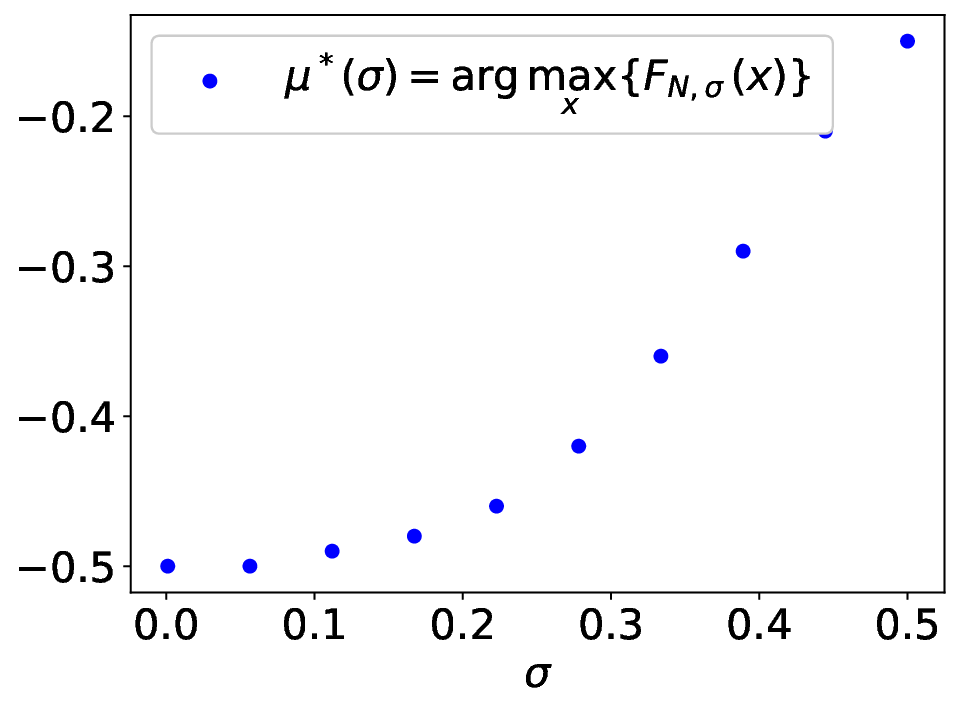} }
        &
        \subfloat[]{\includegraphics[scale=0.3]{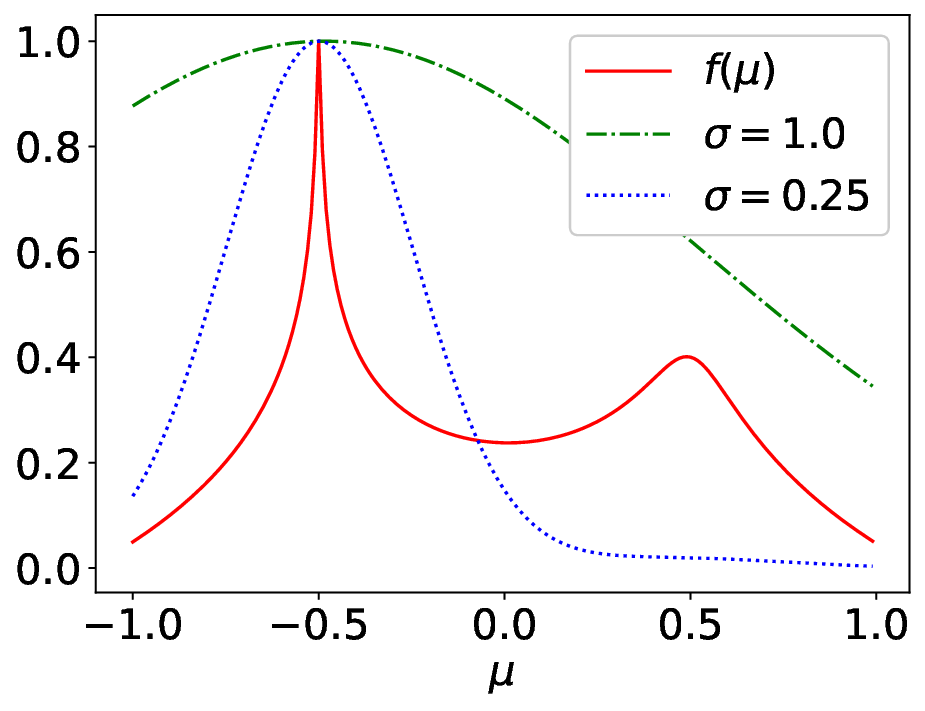}}
\end{tabular}
       	\caption{(a) $\mu^*(\sigma)$ denotes the estimated global $x$-maximizer of $F_{N,\sigma}(x)$, the Gaussian smoothed version of an example function $f:\mathbb{R}\rightarrow\mathbb{R}$ with a global maximizer $x^*=-0.5$ (graph in (b)), where $N=0.2$. It shows that  as $\sigma$ decreases from $0.5$ to $0.001$, the global maximizer $\mu^*(\sigma)$ approaches $x^*$. (b) The graphs of an example $f$ and its Gaussian smoothed function $F_{N,\sigma}$, with $N=1$ and $\sigma$ assuming different values. All three function graphs in (b) are scaled to have a maximum of 1 for easier comparisons.} 
\label{sigma-effect}
\end{figure}

The $\sigma$-decaying mechanism designed in our approach is also inspired from an experimental finding: when $N$ is large, a smaller $\sigma$ tends to amplify $\|\nabla F_{N,\sigma}(\bm{\mu})\|$ within the region $\mathcal{S}$ near $f$'s global maximizer(s), while rendering $\|\nabla F_{N,\sigma}(\bm{\mu})\|$ negligible in regions outside $\mathcal{S}$. This is illustrated with a simple example in Figure \ref{sigma-effect}(b). This provides another motivation for initially setting $\sigma$ to a relatively large value and then gradually decreasing it in the solution update iterations of $\bm{\mu}_{t+1}=\bm{\mu}_{t+1}+\alpha_t \hat{\nabla}F_{N,\sigma}(\bm{\mu}_t)$. 

Below, we provide an intuitive explanation on the finding, under the assumption that $f$ has a unique global maximum $\bm{x}^*$. However, as evidenced by the experimental results in Subsection \ref{Experiment-IAA}, our approach remains effective for objectives with potentially many global optimums.

The gradient admits the representation 
$$ \nabla F_{N,\sigma}(\bm{\mu})
= \frac{1}{\sigma^{2}}
\mathbb{E}_{\bm{X} \sim \mathcal{N}(\bm{\mu}, \sigma^2 I_d)}
\left[
(\bm{X}-\bm{\mu}) e^{N f(\bm{X})}
\right]. $$
Thus, it can be interpreted as a Gaussian-weighted average of the vector field $(\bm{x}-\bm{\mu}) e^{N f(\bm{x})}$ over $\mathbb{R}^d$.

Let $\bm{x}^*$ denote a global maximizer of $f$, and let $\mathcal{S}_{\bm{x}^*}$ be a sufficiently small neighborhood of $\bm{x}^*$. For fixed and sufficiently large $N$, the factor $e^{N f(\bm{x})}$ concentrates sharply around $\bm{x}^*$, so that the dominant contribution to the integral arises from $\bm{x} \in \mathcal{S}_{\bm{x}^*}$. The effective contribution of this neighborhood depends on the overlap between $\mathcal{S}_{\bm{x}^*}$ and the Gaussian density centered at $\bm{\mu}$ with variance $\sigma^2 I_d$.

When $\bm{\mu}$ lies near $\mathcal{S}_{\bm{x}^*}$, choosing a smaller $\sigma$ increases the concentration of the Gaussian kernel around $\bm{\mu}$ and hence around $\bm{x}^*$, thereby amplifying the contribution from $\mathcal{S}_{\bm{x}^*}$. In contrast, when $\bm{\mu}$ is far from $\bm{x}^*$, a larger $\sigma$ increases the Gaussian mass assigned to $\mathcal{S}_{\bm{x}^*}$, enhancing its contribution to the integral and resulting in a larger gradient magnitude in the direction of $\bm{x}^*$.

\section{Motivations for GS-PowerHP}
The exploration--refinement tradeoff discussed in the previous section motivates the design of GS-PowerHP. Specifically, we adopt a continuation strategy for the smoothing radius $\sigma$: starting with a relatively large value to encourage global exploration and maintain sufficient interaction with the basin of attraction of $\bm{x}^*$ when $\bm{\mu}$ is far away, and then gradually decreasing $\sigma$ as the iterates ${\bm{\mu}_t}$ approach high-value regions. This adaptive schedule allows the algorithm to combine informative long-range search in the early stage with more accurate local refinement near the maximizer, thereby mitigating the limitation of using a fixed smoothing radius throughout optimization. It also helps the gradient estimator $\hat{\nabla} F_{N,\sigma}(\bm{\mu}_t)$ remain sufficiently large even when $\bm{\mu}_t$ is far from $f$'s global maximizer(s).

\section{The Proposed Method: GS-PowerHP}
GS-PowerHP aims to solve the deterministic optimization problem of
$\max_{\bm{w}\in\mathbb{R}^d} f(\bm{w})$. With the pre-selected $N$, at each iteration $t$, GS-PowerHP performs a one-step stochastic gradient ascent to solve 
\begin{equation}
\label{FN}
\max_{\bm{\mu}}F_{N,\sigma_t}(\bm{\mu}):=\mathbb{E}_{\bm{x}\sim\mathcal{N}(\bm{\mu},\sigma^2_t I_d)}[f_N(\bm{x})].
\end{equation}
Specifically, the rule for updating the solution candidate is, for all $t\geq 1$,
\begin{equation}
\begin{split}
\label{zo-slph}
\text{GS-PowerHP}:\qquad & \sigma_{t+1} = \beta^{t+1}\sigma_0 + b,\\
&\bm{\mu}_{t+1} = \bm{\mu}_{t} + \alpha_t \hat{\nabla} F_{N,\sigma_{t+1}}(\bm{\mu}_t),\\
\end{split}
\end{equation}
where $\beta\in (0,1)$ is the discounter factor, $b>0$ is a lower bound for $\sigma_t$, $\hat{\nabla} F_{N,\sigma_{t+1}}(\bm{\mu}_t):=\frac{1}{K} \sum_{k=1}^K (\bm{w}_k-\bm{\mu}_t)f_N(\bm{w}_k)$, $\{\bm{w}_k\}_{k=1}^K$ are independently sampled from the multivariate Gaussian distribution $\mathcal{N}(\bm{\mu}_t,\sigma_{t+1}^2 I_d)$, and $f_N(\bm{w})$ is defined as
\begin{equation}
\label{fN}
f_N(\bm{w})=e^{Nf(\bm{w})}.
\end{equation}
The values of $\beta$, $\sigma_0$, $b$, $\bm{\mu}_0\in\mathbb{R}^d$, and $K$ are hyper-parameters. Note that $\hat{\nabla} F_{N,\sigma}(\bm{\mu}_t)$ is an unbiased sample estimate of $\sigma^{2}\nabla F_{N,\sigma}(\bm{\mu}_t)$, since 
\begin{equation}
\begin{split}
\label{nablaF-est}
&\sigma^2\nabla F_{N,\sigma}(\bm{\mu}_t) =\sigma^2\nabla_{\bm{\mu}}  \mathbb{E}_{\bm{w}\sim \mathcal{N}(\bm{\mu}_t,\sigma^2 I_d)}[f_N(\bm{w})] \\
&=  \frac{1}{(\sqrt{2\pi}\sigma)^{d}} \int_{\bm{x}\in \mathbb{R}^d} (\bm{w}-\bm{\mu}_t)f_N(\bm{w}) e^{-\frac{\lVert \bm{w} - \bm{\mu}_t \rVert^2}{2\sigma^2}} d\bm{x}\\
&= \mathbb{E}_{\bm{w}\sim \mathcal{N}(\bm{\mu}_t,\sigma^2 I_d)}[(\bm{w}-\bm{\mu}_t)f_N(\bm{w})] = \mathbb{E}[\hat{\nabla}F_{N,\sigma}(\bm{\mu}_t)],
\end{split}
\end{equation}
where the interchange of differentiation and integral in the second line can be justified by Lebesgue's dominated convergence theorem.

\begin{algorithm}[h]
\caption{GS-PowerHP for Solving $\max_{\bm{x} \in \mathbb{R}^d} f(\bm{x})$.}
\label{alg:gspowerhp}
\begin{algorithmic}[1]
\REQUIRE Power $N > 0$, initial smoothing radius $\sigma_0 > 0$, the minimum smoothing radius $b>0$, decay rate $\beta \in (0,1)$, objective $f$, initial point $\bm{\mu}_0 \in \mathbb{R}^d$, samples per iteration $K$, total iterations $T$, learning rates $\{\alpha_t\}_{t=0}^{T-1}$.
\FOR{$t = 0$ to $T-1$}
    \STATE Set $\sigma_{t+1} = \sigma_0 \beta^{t+1}+b$
    \STATE Sample i.i.d. $\{\bm{x}_k\}_{k=1}^K$ from $\mathcal{N}(\bm{\mu}_t, \sigma_{t+1}^2 I_d)$
    \STATE Compute the gradient estimator:
    \[
     \hat{\nabla} F_{N,\sigma_{t+1}}(\bm{\mu}_t):=\frac{1}{K} \sum_{k=1}^K (\bm{x}_k-\bm{\mu}_t)e^{Nf(\bm{x}_k)}
    \]
    \STATE Update:
    \[
    \bm{\mu}_{t+1} = \bm{\mu}_t + \alpha_t \hat{\nabla} F_{N,\sigma_{t+1}}(\bm{\mu}_t)
    \]
\ENDFOR
\STATE \textbf{Return} $\bm{\mu}^* = \arg\max_{\bm{\mu} \in \{\bm{\mu}_1, \dots, \bm{\mu}_T\}} f(\bm{\mu})$.
\end{algorithmic}
\end{algorithm}

\begin{remark}[Techniques to avoid computational overflow]
To prevent computational overflow potentially caused by large exponents, we can replace $e^{N(f(\bm{x}_k))}$ with $e^{N(f(\bm{x}_k)-f(\bm{\mu}_t))}$ for size control. It does not affect the value of the normalized gradient. Also, we can choose to output the optimal sample $\bm{x}$ found in the whole process instead of $\bm{\mu}^*$.

Another practical approach is to shift the objective by a constant so that it becomes negative. If an upper bound $M$ of $f$ is known, one may replace $f$ by $\tilde f(\bm{x})=f(\bm{x})-M$,
so that $\tilde f(\bm{x})\leq 0$ on the domain. This shift does not change the maximizers of the original objective, and it only rescales the corresponding power-smoothed objective by the constant factor $e^{-NM}$. Such bounded or normalized objectives are common in machine-learning applications; see Remark~\ref{rmk:shift-f} for a representative example.
\end{remark}

\section{Convergence Analysis}
\label{sec-convergence-analysis}
In this section, we establish in Corollary \ref{summary-cor} that, under mild assumptions (Assumptions~\ref{coercivity}, \ref{f-Lipschitz}, \ref{lr}, and \ref{mu-bounded}), GS-PowerHP converges in expectation to an arbitrarily small neighborhood of the global maximizer $\bm{x}^* = \arg\max_{\bm{x} \in \mathbb{R}^d} f(\bm{x})$, provided that the power $N > 0$ is sufficiently large. With an appropriately chosen learning rate, the iteration complexity approaches $O_N((d^2\varepsilon^{{-1}}b^{-2})^{\frac{2}{1-2\gamma}})$.

The proof consists of two parts. The first one shows in Corollary \ref{iter-complexity} that the solutions $\{\bm{\mu}_t\}$ produced by GS-PowerHP's iteration process (\ref{zo-slph}) converges in expectation to a stationary point of $F_{N,\sigma}$. The second part shows in Section \ref{global-sec} that, for any given positive $\delta$ and $b$, there exists a threshold $N_{\delta,b}$ such that for any $N>N_{\delta,b}$, all the stationary points of $F_{N,\sigma}$ lies in the $\bm{x}^*$-neighborhood of $\mathcal{S}_{\bm{x}^*,\delta}:=\{\bm{\mu}\in\mathbb{R}^d : |\mu_i-x_i^*|\leq \delta, i\in\{1,2,...,d\}\}$, as long as $\sigma>b$. Here, $\mu_i$ and $x_i^*$ denote the $i^{th}$ entry of $\bm{\mu}$ and $\bm{x}^*$, respectively. Finally, our stated convergence result in Corollary \ref{summary-cor} is implied by the two parts collectively. The full proofs for all our theoretical results can be found in the appendix (Appendix \ref{appendix-proof}). 


\subsection{Assumptions}
We make three assumptions that are standard in optimization researches for machine learning.
\begin{assumption}
\label{coercivity}
Assume that the maximization objective $f:\mathbb{R}^d\rightarrow \mathbb{R}$ is continuous, $\lim_{\|\bm{x}\|\rightarrow +\infty} f(\bm{x})=-\infty$, and has a unique global maximum point $\bm{x}^*=\arg\max_{\bm{x}\in\mathbb{R}^d}f(\bm{x})$. 
\end{assumption}
\begin{remark}
The counterpart, $\lim_{\|\bm{x}\|\rightarrow +\infty} f(\bm{x})=+\infty$ for $\min_{\bm{x}} f(\bm{x})$ is the commonly seen coercivity assumption.
\end{remark}

\begin{assumption}
\label{f-Lipschitz}
Assume that $f$ is Lipschitz in $\mathbb{R}^d$ with its Lipschitz constant equal to $L_f$.
\end{assumption}

\begin{assumption}
\label{lr}
$\alpha_t>0$, $\sum_{t=0}^{\infty}\alpha^2_t<\infty$, and $\sum_{t=0}^{\infty}\alpha_t=\infty$.
\end{assumption}

\subsection{Convergence of $\mu_t$ to a stationary point of $F_{N,\sigma}$}
The differentiability and the Lipschitz constant of the Gaussian-smoothed objective $F_{N,\sigma}$ is given in the following lemma.

\begin{lemma}
\label{F-property}
Under Assumption \ref{coercivity}, given any $N>0$ and $\sigma>0$, (1) both $F_{N,\sigma}(\bm{\mu})$ and $\nabla F_{N,\sigma}(\bm{\mu})$ are well-defined and Lipschitz in $\mathbb{R}^d$; (2) The Lipschitz constant for $\nabla F_{N,\sigma}$ is $L=2d\sigma^{-2}e^{Nf(\bm{x}^*)}$, (3) $F_{N,\sigma}$ has at least one global maximum $\bm{\mu}^*\in\mathbb{R}^d$, and (4) $\mathbb{E}[\|\hat{\nabla}F_{N,\sigma}(\bm{\mu})\|^2]\leq G=d\sigma^2 e^{2Nf(\bm{x}^*)}$ for all $\bm{\mu}\in\mathbb{R}^d$.
\end{lemma}

With the Lipschitz constant given in Lemma \ref{F-property}, we bound the expected sum of gradients, $\sum_{t=0}^{T-1}\alpha_t\sigma_{t+1}^2\mathbb{E}[\|\nabla F_{N,\sigma_{t+1}}(\bm{\mu}_t) \|^2]$. Since $\lim_{T\rightarrow\infty}\sum_{t=0}^{T-1}\alpha_t=\infty$, by Assumption \ref{lr}, the boundedness implies $\lim_{T\rightarrow \infty}\mathbb{E}[\|\nabla F_{N,\sigma_{T}}(\bm{\mu}_{T-1}) \|^2]=0$.

\begin{theorem}
\label{main-theorem}
For any deterministic $\bm{\mu}_0\in\mathbb{R}^d$, $N,b,\sigma_0>0$, and positive integer $T$, let $\{(\bm{\mu}_t,\sigma_{t})\}_{t=1}^{T}$ be generated by the GS-PowerHP iterations (\ref{zo-slph}). Under Assumption \ref{coercivity} and \ref{f-Lipschitz}, we have
\begin{align*} \sum_{t=0}^{T-1}\alpha_t&\sigma_{t+1}^2\mathbb{E}[\|\nabla F_{N,\sigma_{t+1}}(\bm{\mu}_t) \|^2]\leq e^{Nf(\bm{x}^*)} -F_{N,\sigma_0}(\bm{\mu}_0)\\
&+ 2d^2e^{3N f(\bm{x}^*)} \sum_{t=0}^{T-1}\alpha_t^2 + \sigma_0C_{N,d,f}(1-\beta)\sum_{t=0}^{T-1}\beta^t,  
\end{align*}
where $C_{N,d,f}=N e^{Nf(\bm{x}^*)} L_f \sqrt{d}$.
\end{theorem}

With the bound in Theorem \ref{main-theorem}, we derive the convergence rate and iteration complexity for $\bm{\mu}_t$ converging to stationary points of $F_{N,\sigma}$.

\begin{corollary}[Convergence Rate]
\label{convergence-rate} Assume the conditions in Theorem \ref{main-theorem} and Assumption \ref{lr}. For any positive integer $T$, define 
$$\nu_T:= \min_{\tau\in\{0,1,...,T\}}\mathbb{E}[\|\nabla F_{N,\sigma_{\tau+1}}(\bm{\mu}_\tau)\|^2],$$ where $\{\bm{\mu}_{\tau}\}$ are generated from GS-PowerHP (\ref{zo-slph}) with $\beta\in [0,1)$. Then, we have $\nu_T=\mathcal{O}_N\left( (\sigma_0\beta^T+b)^{-2}(\sum_{t=0}^{T-1} \alpha_t )^{-1}\right)$. Furthermore, the dependence of the $\mathcal{O}_N$-factor on $N$ can be removed under the additional assumption of $f(\bm{x})< 0$ for all $\bm{x}\in\mathbb{R}^d$.
\end{corollary}
\begin{remark}
This corollary clarifies this advantage by providing a first-order stationarity convergence rate of GS-PowerHP. This rate improves when the current smoothing scale $\sigma_0\beta^T+b$ is larger, explaining why using the small refinement scale too early can slow finite-time stationarity. Thus, compared with a GS-PowerOpt that uses the fixed small refinement scale $b$ from the beginning (equivalent to GS-PowerHP with $\beta=0$), GS-PowerHP ($\beta>0$) benefits from a \textbf{faster} stationarity convergence rate during the finite-time optimization process while still approaching the same refinement scale $b$ as $T\to\infty$, which supports late-stage \textbf{accuracy}. The experimental results in Table \ref{Tab:effects-b-beta} illustrates the expected tradeoff: smaller late-stage smoothing, controlled by $\beta$ and $b$, often improves final fitness but increases the iterations needed to reach the gradient threshold. The extreme $\beta=0$ case uses the small scale from the beginning and fails due to poor exploration.
\end{remark}

\begin{corollary}[Iteration Complexity]
\label{iter-complexity}
Assume the conditions in Theorem \ref{main-theorem} and Assumption \ref{lr}. Let $\gamma\in (0,1/2)$ and $\alpha_t = (t+1)^{-(1/2+\gamma)}$. For any $\varepsilon\in(0,1)$, after $T> (C_2C_1d^2\varepsilon^{-1})^{\frac{2}{1-2\gamma}}$ $=\mathcal{O}_N((d^2\varepsilon^{{-1}}b^{-2})^{\frac{2}{1-2\gamma}})$ times of parameter updating, we have
$$\min_{t\in\{0,1,2,...,T\}}\mathbb{E}[ \| \nabla F_{N,\sigma_{t+1}}(\bm{\mu}_t) \|^2 ]<\varepsilon,$$
where $C_0:=f_N(\bm{x}^*) - F_{N,\sigma_0}(\bm{\mu}_{0}) + 2f_N^3(\bm{x}^*) \sum_{t=1}^\infty t^{-(1+2\gamma)} + \sqrt{2}\sigma_0Nf_N(\bm{x}^*)L_f$, $C_1 = b^{-2}C_0(1-2\gamma)$, and $C_2=\max\{1,1/C_1\}$. Furthermore, the dependence of the $\mathcal{O}_N$-factor on $N$ can be removed under the additional assumption of $f(\bm{x})< 0$ for all $\bm{x}\in\mathbb{R}^d$.
\end{corollary}
\begin{remark}
\label{rmk:shift-f}
The additional assumption that $f(\bm{x})<0$ is satisfied by a broad class of optimization problems in machine learning, particularly when a loss-minimization problem is reformulated as an equivalent maximization problem. For example, in supervised regression, including neural-network regression, the training loss $L(\bm{x})$ is often defined as the sum of squared prediction errors and is therefore nonnegative. Maximizing the negative loss is equivalent to minimizing the original loss. Thus, one may define $f(\bm{x}):= -L(\bm{x}) - 0.1$,
where the constant shift $-0.1$ is introduced solely to ensure that $f(\bm{x})<0$ for all $\bm{x}$. 
\end{remark}

\subsection{Convergence to the Global Maximum of $f$}
\label{global-sec}
In this section, under an additional mild assumption (i.e., Assumption \ref{mu-bounded}), we show that for any given $\delta>0$ and $\sigma>0$, there exists a threshold $N_{\delta,\sigma}>0$ such that, whenever $N>N_{\delta,\sigma}$ all the stationary points of $F_{N_\delta,\sigma}$ lies in the $\delta$-neighborhood $\mathcal{S}_{\bm{x}^*,\delta}$, as long as $\sigma>b$.

\begin{theorem}
\label{thm2}
Suppose Assumption \ref{coercivity} holds. Given any positive numbers $N$, $\sigma$, $\delta$, and $M$ such that $\delta<M$, there exists $N_{\delta,\sigma,M}>0$ such that whenever $N>N_{\delta,\sigma,M}$ the following statement holds true: if $\bm{\mu}\in\mathbb{R}^d$ and $\|\bm{\mu}\|<M$, then for any $i\in\{1,2,...,d\}$
\begin{enumerate}
 \item $\mu_i>x^*+\delta \Rightarrow \frac{\partial F_{N,\sigma}(\bm{\mu})}{\partial \mu_i}<0$; and
 \item $\mu_i<x^*-\delta \Rightarrow \frac{\partial F_{N,\sigma}(\bm{\mu})}{\partial \mu_i}>0$,
\end{enumerate}
where $\mu_i$ and $x^*_i$ denotes the $i^{th}$ entry of $\bm{\mu}$ and $\bm{x}^*$, respectively.
\end{theorem}

\begin{assumption}
\label{mu-bounded}
Given any $\delta,b>0$ and any $\bm{\mu}_0\in\mathbb{R}^d$, assume that there is some $N_{\delta,b,\bm{\mu}_0}>0$ such that for all $N>N_{\delta,b,\bm{\mu}_0}$, $\{\bm{\mu}_t\}$ generated by the iteration process (\ref{zo-slph}) are bounded by some constant $M_{(\delta,\sigma,\bm{\mu}_0)}$. 
\end{assumption}

Assumption~\ref{mu-bounded} is consistent with the qualitative behavior suggested by Theorem~\ref{thm2} when $N$ is sufficiently large. In particular, Theorem~\ref{thm2} shows that outside the $\delta$-neighborhood $\mathcal{S}_{\bm{x}^*,\delta}$, each coordinate of $\bm{\mu}$ experiences a directional drift toward this neighborhood. This indicates that, in regions far from $\bm{x}^*$, the gradient field points inward in a componentwise sense.

While this observation does not by itself constitute a formal boundedness proof for the stochastic iteration, it suggests that the deterministic drift of the dynamics discourages escape to infinity. When combined with a decreasing sequence $\{\sigma_t\}$ that preserves this directional property for large $N$, it is therefore reasonable to expect that the iterates remain in a bounded region.

Such bounded-iterate assumption is also a technical condition used in stochastic approximation and non-convex SGD analyses, as noted in \cite{Mertikopoulos2020}. It persists in the literature (e.g., \cite{Borkar2022,bonnabel2013stochastic}).

\begin{corollary}
\label{summary-cor}
Suppose Assumption \ref{coercivity}, \ref{f-Lipschitz}, \ref{lr}, and \ref{mu-bounded} hold. For any $\delta>0$, $b>0$, and $\bm{\mu}_0\in\mathbb{R}^d$, let $N_{\delta,b,\bm{\mu}_0}$ and $M_{(\delta,b,\bm{\mu}_0)}$ be the threshold and $\|\bm{\mu}\|$-bound in Assumption \ref{mu-bounded}, respectively. Let $N_{\delta,b,M_{(\delta,b,\bm{\mu}_0)}}$ be the threshold stated in Theorem \ref{thm2}. Then, whenever $N>\max\{N_{\delta,b,\bm{\mu}_0},N_{\delta,b,M_{(\delta,b,\bm{\mu}_0)}}\}$, $\{\bm{\mu}_t\}$ generated by (\ref{zo-slph}) converges in expectation to $\mathcal{S}_{\bm{x}^*,\delta}$ with the iteration complexity of $O((d^2\varepsilon^{-1}b^{-2})^{2/(1-2\gamma)})$, given the learning rate of $\alpha_t = (t+1)^{-(1/2+\gamma)}$ for some $\gamma\in (0,1/2)$.
\end{corollary}

\section{Experiments}
\label{sec-experiments}
In this section, we test the capacity of GS-PowerHP with extensive experiments. We use a normalized implementation variant of Algorithm~\ref{alg:gspowerhp}: 
$
\bm{\mu}_{t+1} = \bm{\mu}_t + \alpha_t \frac{\hat{\nabla} F_{N,\sigma_{t+1}}(\bm{\mu}_t)}{\|\hat{\nabla} F_{N,\sigma_{t+1}}(\bm{\mu}_t)\|}.
$
This normalization stabilizes the step size when the magnitude of the zeroth-order estimator varies substantially across iterations. The normalized update preserves the stochastic search direction but changes the step length. Therefore, the convergence guarantees in Section~\ref{sec-convergence-analysis} should be interpreted as applying to the unnormalized theoretical update, while the normalized update is used as an implementation variant.

The algorithms evaluated in our experiments can be categorized into two main groups: smoothing-based zeroth-order optimization methods and the state-of-the-art evolutionary algorithm, CMA-ES \cite{hansen2016cma}. The former category comprises GS-PowerHP, ZOSGD \cite{ghadimi2013}, ZO-AdaMM \cite{Chen2019zo}, GS-PowerOpt \cite{GS-PowerOpt},\footnote{GS-PowerOpt can be implemented in two variants: PGS, which transforms the objective function $f$ to $f^N$, and EPGS, which transforms $f$ to $e^{Nf}$. In this work, GS-PowerOpt specifically refers to the EPGS variant.} the standard homotopy method (STD-Homotopy), and ZOSLGHd and ZOSLGHd \cite{Iwakiri2022}. A brief overview of some of these zeroth-order methods is provided in Appendix D of \cite{GS-PowerOpt}.

Square Attack~\cite{Andriushchenko2020}, a powerful black-box adversarial attack designed specifically for image classification, is also included but reported separately in Appendix~\ref{appendix:square-attack}. We include it as a specialized attack reference rather than as an objective-matched optimization baseline, since it uses a different attack objective and stopping rule from the general-purpose zeroth-order optimization methods considered in the main experiments.

The hyper-parameters' values for each algorithm are determined through a grid search, with the search space and selection procedure detailed in Appendix~\ref{appendix-hyperparameters}.

\subsection{Performances on Optimizing Benchmark Objectives}
\label{benchmark}
\begin{figure}[htb!]
	\centering
      \begin{tabular}{cc}
        \subfloat[]{\includegraphics[scale=0.25]{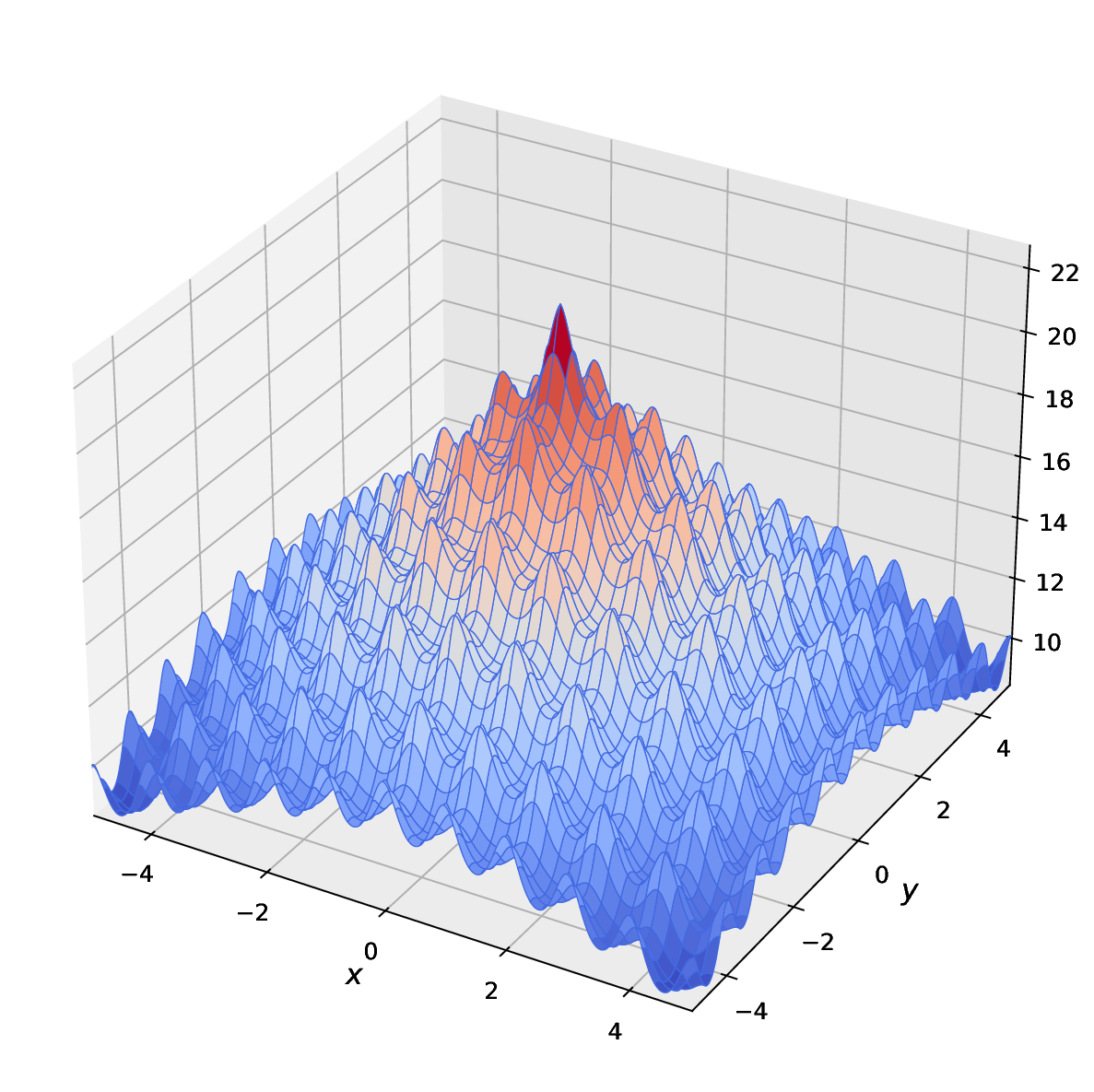} }
        &
        \subfloat[]{\includegraphics[scale=0.25]{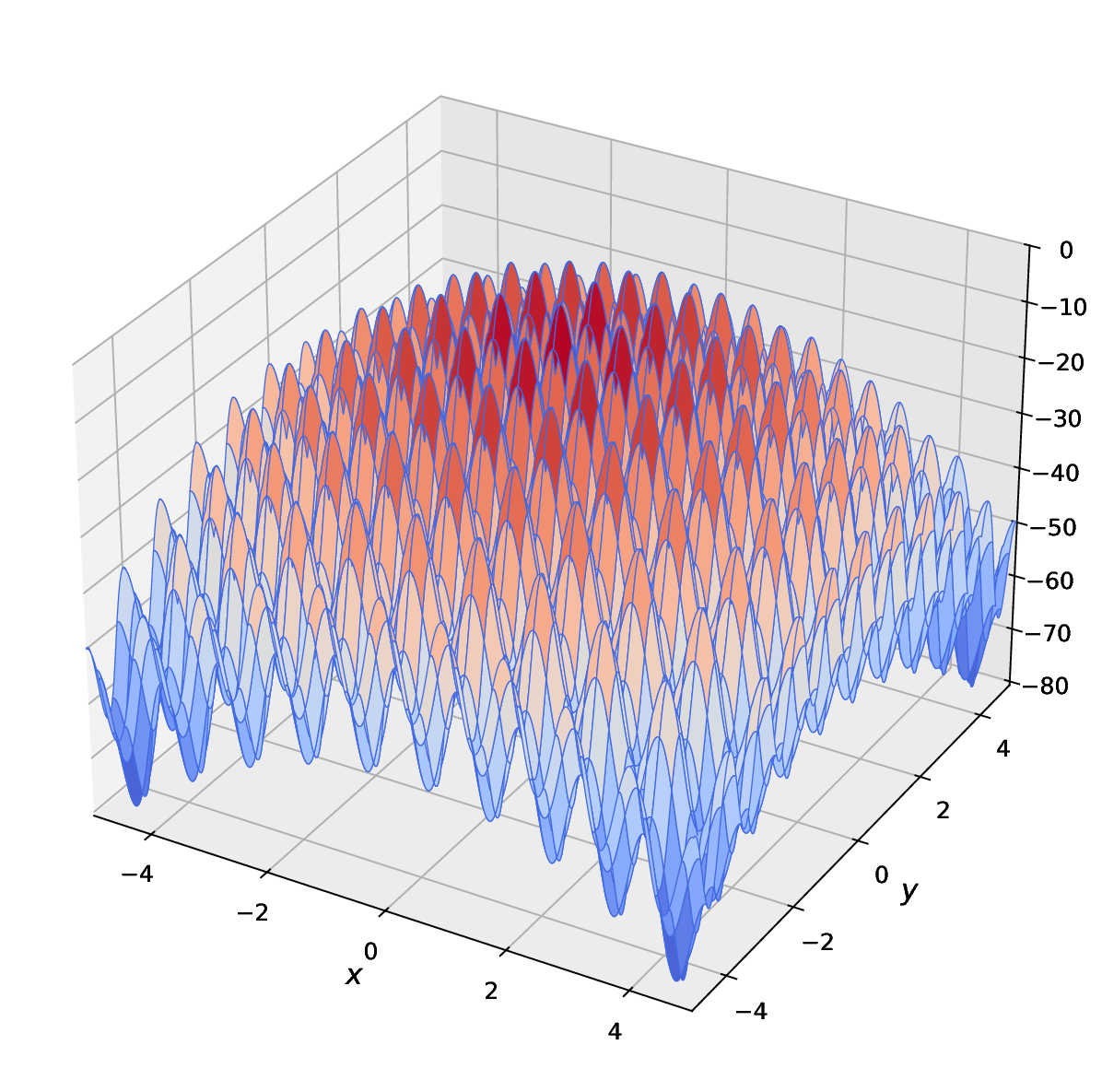}}
\end{tabular}
       	\caption{(a) Maximize-Version of Ackley ($d=2$, figure borrowed from \cite[Figure 3]{GS-PowerOpt}). (b) Maximize-Version of Rastrigin ($d=2$) .} 
\label{fig-Ackley-Rastrigin}
\end{figure}

We test GS-PowerHP on two popular benchmark objectives for non-convex optimization methods: the Ackley and Rastrigin functions. We use their maximization versions. For $\bm{x}=(x_1,\ldots,x_d)\in\mathbb{R}^d$, their functional forms are
\begin{equation}
\begin{split}
\text{Ackley:}\quad
&f(\bm{x})
=
20\exp\left(
-0.2\sqrt{\frac{1}{d}\sum_{i=1}^d x_i^2}
\right)
+
\exp\left(
\frac{1}{d}\sum_{i=1}^d \cos(2\pi x_i)
\right),\\
\text{Rastrigin:}\quad
& f(\bm{x})
=
-\left[
10d+\sum_{i=1}^d\left(x_i^2-10\cos(2\pi x_i)\right)
\right].
\end{split}
\end{equation}
These two objectives are non-trivial to maximize because they are highly non-convex and multimodal. In particular, the Ackley function contains a broad flat region away from the global maximizer, while the Rastrigin function has many regularly distributed local maxima and minima. Their 2D-version landscapes are shown in Figure~\ref{fig-Ackley-Rastrigin}. In our experiments, we set $d=100$, which represents a relatively high-dimensional benchmark setting for Ackley and Rastrigin-type test problems. 

We follow similar procedures of the experiments in \cite[Section 5.2 \& 5.3]{GS-PowerOpt} to test and compare GS-PowerHP's performance with other algorithms. Specifically, for each algorithm, we perform 100 trials to optimize the objective function. In each trial, each algorithm is performed for $T$ iterations. Let $\bm{\mu}^*_n$ denote the candidate solution in the $n^{th}$ trial with the largest $f$ value and let $t^*_n$ denote the number of solution updates taken to achieve $\bm{\mu}^*_n$. Then, for each algorithm, in Table \ref{tab:combined-ackley-rastrigin} we report $\bar{t}^*:=\sum_{n=1}^{100}t_n^*/100$, $\bar{f}(\bm{\mu}^*):=\sum_{n=1}^{100}f(\bm{\mu}_n^*)/100$, and the average mean square error $\sum_{n=1}^{100}\| \bm{\mu}^*_n-\bm{x}^*\|^2/(100d)$.
\begin{table}[h]
\caption{Comparison of optimization performance on Ackley and Rastrigin maximization problems. Both functions attain their global maximum at $\bm{x}^*=\bm{0}$. $\bar{t}^*$, $\bar{\bm{\mu}}^*$, and $\bar{f}(\bm{\mu}^*)$  are defined in Section \ref{benchmark}. All values are rounded to at most 2 decimal places. Total generations $T$ in each trial: $1,000$ for Ackley and $2,000$ for Rastrigin. The popolation size in each population is $500$. The reported statistics are averages, with standard deviation in the parentheses. The top performance (i.e., highest $\bar{f}$ value) among the smoothing-based algorithms is bolded (note that CMA-ES is not smoothing-based).}
\label{tab:combined-ackley-rastrigin}
\centering
\small
\begin{tabular}{l p{0.8cm}  p{2.3cm}  p{2cm} |  p{0.8cm}  p{2.3cm}  p{2cm}}
\toprule
& \multicolumn{3}{c}{Ackley ($d=100$)} 
& \multicolumn{3}{c}{Rastrigin  ($d=100$)} \\
\cmidrule(lr){2-4} \cmidrule(lr){5-7}
Algorithm & $\bar{t}^*$ & $\|\bm{\mu}^*-\bm{x}^*\|^2/d$ & $\bar{f}(\bm{\mu}^*)$ 
          & $\bar{t}^*$ & $\|\bm{\mu}^*-\bm{x}^*\|^2/d$ & $\bar{f}(\bm{\mu}^*)$ \\
\midrule
\textit{Our Algo.} & $930$ & $0.0\,(0.0)$ & $\mathbf{22.55}\,(0.01)$ 
                & $608$ & $0.0\,(0.0)$ & $\mathbf{-4.20}\,(0.39)$ \\
GS-PowerOpt     & $999$  & $0.01\,(0.0)$  & $22.03\,(0.05)$ 
                & $1815$  & $0.0\,(0.0)$   & $-78.14\,(10.0)$ \\
ZOSLGHd         & $785$  & $0.0\,(0.0)$  & $22.47$ 
                & $905$  & $0.99\,(0.0)$  & $-99.51\,(0.0)$ \\
ZOSLGHr         & $966$  & $0.0\,(0.0)$  & $22.40\,(0.02)$ 
                & $1992$  & $0.83\,(0.08)$   & $-83.79\,(7.77)$ \\
ZOAdaMM         & $915$  & $0.02\,(0.0)$  & $21.53\,(0.06)$ 
                & $1790$  & $0.99\,(0.0)$   & $-101.75\,(0.28)$ \\
ZOSGD           & $939$  & $21.06\,(0.56)$  & $10.21\,(0.10)$ 
                & $1784$   & $0.99\,(0.0)$   & $-101.85\,(0.28)$ \\
STD-Htp         & $990$  & $4.92\,(1.09)$   & $13.93$ 
                & $1706$  & $0.99\,(0.0)$   & $-103.50\,(0.68)$ \\
CMA-ES          & $980$  & $0.0\, (0.0)$  & $22.72\,(0.0)$ 
                & $803$   & $0.06$ (0.02)   & $-6.52\, (2.42)$ \\
\bottomrule
\end{tabular}
\end{table}

\subsection{Least-likely Targeted Black-Box Image Adversarial Attacks (IAA)}
\label{Experiment-IAA}
Given an image classifier $\mathcal{G}(\cdot)$ that outputs a probability distribution over classes, let $\mathcal{T}$ denote the least-likely class predicted for a clean image $\bm{a}$ (i.e., $\mathcal{T} = \arg\min \mathcal{G}(\bm{a})$). The goal of a least-likely targeted black-box IAA is to find a perturbation $\bm{x}$ with the \textit{minimal} norm such that the adversarial image $\bm{a} + \bm{x}$ is successful. Here, \textit{success} implies that $\mathcal{G}(\bm{a} + \bm{x})$ assigns the highest probability to the target class $\mathcal{T}$.

This problem can be formulated as an optimization problem. Suppose the pixels in the image $\bm{a}$ are normalized to $[-1,1]$ and the perturbation is transformed to $\bm{y}=\tanh(\bm{x})$ to restrict size. We define the loss in the following way, which is similar to the one used in \cite{CarliniWagner2017}.
\begin{equation}
 L(\bm{x}):= \max( \max_{i\neq\mathcal{T}} \mathcal{G}(\bm{a+y})_i -  \mathcal{G}(\bm{a+y})_{\mathcal{T}}, \; \kappa) + \lambda \|\bm{y}\|,
\end{equation}
where $\mathcal{G}(\bm{a+y})_i$ denotes the probability assigned to class $i$, and $\kappa<0$ is a hyper-parameter used to prevent excessive effort on attack confidence, which we set to $-0.001$. An attack is called \textit{successful} if it produces a successful perturbation, while a perturbation $\bm{x}$ is called successful if $\max_{i\neq\mathcal{T}} \mathcal{G}(\bm{a+y})_i - \mathcal{G}(\bm{a+y})_{\mathcal{T}}<\kappa$.

To test its capacity for large-dimension problems, we apply GS-PowerHP and other compared algorithms to solve $\max_{\bm{x}\in\mathbb{R}^d} f(\bm{x}):=-L(\bm{x})$ 
on images randomly selected from three popular image sets, the MNIST hand-written figures \cite{mnist1998}, CIFAR-10 \cite{krizhevsky2009}, and ImageNet \cite{imagenet-JiaDeng2009}. The images from MNIST, CIFAR-10, and ImageNet have $28\times28=784$, $32\times32\times3=3,072$, and $224\times224\times3=150,528$ pixels, respectively. These values correspond to the dimensionality $d$ of the input space for each dataset. We normalize the pixel range to $[-1,1]$ for all involved images before performing experiments.
 
For each algorithm and each dataset, we randomly pick $100$ images to attack. If the attack on the $m^{th}$ image $\bm{a_m}$ is successful, let $\bm{\mu}^*_m$ denote the one with the least $L_2$ norm among all the successful perturbations produced during the attack process, let $T_m$ denote the number of solution updates (i.e., number of iterations) to reach $\bm{\mu}^*_m$, and let $R_m^2:=1-\frac{\sum_{j=1}^d(\mu_j^*-a_j)}{(a_j-\bar{a})^2}$, where $\mu_j^*$ and $a_j^*$ denote the $j^{th}$ entry of $\bm{\mu}_m^*$ and $\bm{a_m}$, respectively.

For each algorithm, let $\mathcal{S}$ denote the set of indices corresponding to successful attacks out of the $100$ performed, and let $|\mathcal{S}|$ denote the number of successful attacks. The performance statistics we will report are the sucess rate SR$:= |\mathcal{S}|/100$, the average $R^2$ score $\overline{R^2}:=\sum_{m\in \mathcal{S}} R_m^2/|\mathcal{S}|$, the average perturbation $L_2$ norm $\overline{\|\bm{\mu}^*\|_2}:=\sum_{m\in \mathcal{S}} \|\bm{\mu}_m^*\|/|\mathcal{S}|$, and the average number of iterations taken to reach the optimal $\overline{T}:=\sum_{m\in \mathcal{S}} T_m/|\mathcal{S}|$.

For each smoothing-based method, we use $K=10$ random samples for every solution update. The hyper-parameters values are selected by grid search and are reported in Appendix \ref{appendix-hyperparameters}.

The classifiers for MNIST and CIFAR-10 are neural networks trained using the distillation technique (\cite{CarliniWagner2017}). The test scores is $98\%$ for the MNIST classifer and $86\%$ for the CIFAR-10 classifier. For ImageNet images, we use a pretrained ResNet-50 (\cite{ResNet, tfresnet50}) as the classifier, which is a common choice in this research field (\cite{Andriushchenko2020}).

\begin{table}[t]
\caption{Per-image adversarial attacks on 100 images from MNIST and CIFAR-10. 
MNIST: $\bm{\mu}_0 = \bm{0}$, $T_{total} = 2{,}500$; GS-PowerHP hyperparameters: $N=0.5$, $\sigma_0=0.05$, $\alpha=0.065$, $\beta=0.999$. 
CIFAR-10: $\bm{\mu}_0 = \bm{0}$, $T_{total} = 2{,}500$; GS-PowerHP hyperparameters: $N=0.15$, $\sigma_0=0.03$, $\alpha=0.06$, $\beta=0.999$. 
All values are mean (standard deviation) except SR. SR = success rate; $\bar{R}^2$ measures similarity between the original and perturbed images; $\overline{\|\bm{\mu}^*\|_2}$ = average $\ell_2$ perturbation norm; $\bar{T}$ = average number of iterations. The top performance (i.e., highest SR and lowest $\overline{\|\bm{\mu}^*\|_2}$) among the smoothing-based algorithms is bolded (note that CMA-ES is not smoothing-based).}
\label{tab:per-image-attacks-mnist-cifar}
\centering
\small
\begin{tabular}{l 
  c c c c 
  | 
  c c c c}
\toprule
& \multicolumn{4}{c}{MNIST (100 images)} 
& \multicolumn{4}{c}{CIFAR-10 (100 images)} \\
\cmidrule(lr){2-5} \cmidrule(lr){6-9}
Algorithm & SR & $\bar{R}^2$ & $\overline{\|\bm{\mu}^*\|_2}$ & $\bar{T}$ 
          & SR & $\bar{R}^2$ & $\overline{\|\bm{\mu}^*\|_2}$ & $\bar{T}$ \\
\midrule
\textit{Our Algo.} 
& $\mathbf{100\%}$ & $92\%$ & $\mathbf{4.68(0.97)}$ & $1247(314)$ 
& $\mathbf{100\%}$ & $99\%$ & $\mathbf{1.67(0.36)}$ & $597(212)$ \\

GS-PowerOpt 
& $100\%$ & $92\%$ & $4.73(0.93)$ & $1158(448)$ 
& $100\%$ & $99\%$ & $1.76(0.40)$ & $665(262)$ \\

ZOSGD 
& $100\%$ & $90\%$ & $5.25(0.88)$ & $2318(460)$ 
& $80\%$  & $99.6\%$ & $1.20(0.18)$ & $902(576)$ \\

ZOSLGHd 
& $100\%$ & $86\%$ & $6.09(0.97)$ & $2461(85)$ 
& $100\%$ & $99\%$ & $1.68(0.31)$ & $2401(443)$ \\

ZOSLGHr 
& $100\%$ & $83\%$ & $6.74(1.03)$ & $2374(532)$ 
& $100\%$ & $98\%$ & $2.48(0.60)$ & $846(68)$ \\

ZOAdaMM 
& $100\%$ & $71\%$ & $8.84(1.31)$ & $92(42)$ 
& $100\%$ & $61\%$ & $12.80(2.09)$ & $71(41)$ \\

STD-Htp 
& $48\%$  & $68\%$ & $8.58(1.02)$ & $917(589)$ 
& $50\%$  & $87\%$ & $7.91(1.61)$ & $811(574)$ \\

CMAES 
& $100\%$ & $93\%$ & $4.37(0.98)$ & $2485(22)$ 
& $100\%$ & $82\%$ & $8.61(1.21)$ & $2025(977)$ \\

\bottomrule
\end{tabular}
\end{table}

\textbf{MNIST and CIFAR-10} The experimental results on MNIST and CIFAR-10 are reported in Table \ref{tab:per-image-attacks-mnist-cifar}. They show that GS-PowerHP achieves a $100\%$ success rate and its perturbation size is the smallest among the smoothing-based algorithms.

\begin{table}[tbh!]
\caption{Per-image Attack on 100 Images from ImageNet. For each image attack, we set $\bm{\mu}_0=\bm{0}$ and $T_{total}=3,000$. The hyper-parameters for GS-PowerHP are selected as $N=8,\sigma_0=0.1$, $\alpha=1.0$, and $\beta=0.999$. Example perturbed images generated by GS-PowerHP can be found in Appendix \ref{appendix-adversarial}.}
\label{imagenet-attack}
\centering%
\begin{tabular}{ >{\centering\arraybackslash} p{3.5cm} | >{\centering\arraybackslash}p{0.9cm} >{\centering\arraybackslash}p{1.90cm}>{\centering\arraybackslash}p{1.9cm} >{\centering\arraybackslash}p{1.5cm} }  
\midrule
Algorithm & SR & $\bar{R}^2$ & $\overline{\|\bm{\mu}^*\|}$ & $\bar{T}$ \\
\toprule
\textit{Our Algo.} & $\mathbf{78}\%$      &  $95\%(7\%)$     & $\mathbf{35.8}(6.5)$          & $1389(542)$\\ 
ZOSLGHd     & $\mathbf{67}\%$  &$97\%(4\%)$     & $\mathbf{28.1}(3.3)$    & $2330(933)$\\
GS-PowerOpt     & $47\%$  &$95\%(3\%)$     & $36.5(2.8)$    & $1780(644)$\\
ZOSLGHr     & $43\%$  &$69\%(16\%)$     & $92.8(3.4)$    & $2167(801)$\\
ZOAdaMM  & $61\%$  &$61\%(27\%)$     & $104.9(13.6)$   &  $2896(454)$\\
ZOSGD  & $80\%$  &$-25\%(1.46)$     & $180(22)$    &  $3000(0)$\\
STD-Htp  &$0\%$   &NA  & NA    &NA\\
\bottomrule
\end{tabular}
\end{table}

\textbf{IMAGENET} The results on attacking images from the ImageNet are reported in Table \ref{imagenet-attack}. We do not test CMA-ES since it requires too much more computing resource than other algorithms for this huge-dimension task\footnote{According to the reported Python message, CMAES requires a RAM of 165 GB each iteration, which is beyond the capacity of our computer.}. 

The ImageNet results show that GS-PowerHP achieves the best overall tradeoff among the general-purpose smoothing-based zeroth-order methods. It substantially improves the success rate over GS-PowerOpt, from $47\%$ to $78\%$, while maintaining a comparable perturbation size and image similarity. Compared with ZOSLGHd, GS-PowerHP attains a higher success rate and requires fewer iterations, although ZOSLGHd produces slightly smaller perturbations. These results suggest that GS-PowerHP is particularly effective in the low-distortion regime for high-dimensional black-box adversarial attacks.

The superior performance of GS-PowerHP over other smoothing-based algorithms carries over to \textit{untargeted} image adversarial attacks against an adversarially robust classifier trained on ImageNet (ResNet18-based, \cite{salman2020adversarially}), with detailed results provided in Appendix \ref{sec-untargeted-IAA}.

\subsection{Ablation Study: Effects of Decaying-$\sigma$}
\label{sec-ablation-study}
We test effects of the decaying-$\sigma$ mechanism by comparing our algorithm with GS-PowerOpt through two kinds of experiments, one on a synthetic objective, and the other one on targeted image adversarial attacks.

\textbf{Synthetic Objective.} The maximization objective is set as
\begin{equation}
\label{two-log}
f(\bm{x}) = -\log(\|\bm{x}-\bm{m}_1 \|^2+10^{-5})-\log(\|\bm{x}-\bm{m}_2\|^2+10^{-2}),
\end{equation}
where all entries of $\bm{m}_1\in\mathbb{R}^d$ equal $-0.5$ and all entries of $\bm{m}_2\in\mathbb{R}^d$ equal $0.5$. Note that $\bm{x}^*=\bm{m}_1$ for this objective.

For each $\sigma$ from a pre-selected set $\{0.1,0.5,1.0,2.0,3.0\}$, we perform 100 trials of GS-PowerOpt to maximize $f(\bm{x})$ in (\ref{two-log}). Then, in Table \ref{sigma-effect-table}, we report the average of $\{f(\bm{\mu}^*_n)\}_{n=1}^{100}$ and the average of $\{\text{MSE}(\bm{\mu}^*_n)\}_{n=1}^{100}$, where $\text{MSE}(\bm{\mu}^*)$ denotes the $L_2$ norm $\| \bm{\mu}^*-\bm{x}^*\|^2/d$. For GS-PowerHP, we perform the same experiment, with $\sigma_0=3$, $b=0$, and the decaying factor $\beta$ being the solution to $3\beta^{1000}=0.1$.

From the table, we see that GS-PowerHP is able to produce a solution $\bm{\mu}^*$ that is closer to $\bm{x}^*$ than all the solutions found by GS-PowerOpt, indicating the advantage of the $\sigma$-decreasing mechanism of GS-PowerHP over the fixed-$\sigma$ method of GS-PowerOpt. 

Note that although $\sigma=0.1$ is associated with a comparatively large MSE, the corresponding average fitness is good. This is because GS-PowerOpt's performance is volatile with this $\sigma$ value. Some trials produces $\bm{\mu}^*$ close to $\bm{x}^*$ and in turn lead to a high fitness.

\begin{table}[t]
\caption{The performances of GS-PowerOpt and GS-PowerHP on optimizing $f(\bm{x})$ in (\ref{two-log}). We set $N=1$ and $T=1,000$. The reported $\sigma$ value for GS-PowerHP is the initial $\sigma$ value (i.e., $\sigma_0$).}
\label{sigma-effect-table}
\centering%
\begin{tabular}{ p{2.5cm}| p{1.0cm}  p{1.0cm}  p{1.0cm} | p{1.0cm}  p{1.0cm}  p{1.0cm}  }
\toprule
 \multirow{2}{*}{Algo.}& \multicolumn{3}{c|}{$d=3$} & \multicolumn{3}{c}{$d=5$} \\
 & $\sigma$ & $f(\bm{\mu}^*)$ & MSE & $\sigma$ & $f(\bm{\mu}^*)$ & MSE \\
\midrule
\multirow{5}{*}{\begin{tabular}{@{}l@{}}GS-PowerOpt\end{tabular}}&3.0  &1.54 &0.33 &3.0   & $-$0.11  &  0.31  \\
&2.0  &1.72 &0.29 & 2.0   &$-$0.11       &  0.28   \\
&1.0  &4.18 &0.01 &1.0  &0.50       &  0.20  \\
&0.5  &6.23 &0.42 &0.5    &4.09        &  0.49  \\
&0.1  &5.74 &0.56 &0.1  &4.69        &0.56   \\
\midrule
Our Algo. &3.0  &\textbf{7.26} &\textbf{0.23} &0.1   &\textbf{5.44}         &\textbf{0.29} \\
\bottomrule
\end{tabular}
\end{table}
\textbf{Image Adversarial Attacks.}
\label{sample-efficiency-exp}
To isolate the effect of the decaying-$\sigma$ mechanism, we compare GS-PowerHP with GS-PowerOpt under identical experimental conditions. In particular, both algorithms share the same initial point $\bm{\mu}_0=\mathbf{0}$, learning rate $\alpha$, power parameter $N$, and adversarial attack setting described in Section~\ref{Experiment-IAA}. The only difference between the two methods is the use of a decaying smoothing radius in GS-PowerHP versus a fixed smoothing radius in GS-PowerOpt. For fairness, the fixed $\sigma$ used in GS-PowerOpt is set equal to the initial smoothing parameter $\sigma_0$ of GS-PowerHP (where we set $b=0$).

For each of the 100 randomly selected images, we independently and uniformly sample the power parameter $N$, learning rate $\alpha$, and smoothing parameter $\sigma$ from pre-specified intervals. The bounds of these intervals coincide with the hyperparameter values used in the targeted IAA in Section~\ref{Experiment-IAA}. For example, on MNIST we sample $N\in[0.1,0.5]$ and $\sigma\in[0.05,0.1]$. Additional implementation details and full hyperparameter ranges are provided in the code.

The results are summarized in Table~\ref{decay-effect-exp}. Across all datasets, GS-PowerHP consistently requires fewer iterations (i.e., a smaller average iteration count $\bar{T}$) to reach a solution of higher quality than GS-PowerOpt. Since both methods use the same number of Monte Carlo samples per iteration, these results provide empirical evidence that the decaying-$\sigma$ mechanism in GS-PowerHP leads to more efficient optimization dynamics compared to a fixed-$\sigma$ scheme.

\begin{table}[t]
\caption{Per-image adversarial attack results on 100 images (least likely targeted). For each image, GS-PowerOpt and GS-PowerHP share the same initial point $\boldsymbol{\mu}_0=\mathbf{0}$, power parameter $N$, and fixed learning rate $\alpha$, which are uniformly drawn from a pre-specified range. All other experimental settings follow Subsection~\ref{Experiment-IAA}.}
\label{decay-effect-exp}
\centering%
\begin{tabular}{ >{\centering\arraybackslash} p{1.7cm} | >{\centering\arraybackslash}p{2.5cm} >{\centering\arraybackslash}p{1.0cm}>{\centering\arraybackslash}p{1.6cm} >{\centering\arraybackslash}p{2.0cm} }  
\midrule
Imageset & Algorithm& SR &  $\overline{\|\bm{\mu}^*\|_2}$ & $\bar{T}$ \\
\toprule
\multirow{2}{*}{MNIST} &\textbf{Our Algo.} &  $\mathbf{100}\%$     & $\mathbf{4.86}(0.9)$          & $\mathbf{979}(269)$\\ 
  &GS-PowerOpt   & $100\%$      & $4.93(0.9)$    & $1,036(304)$\\
\hline
\multirow{2}{*}{CIFAR10} &\textbf{Our Algo.} &  $\mathbf{100}\%$     & $\mathbf{2.22}(0.4)$          & $\mathbf{516}(162)$\\ 
  &GS-PowerOpt   & $98\%$      & $2.50(0.5)$    & $651(248)$\\
\hline
\multirow{2}{*}{ImageNet} &\textbf{Our Algo.} &  $\mathbf{80}\%$     & $\mathbf{40.5}(7.2)$          & $\mathbf{1,114}(358)$\\ 
  &GS-PowerOpt   & $61\%$      & $45.3(7.3)$    & $1,362(370)$\\
\bottomrule
\end{tabular}
\end{table}

\subsection{Hyper-parameter Analysis}


\subsubsection{Effects of Varying $\beta$ and $b$}
\label{subsec-beta-b-sensitivity}

GS-PowerHP uses the smoothing schedule
\[
\sigma_t=\sigma_0\beta^t+b,
\]
where $\beta\in(0,1)$ controls the decay speed and $b\geq 0$ specifies the limiting smoothing scale. This schedule separates the two roles of $\sigma$: larger early-stage values of $\sigma_t$ support global exploration, whereas smaller late-stage values support local refinement.

Corollary~\ref{convergence-rate} provides a theoretical explanation for this tradeoff. The finite-time stationarity rate improves when the current smoothing scale $\sigma_0\beta^T+b$ is larger, which explains why using a small smoothing scale too early may slow down the finite-time optimization process. In contrast, as $T$ increases, the schedule gradually approaches the limiting scale $b$, thereby supporting more accurate late-stage refinement. Hence, compared with a fixed small-$\sigma$ variant that uses the limiting scale $b$ from the beginning, GS-PowerHP with $\beta>0$ can benefit from larger smoothing radii during finite-time optimization while still approaching the same limiting smoothing scale.

Table~\ref{Tab:effects-b-beta} illustrates this exploration--refinement tradeoff. The extreme case $\beta=0$ uses the small smoothing scale almost from the beginning and fails to reach the gradient threshold, indicating insufficient exploration. In contrast, $\beta=1$ corresponds to a fixed large smoothing scale, which reaches the gradient threshold quickly but attains poor final fitness. Intermediate decay rates achieve substantially better final fitness, suggesting that a gradual decrease of $\sigma_t$ provides a more effective balance between exploration and refinement. The sensitivity with respect to $b$ shows a related pattern: smaller limiting scales can improve final solution quality, while larger values of $b$ tend to reduce the number of iterations needed to reach the gradient threshold.

\begin{table}
\caption{Sensitivity analysis of $\beta$ and $b$ in GS-PowerHP for maximizing the objective in Eq.~(\ref{two-log}). In the left panel, $\beta$ is varied while other parameters are fixed; in the right panel, $b$ is varied while other parameters are fixed. The cases $\beta=0$ and $\beta=1$ correspond to fixed-smoothing variants: $\beta=0$ uses the limiting scale $b$ after the initial step, while $\beta=1$ uses the fixed scale $\sigma_0+b$. Here, $\bar{f}(\bm{\mu}^*)$ denotes the average fitness of the best solution found over 100 trials, and $\bar{t}{\epsilon}$ denotes the average number of iterations required for the scaled gradient estimate to satisfy $\|\hat{\nabla}F_{N,\sigma_{t+1}}(\bm{\mu}_t)/\sigma^2_{t+1}\|\leq \epsilon$, where $\epsilon=10^{-6}$. ``NA'' means that the threshold is not reached within the prescribed iteration budget.}
\label{Tab:effects-b-beta}
\centering%
\small
\begin{tabular}{ p{0.8cm}  p{0.8cm} p{0.8cm} | p{0.8cm} p{0.8cm} p{0.8cm}}  
\midrule
$\beta$ & $\bar{f}(\bm{\mu}^*)$ & $\bar{t}_{\epsilon}$ & $b$ & $\bar{f}(\bm{\mu}^*)$ & $\bar{t}_{\epsilon}$\\
\hline
0.0      &   4.12   &    NA      &	0	&	8.75 	&	153 	\\     
0.997	&	8.93 	&	241 	&	0.05	&	8.65 	&	137 	\\
0.999	&	8.54 	&	208 	 	&   0.1	 	&8.64  	&	93 	\\
1.0	&	4.44 	&	54 			&  0.15  & 	8.38 	&	77 	\\
\bottomrule
\end{tabular}
\end{table}

\subsubsection{Hyper-parameter Sensitivity Analysis}
We evaluate the sensitivity of GS-PowerHP to its hyper-parameters using the targeted adversarial attack experiment on CIFAR-10. Specifically, for a fixed set of 50 randomly selected CIFAR images, we uniformly sample 20 hyper-parameter configurations from the ranges $N\in [0.1,0.3]$, $\beta\in[0.998,0.9995]$, and $b\in[0.0,0.005]$. For each configuration, we repeat the same attack procedure as in Table~\ref{tab:per-image-attacks-mnist-cifar}. Averaged over the 20 configurations, GS-PowerHP achieves $SR=100\%$ and $\overline{\|\bm{\mu}^*\|_2}=1.70$, which is close to the result reported in Table~\ref{tab:per-image-attacks-mnist-cifar}. These results suggest that GS-PowerHP is reasonably stable with respect to hyper-parameter choices within the tested ranges. The standard deviation of perturbation norm across the 20 configurations is small (i.e., 0.026), further indicating that the method is not highly sensitive to moderate hyper-parameter variations.

\subsection{Summary of Experimental Results}
Among the smoothing-based methods evaluated, GS-PowerHP consistently outperforms the competing smoothing-based approaches in our experiments (Table \ref{tab:combined-ackley-rastrigin}, \ref{tab:per-image-attacks-mnist-cifar}, and \ref{imagenet-attack}). The ablation study's results (Table \ref{sigma-effect-table} and \ref{decay-effect-exp}) indicate that the decaying-$\sigma$ mechanism in GS-PowerHP effectively improves both the speed and quality of the solution searching process. The hyper-parameter sensitivity analysis indicates that GS-PowerHP is reasonably stable with respect to hyper-parameter choices within the tested ranges.

Although CMA-ES, a strong non-smoothing-based evolutionary method, achieves better performance on the Ackley and MNIST tasks, GS-PowerHP performs better on the Rastrigin and CIFAR-10 tasks, suggesting that its advantages become more pronounced in challenging high-dimensional settings. In addition, compared with CMA-ES, which typically incurs higher computational and memory costs due to covariance adaptation, GS-PowerHP uses a simpler stochastic zeroth-order update based on sampled perturbations and avoids maintaining a full covariance matrix. This makes GS-PowerHP computationally attractive for high-dimensional black-box optimization problems, especially when scalability is a major concern.

\section{Conclusion}
In this work, we introduce GS-PowerHP, a novel zeroth-order method for general non-convex optimization. By combining a power-transformed Gaussian smoothed objective with an principled decreasing scaling parameter $\sigma$, our approach delivers superior theoretical convergence guarantees and empirical performance compared to existing smoothing-based methods. Extensive experiments demonstrate that GS-PowerHP consistently outperforms other smoothing-based competing algorithms, even in extremely high-dimensional settings such as adversarial image attacks on ImageNet.

\clearpage  

\section*{Acknowledgements}
None.

%
%
\bibliographystyle{splncs04}
\bibliography{Literature}

\newpage
\appendix
\onecolumn
\section{Appendix}

\subsection{Convergence Guarantees of Existing Smoothing-based Zeroth-Order Methods}
\begin{table}[!htb]
\caption{Convergence Guarantees}
\label{guarantees}
\centering%
\begin{tabular}{ p{4cm} | p{3cm} | p{3cm} }
\midrule
Methods & Convergence Target & Theoretical Iteration Complexity \\
\toprule
Standard Homotopy \cite{Hazan2016} & Global & $\mathcal{O}(d^2\sigma^{-2}\epsilon^{-4})$\\
ZO-SLGH \cite{Iwakiri2022}   & Local & $\mathcal{O}(d\epsilon^{-2})$\\
GS-PowerOpt \cite{GS-PowerOpt}   & Global & $\mathcal{O}((d^2\varepsilon^{{-1}}\sigma^{-2})^{\frac{2}{1-2\gamma}})$\\
GS-PowerHp (Ours. )  & Global & $\mathcal{O}((d^2\varepsilon^{{-1}}b^{-2})^{\frac{2}{1-2\gamma}})$\\
\bottomrule
\end{tabular}
\end{table}

\subsection{Proof Details}
\label{appendix-proof}

\textbf{Theorem }\ref{sigma-iter-complexity}.
\textit{Under the conditions specified in \cite[Corollary 3.9]{GS-PowerOpt}, for some constant $U>0$ that depends on $e^{Nf(\bm{x}^*)}$ and the learning rate parameter $\gamma$, whenever $\sigma\in (0, U)$, the iteration complexity of GS-PowerOpt is $O(\sigma^{-2}d^2\varepsilon^{-1})^{2/(1-2\gamma)}$. }

\begin{proof}
From \cite[Corollary 3.9]{GS-PowerOpt}, the iteration complexity of GS-PowerOpt is $T=(C_1C_2 d^2 \varepsilon^{-1})^{2/(1-2\gamma)}$, where $\gamma\in(0,1/2)$, $C_2=\max\{1, 2/C_1\}$, $C_1=C_0(1-2\gamma)\sigma^{-2}$, and 
$$C_0=e^{Nf(\bm{x}^*)} - F_{N,\sigma}(\bm{\mu}_0) + 2e^{3Nf(\bm{x}^*)}\sum_{t=1}^\infty t^{-(1+2\gamma)}.$$ 
Therefore, when $\sigma\in(0, \sqrt{C_0(1-2\gamma)/2})$, $2/C_1 = 2(C_0(1-2\gamma))^{-1}\sigma^2 < 1$ and $C_2=1$. Hence, 
$$T=(C_1 d^2 \varepsilon^{-1})^{2/(1-2\gamma)} =O((\sigma^{-2} d^2 \varepsilon^{-1})^{2/(1-2\gamma)}).$$
\end{proof}

\noindent\textbf{Proposition }\ref{unique-max}. \textit{Assume that $\bm{y}^*$ is any maximizer of $f:\mathbb{R}^d\to\mathbb{R}$, $f$ is second-order continuously differentiable, and there exists positive numbers $R$, $m$, $M$, and $L_f$ such that 
\begin{enumerate}
\item $|f(\bm{x}_1)-f(\bm{x}_2)|\leq L_f \|\bm{x}_1-\bm{x}_2\|$ for any $\bm{x}_1,\bm{x}_2\in\mathbb{R}^d$; 
\item $\bm{v}'\nabla^2 f(\bm{x})\bm{v} <-m \|\bm{v}\|^2$ for all $\bm{v}\in\mathbb{R}^d$ and any $\bm{x}\in B_{R}(\bm{y^*}):=\{\bm{z}\in\mathbb{R}^d | \|\bm{z}-\bm{y}^*\|<R\}$;
\item $f(\bm{x})< M$ for all $\bm{x}\in\mathbb{R}^d$; 
\item and $\|\nabla^2 f(\bm{x})\|<M$ for all $\bm{x}\in\mathbb{R}^d$, where $\|\nabla^2 f(\bm{x})\|:=\max_{\|\bm{v}\|=1}\|\nabla^2 f(\bm{x})\bm{v}\|$ denotes the spectral norm.
\end{enumerate}
Then, for any $N>0$, we have the following results.
\begin{itemize}
\item There exists $R_1\in (0,R)$ such that $\sup_{\|\bm{x}-\bm{y}^*\|<R_1}\| \nabla f(\bm{x}) \|<\sqrt{m/2N}$.
\item For any $r\in (0,R_1/2)$, $F_{N,\sigma}$ has a unique maximizer $\bm{z}^*$ in $B_{r}(\bm{y^*})$ such that $F_{N,\sigma}(\bm{z}^*)>F_{N,\sigma}(\bm{x})$ for any $\bm{x}\in B_{r}(\bm{y}^*)$ and $\bm{x}\neq\bm{z}^*$, if 
$$0<\sigma< \min\left\{\frac{R_1}{2\sqrt{d}+2\sqrt{2\log((m_1+M_1)/m_1)}}, \frac{\Delta_r}{2 e^{NM} N L_f  \mathbb{E}_{\bm{\epsilon}\sim\mathcal{N}(\bm{0}, I_d)} \| \bm{\epsilon}\|} \right\},$$
where $m_1:= Ne^{Nf_c}(m/2)>0$, $f_c$ is a lower bound of $f$ on the closed set $\overline{B_{R_1}(\bm{y}^*)}\in\mathbb{R}^d$, $M_1:= N^2 e^{NM} L_f^2 + Ne^{NM}M$, and $\Delta_r:= f_N(\bm{y^*})- \max_{\bm{x}\in\partial B_r(\bm{y}^*)} f_N(\bm{x})>0$.
\end{itemize}
}


\begin{proof}
Our proof consists of three parts: (1) We prove that $f_N(\bm{x}):=e^{Nf(\bm{x})}$ has properties similar to Point 1 and 2 of $f$, with which we prove (2) $F_{N,\sigma}$ is strictly concave in $B_{R_1/2}(\bm{y}^*)$ for some $R_1\in(0,R)$, if $\sigma>0$ is less than a threshold. (3) We prove that for any $r\in (0,R_1/2)$, if $\sigma>0$ is less than some threhold, then $F_{N,\sigma}$ has a unique maximum $\bm{z}^*$ in $B_{r}(\bm{y^*})$ such that $F_{N,\sigma}(\bm{z}^*)>F_{N,\sigma}(\bm{x})$ for any $\bm{x}\in B_{r}(\bm{y}^*)$ and $\bm{x}\neq\bm{z}^*$.

\textit{Proof Part (1)}. Since $f$ is assumed to be second-order differentiable, its gradient $\nabla f$ is continuous at $\bm{y}^*$. Hence, for $\varepsilon=\sqrt{m/2N}$, there exists $R_1\in(0, R)$ such that whenever $\|\bm{x}-\bm{y}^*\|<R_1$, $\sqrt{m/2N}>\|\nabla f (\bm{x}) - \nabla f(\bm{y}^*) \|=\|\nabla f (\bm{x})\|$.

For any $\bm{x}\in B_{R_1}(\bm{y}^*)$,
\begin{equation}
\begin{split}
\bm{v}'\nabla^2 f_N(\bm{x}) \bm{v} &= \bm{v}'\cdot \nabla (Ne^{Nf(\bm{x})}\nabla f(\bm{x}))\bm{v},\quad \bm{v}' \text{ denotes transpose of }\bm{v},\\
 &= N^2e^{Nf(\bm{x})}\bm{v}'\nabla f (\bm{x})(\nabla f(\bm{x}))'\bm{v}+Ne^{Nf(\bm{x})}\bm{v}'\nabla^2 f(\bm{x})\bm{v} \\
&\leq N^2e^{Nf(\bm{x})}\|\bm{v}'\nabla f(\bm{x})\|^2+Ne^{Nf(\bm{x})}(-m)\|\bm{v}\|^2\\
&\leq N^2e^{Nf(\bm{x})}\|\nabla f(\bm{x})\|^2\|\bm{v}\|^2+Ne^{Nf(\bm{x})}(-m)\|\bm{v}\|^2\\
&\leq -Ne^{Nf_c}(m/2)\|\bm{v}\|^2\quad\text{since } \sqrt{m/2N}>\|\nabla f (\bm{x})\|,\\
&= -m_1\|\bm{v}\|^2,
\label{R1-property}
\end{split}
\end{equation}
where $f_c$ is a lower bound of $f$ on the closed and bounded set $\overline{B_{R_1}(\bm{y}^*)}\in\mathbb{R}^d$ and 
$$m_1:= Ne^{Nf_c}(m/2)>0.$$

Since $f$ is Lipschitz, its gradient $\nabla f$ has a bounded norm over $\mathbb{R}^d$ and we denote this bound with $L_f$. For any $\bm{x}\in\mathbb{R}^d$, 
\begin{align*}
\| \nabla^2 f_N(\bm{x}) \| &= \|  N^2e^{Nf(\bm{x})}\nabla f (\bm{x})(\nabla f(\bm{x}))'+Ne^{Nf(\bm{x})}\nabla^2 f(\bm{x}) \|\\
&\leq N^2e^{Nf(\bm{x})} \| \nabla f (\bm{x})(\nabla f(\bm{x}))'\|+Ne^{Nf(\bm{x})}\|\nabla^2 f(\bm{x}) \|\\
&= N^2e^{Nf(\bm{x})} \| \nabla f (\bm{x})\|^2+Ne^{Nf(\bm{x})}\|\nabla^2 f(\bm{x}) \|\\
&\leq N^2 e^{NM} L_f^2 + Ne^{NM}\| \nabla^2 f(\bm{x}) \|\\
&\leq N^2 e^{NM} L_f^2 + Ne^{NM} M\\
& = M_1,\quad\text{where }M_1:= N^2 e^{NM} L_f^2 + Ne^{NM}M>0,
\end{align*}
which implies for all $\bm{v}\in\mathbb{R}^d$ that
\begin{equation}
\begin{split}
\bm{v}' \nabla^2 f_N(\bm{x}) \bm{v} & \leq \| \bm{v}' \nabla^2 f_N(\bm{x}) \bm{v} \| \\
& \leq  \| \bm{v}'\|\cdot \left\| \nabla^2 f_N(\bm{x}) \frac{\bm{v}}{\|\bm{v}\|}\right\| \cdot \| \bm{v}\| \\
& \leq  M_1\|\bm{v} \|^2.
\label{Rd-property}
\end{split}
\end{equation}
\textit{Proof Part (2)}. We next show that $F_{N,\sigma}$ is strictly concave in $B_{R_1/2}(\bm{y}^*)$. For any $\bm{\mu}\in B_{R_1/2}(\bm{y}^*)$ and any $\bm{v}\in\mathbb{R}^d$,
\begin{align*}
\bm{v}'\nabla^2 F_{N,\sigma}(\bm{\mu}) \bm{v}= & \mathbb{E}_{\bm{\epsilon}\sim\mathcal{N}(0,I_d)}[\bm{v}'\nabla^2f_N(\bm{\mu}+\sigma\bm{\epsilon})\bm{v}] \text{ by dominated convergence,}\\
\leq& \mathbb{P}\{\bm{\mu}+\sigma\bm{\epsilon}\in B_{R_1}(\bm{y}^*)\} (-m_1\|\bm{v} \|^2),\quad\text{by (\ref{R1-property})}, \\
& +  \mathbb{P}\{\bm{\mu}+\sigma\bm{\epsilon}\notin B_{R_1}(\bm{y}^*)\} M_1\|\bm{v} \|^2,\quad\text{by (\ref{Rd-property})},\\ 
= &(1-p_{\sigma})  (-m_1\|\bm{v} \|^2) + p_{\sigma} M_1\|\bm{v} \|^2\\
= &-m_1\|\bm{v} \|^2 + p_{\sigma} (m_1+M_1)\|\bm{v} \|^2,
\end{align*}
where $p_{\sigma}:=\mathbb{P}\{\bm{\mu}+\sigma\bm{\epsilon}\notin B_{R_1}(\bm{y}^*)\}\leq \mathbb{P}\{\|\bm{\epsilon}\|\geq R_1/(2\sigma)\}$ since the former event implies the latter event. Also, $\lim_{\sigma\to 0}\mathbb{P}\{\|\bm{\epsilon}\|\geq R_1/(2\sigma)\}=0$ since $\bm{\epsilon}$ is standard Gaussian. This implies that if $\sigma$ is sufficiently small, $\bm{v}'\nabla^2 F_{N,\sigma}(\bm{\mu}) \bm{v}<0$, which further implies that $F_{N,\sigma}(\bm{\mu})$ is strictly concave in $B_{R_1/2}(\bm{y}^*)$. Specifically, a sufficient upper bound for $\sigma$ is $\frac{R_1}{2\sqrt{d}+2\sqrt{2\log((m_1+M_1)/m_1)}}$ (proof postponed to the last). 

\textit{Proof Part (3)}. Now, let $r$ be any number in $(0,R_1/2)$. We prove the existence of a maximizer of $F_{N,\sigma}$ in $B_r(\bm{y}^*)$. Since $F$ is continuous, it attains a maximum in the compact set $\overline{B_{r}(\bm{y}^*)}$. That is, there exists some $\bm{z}^*\in \overline{B_{r}(\bm{y}^*)}$ such that
$$ F_{N,\sigma}(\bm{z}^*) = \max_{\bm{x}\in \overline{B_{r}(\bm{y}^*)}} F_{N,\sigma}(\bm{x}).$$
We prove that $\bm{z}^*$ does not fall on the boundary $\partial B_r(\bm{y}^*)$, which implies it must be in $B_r(\bm{y}^*)$ and implies that $\bm{z}^*$ is a maximizer of $F_{N,\sigma}$.

In $B_R(\bm{y}^*)$, since $\nabla^2 f$ is strictly negative definite, $f$ must be strictly concave. Hence, since $r<R$, $f(\bm{y^*})>\max_{\bm{x}\in\partial B_r(\bm{y}^*)} f(\bm{x})$, where the latter exists since $f$ is continuous on the compact set $\partial B_{r}(\bm{y}^*):=\{ \bm{x}\in\mathbb{R}^d | \|\bm{x}-\bm{y}^*\|=r\}$. Hence,
$$\Delta_r:= f_N(\bm{y^*})- \max_{\bm{x}\in\partial B_r(\bm{y}^*)} f_N(\bm{x})>0.$$
Also, for any $\bm{\mu}\in\mathbb{R}^d$,
\begin{align*}
 |F_{N,\sigma}(\bm{\mu})-e^{Nf(\bm{\mu})}| & =| \mathbb{E}_{\bm{\epsilon}\sim\mathcal{N}(\bm{0},I_d)} [e^{Nf(\bm{\mu}+\sigma\bm{\epsilon})}-e^{Nf(\bm{\mu})}]| \\
& \leq e^{NM} N \mathbb{E}_{\bm{\epsilon}\sim\mathcal{N}(\bm{0}, I_d)} |f(\bm{\mu}+\sigma \bm{\epsilon})-f(\bm{\mu})|,\quad\text{by MVT on }e^{x},\\
& \leq e^{NM} N L_f \sigma \mathbb{E}_{\bm{\epsilon}\sim\mathcal{N}(\bm{0}, I_d)} \| \bm{\epsilon} \|.
\end{align*}
This implies
\begin{enumerate}
\item $F_{N,\sigma}(\bm{y}^*)\geq f_N(\bm{y}^*)-e^{NM} N L_f \sigma \mathbb{E}_{\bm{\epsilon}\sim\mathcal{N}(\bm{0}, I_d)} \| \bm{\epsilon}\|$;
\item $\max_{\bm{\mu}\in \partial B_r(\bm{y}^*)}F_{N,\sigma}(\bm{\mu})\leq  \max_{\bm{\mu}\in \partial B_r(\bm{y}^*)} f_N(\bm{\mu}) + e^{NM} N L_f \sigma \mathbb{E}_{\bm{\epsilon}\sim\mathcal{N}(\bm{0}, I_d)} \| \bm{\epsilon} \|$.
\end{enumerate}
Hence, subtracting the second inequality from the first gives
\begin{align*}
&F_{N,\sigma}(\bm{y}^*)-\max_{\bm{\mu}\in \partial B_r(\bm{y}^*)}F_{N,\sigma}(\bm{\mu})\\
\geq & f_N(\bm{y}^*)-\max_{\bm{\mu}\in \partial B_r(\bm{y}^*)} f_N(\bm{\mu}) - 2 e^{NM} N L_f \sigma \mathbb{E}_{\bm{\epsilon}\sim\mathcal{N}(\bm{0}, I_d)} \| \bm{\epsilon} \|\\
= & \Delta_r - 2 e^{NM} N L_f \sigma \mathbb{E}_{\bm{\epsilon}\sim\mathcal{N}(\bm{0}, I_d)} \| \bm{\epsilon} \|\\
> & 0, \text{ if }\sigma< \frac{\Delta_r}{2 e^{NM} N L_f  \mathbb{E}_{\bm{\epsilon}\sim\mathcal{N}(\bm{0}, I_d)} \| \bm{\epsilon}\|}.
\end{align*}
Hence, $\bm{z}^*\notin \partial B_r(\bm{y}^*)$, and it must be in the open set $B_r(\bm{y}^*)$, which further implies that it is a (possibly local) maximizer of $F_{N,\sigma}$.

In sum, whenever $\sigma<\frac{R_1}{2\sqrt{d}+2\sqrt{2\log((m_1+M_1)/m_1)}}$, $F_{N,\sigma}$ is strictly concave in $B_{R_1/2}(\bm{y}^*)$. For any $r\in (0,R_1/2)$, whenever $\sigma< \frac{\Delta_r}{2 e^{NM} N L_f  \mathbb{E}_{\bm{\epsilon}\sim\mathcal{N}(\bm{0}, I_d)} \| \bm{\epsilon}\|}$, $F_{N,\sigma}$ has a possibly local maximizer in $B_{r}(\bm{y}^*)$. Therefore, whenever $\sigma$ is less than both of the two thresholds, $F_{N,\sigma}$ has a unique maximizer $\bm{z}^*$ in $B_{r}(\bm{y}^*)$ such that $F_{N,\sigma}(\bm{z}^*)>F_{N,\sigma}(\bm{x})$ for any $\bm{x}\in B_{r}(\bm{y}^*)$ and $\bm{x}\neq\bm{z}^*$.

Finally, we prove that an upper bound of $\sigma$ such that $p_{\sigma}<m_1/(M_1+m_1)$ is\\ 
$\frac{R_1}{2\sqrt{d}+2\sqrt{2\log((m_1+M_1)/m_1)}}$. Set
\[
t:=\log\left(\frac{m_1+M_1}{m_1}\right)>0.
\]
Using the standard Gaussian norm concentration inequality,
\[
\mathbb{P}\left\{\|\bm{\epsilon}\|\geq \sqrt{d}+\sqrt{2t}\right\}
\leq e^{-t}.
\]
Since
\[
e^{-t}
=
\frac{m_1}{m_1+M_1},
\]
we obtain
\[
\mathbb{P}\left\{\|\bm{\epsilon}\|\geq 
\sqrt{d}+\sqrt{2\log\left(\frac{m_1+M_1}{m_1}\right)}
\right\}
\leq
\frac{m_1}{m_1+M_1}.
\]
Therefore, if
\[
\frac{R_1}{2\sigma}
>
\sqrt{d}+\sqrt{2\log\left(\frac{m_1+M_1}{m_1}\right)},
\]
or equivalently,
\[
\sigma
<
\frac{R_1}{
2\sqrt{d}+2\sqrt{2\log\left(\frac{m_1+M_1}{m_1}\right)}
},
\]
then
\[
\mathbb{P}\left\{\|\bm{\epsilon}\|\geq \frac{R_1}{2\sigma}\right\}
\leq
\mathbb{P}\left\{\|\bm{\epsilon}\|\geq 
\sqrt{d}+\sqrt{2\log\left(\frac{m_1+M_1}{m_1}\right)}
\right\}
\leq
\frac{m_1}{m_1+M_1}.
\]
Therefore,
$$ p_{\sigma} \leq \mathbb{P}\left\{\|\bm{\epsilon}\|\geq \frac{R_1}{2\sigma}\right\} \leq \frac{m_1}{m_1+M_1}.$$
\end{proof}

\noindent\textbf{Proposition \ref{mu-dist}}. \textit{Assume that $f(\bm{x})<M$ for any $\bm{x}\in\mathbb{R}^d$ and $f$ is Lipschitz with constant $L_f>0$. Let $\bm{y}^*$ be a locally sharp maximum of $f$. That is, there are some $\kappa, p>0$ and a bounded neighborhood $U:=\{\bm{x}\in\mathbb{R}^d | \|\bm{x}-\bm{y}^*\|<r \}$ of $\bm{y}^*$ such that $f(\bm{y}^*) - f(\bm{x})\geq \kappa \|\bm{y}^*-\bm{x}\|^p$ for any $\bm{x}\in U$. For any $N, \sigma>0$, if there exists $\bm{\mu}^*_{N,\sigma}\in U$ such that $F_{N,\sigma}(\bm{\mu}^*_{N,\sigma})>F_{N,\sigma}(\bm{x})$ for any $\bm{x}\in U$ and $\bm{x}\neq \bm{\mu}^*_{N,\sigma}$, then $ \|\bm{y}^*-  \bm{\mu}^*_{N,\sigma} \| \leq \left( \frac{2 e^{N(M-f_U)} L_f \sigma \mathbb{E}_{\bm{\epsilon}\sim\mathcal{N}(\bm{0}, I_d)}\| \bm{\epsilon} \|}{N\kappa} \right)^{1/p}$, where $f_U>-\infty$ denotes the lower bound of $f$ on $U$.} 
\begin{proof}
Let $\bm{\epsilon}\in\mathbb{R}^d$ follow a standard Gaussian distribution. For any $\bm{\mu}$,
\begin{align*}
 |F_{N,\sigma}(\bm{\mu})-e^{Nf(\bm{\mu})}| & =| \mathbb{E}_{\bm{\epsilon}\sim\mathcal{N}(\bm{0},I_d)} [e^{Nf(\bm{\mu}+\sigma\bm{\epsilon})}-e^{Nf(\bm{\mu})}]| \\
& \leq e^{NM} N \mathbb{E}_{\bm{\epsilon}\sim\mathcal{N}(\bm{0}, I_d)} |f(\bm{\mu}+\sigma \bm{\epsilon})-f(\bm{\mu})|,\quad\text{by MVT on }e^{x},\\
& \leq e^{NM} N L_f \sigma \mathbb{E}_{\bm{\epsilon}\sim\mathcal{N}(\bm{0}, I_d)} \| \bm{\epsilon} \|.
\end{align*}
This inequality implies 
\begin{itemize}
\item if $\bm{\mu}=\bm{\mu}^*_{N,\sigma}$, then 
$$F_{N,\sigma}(\bm{\mu}^*_{N,\sigma}) \leq e^{Nf(\bm{\mu}^*_{N,\sigma})} + e^{NM} N L_f \sigma \mathbb{E}_{\bm{\epsilon}\sim\mathcal{N}(\bm{0}, I_d)}\| \bm{\epsilon} \|.$$
\item if $\bm{\mu}=\bm{y}^*$, then 
$$e^{Nf(\bm{y}^*)} - e^{NM} N L_f \sigma \mathbb{E}_{\bm{\epsilon}\sim\mathcal{N}(\bm{0},I_d)}\| \bm{\epsilon} \|\leq F_{N,\sigma}(\bm{y}^*).$$
\end{itemize}
\end{proof}
From the above two inequalities and $F_{N,\sigma}(\bm{y}^*)\leq F_{N,\sigma}(\bm{\mu}^*_{N,\sigma})$, we have
$$ e^{Nf(\bm{y}^*)} - e^{Nf(\bm{\mu}^*_{N,\sigma})} \leq 2 e^{Nf(\bm{x}^*)} N L_f \sigma \mathbb{E}_{\bm{\epsilon}\sim\mathcal{N}(\bm{0}, I_d)}\| \bm{\epsilon} \|,$$
which further implies
\begin{align*}
2 e^{NM} N L_f \sigma \mathbb{E}_{\bm{\epsilon}\sim\mathcal{N}(\bm{0}, I_d)}\| \bm{\epsilon} \| & \geq e^{Nf(\bm{y}^*)} - e^{Nf(\bm{\mu}^*_{N,\sigma})} \\
&\geq e^{Nf_U} N (f(\bm{y}^*)-f(\bm{\mu}^*_{N,\sigma})),\quad\text{from MVT on }e^x,\\
&\geq e^{Nf_U} N \kappa \|\bm{y}^*-  \bm{\mu}^*_{N,\sigma} \|^p.
\end{align*}
This implies that
$$ \|\bm{y}^*-  \bm{\mu}^*_{N,\sigma} \| \leq \left( \frac{2 e^{N(M-f_U)} L_f \sigma \mathbb{E}_{\bm{\epsilon}\sim\mathcal{N}(\bm{0}, I_d)}\| \bm{\epsilon} \|}{N\kappa} \right)^{1/p}. $$

\noindent\textbf{Lemma \ref{F-property}.} \textit{Under Assumption \ref{coercivity}, given any $N>0$ and $\sigma>0$, (1) both $F_{N,\sigma}(\bm{\mu})$ and $\nabla F_{N,\sigma}(\bm{\mu})$ are well-defined and Lipschitz in $\mathbb{R}^d$; (2) The Lipschitz constant for $\nabla F_{N,\sigma}$ is $L=2d\sigma^{-2}e^{Nf(\bm{x}^*)}$, (3) $F_{N,\sigma}$ has at least one global maximum $\bm{\mu}^*\in\mathbb{R}^d$, and (4) $\|\hat{\nabla}F_{N,\sigma}(\bm{\mu})\|^2\leq G=d\sigma^2 e^{2Nf(\bm{x}^*)}$ for all $\bm{\mu}\in\mathbb{R}^d$.}

\begin{proof}
The proofs for Property (1), (2), and (4) are the same as those to their counterpart in \cite{{GS-PowerOpt}}, although in Assumption \ref{coercivity} we assume that the objective $f$ is defined over the whole space $\mathbb{R}^d$. Specifically, the proof for (1) is the same as that for \cite[Lemma 3.3]{GS-PowerOpt}. For (2), the proof that $F_{N,\sigma}(\bm{\mu})$ is Lipschitz is the same as that for Point 1 in the proof to \cite[Proposition 2.3]{GS-PowerOpt}. The proof for the Lipschitz property of $\nabla F_{N,\sigma}$ and $L=2d\sigma^{-2}e^{Nf(\bm{x})}$ is the same as the proof to \cite[Lemma 3.5]{GS-PowerOpt}. The proof to (4) is the same as that to \cite[Lemma 3.6]{GS-PowerOpt}.

Property (3) shares the same proof as \cite[Proposition 2.3]{GS-PowerOpt}, except the part of showing $\lim_{\|\bm{\mu}\|\rightarrow +\infty }F_{N,\sigma}=0$, which is given below. Specifically, for any $\epsilon>0$, let $M>0$ be such that whenever $\|\bm{x}\|>M$, $f(\bm{x})<(\ln\epsilon)/N$. Then, 
\begin{align*}
F_{N,\sigma}(\bm{\mu}) =&\frac{1}{(\sqrt{2\pi}\sigma)^d }\int_{\bm{x}\in\mathbb{R}^d} e^{Nf(\bm{x})}e^{-\frac{\|\bm{x-\mu}\|^2}{2\sigma^2}}d\bm{x}\\
=&\frac{1}{(\sqrt{2\pi}\sigma)^d }\int_{\|\bm{x}\|\leq M} e^{Nf(\bm{x})}e^{-\frac{\|\bm{x-\mu}\|^2}{2\sigma^2}}d\bm{x} \\
&+ \frac{1}{(\sqrt{2\pi}\sigma)^d }\int_{\|\bm{x}\|> M} e^{Nf(\bm{x})}e^{-\frac{\|\bm{x-\mu}\|^2}{2\sigma^2}}d\bm{x} \\
\leq &\frac{1}{(\sqrt{2\pi}\sigma)^d }\int_{\|\bm{x}\|\leq M} e^{Nf(\bm{x})}e^{-\frac{(\|\bm{\mu}\|-M)^2}{2\sigma^2}}d\bm{x} + \epsilon\\
\leq &2\epsilon,\quad \text{when }\|\bm{\mu}\|\text{ is larger than some } b_{\epsilon}>0.
\end{align*}
In sum, for any $\epsilon>0$, there exists $b_\epsilon>0$ such that whenever $\|\bm{\mu} \|>b_\epsilon$, $F(\bm{\mu})<2\epsilon$. Hence, $\lim_{\|\bm{\mu}\|\rightarrow +\infty }F_{N,\sigma}(\bm{\mu})=0$.
\end{proof}

\noindent\textbf{Theorem \ref{main-theorem}}: \textit{For any deterministic $\bm{\mu}_0\in\mathbb{R}^d$, $N,b,\sigma_0>0$, and positive integer $T$, let $\{(\bm{\mu}_t,\sigma_{t})\}_{t=1}^{T}$ be generated by the GS-PowerHP iterations (\ref{zo-slph}). Under Assumption \ref{coercivity} and \ref{f-Lipschitz}, we have
\begin{align*} \sum_{t=0}^{T-1}\alpha_t&\sigma_{t+1}^2\mathbb{E}[\|\nabla F_{N,\sigma_{t+1}}(\bm{\mu}_t) \|^2]\leq e^{Nf(\bm{x}^*)} -F_{N,\sigma_0}(\bm{\mu}_0)\\
&+ 2d^2e^{3N f(\bm{x}^*)} \sum_{t=0}^{T-1}\alpha_t^2 + \sigma_0C_{N,d,f}(1-\beta)\sum_{t=0}^{T-1}\beta^t,  
\end{align*}
where $C_{N,d,f}=N e^{Nf(\bm{x}^*)} L_f \sqrt{d}$.}
\begin{proof}
From the Gradient Mean Value Theorem, there exists $\bm{\nu}_t\in\mathbb{R}^d$ on the line segment joining $\bm{\mu}_{t+1}$ and $\bm{\mu}_{t}$ such that
\begin{equation*}
\begin{split}
F_{N,\sigma_{t+1}}(\bm{\mu}_{t+1})=&F_{N,\sigma_{t+1}}(\bm{\mu}_t) + (\nabla F_{N,\sigma_{t+1}}(\bm{\nu}_t))' (\bm{\mu}_{t+1}-\bm{\mu}_t)\\
=&F_{N,\sigma_{t+1}}(\bm{\mu}_t) + (\nabla F_{N,\sigma_{t+1}}(\bm{\mu}_t))' (\bm{\mu}_{t+1}-\bm{\mu}_t)\\
 &+ (\nabla F_{N,\sigma_{t+1}}(\bm{\nu}_t)-\nabla F_{N,\sigma_{t+1}}(\bm{\mu}_t))'(\bm{\mu}_{t+1}-\bm{\mu}_t)\\
=&F_{N,\sigma_{t+1}}(\bm{\mu}_t) + \alpha_t (\nabla F_{N,\sigma_{t+1}}(\bm{\mu}_t))' \hat{\nabla} F_{N,\sigma_{t+1}}(\bm{\mu}_t) \\
&+ (\nabla F_{N,\sigma_{t+1}}(\bm{\nu}_t)-\nabla F_{N,\sigma_{t+1}}(\bm{\mu}_t))'(\bm{\mu}_{t+1}-\bm{\mu}_t)\\
\geq & F_{N,\sigma_{t+1}}(\bm{\mu}_t) + \alpha_t (\nabla F_{N,\sigma_{t+1}}(\bm{\mu}_t))' (\hat{\nabla} F_{N,\sigma_{t+1}}(\bm{\mu}_t))\\
&-L\|\bm{\nu}_t-\bm{\mu}_t\|\cdot\|\bm{\mu}_{t+1}-\bm{\mu}_t \|\\
\geq & F_{N,\sigma_{t+1}}(\bm{\mu}_t) + \alpha_t (\nabla F_{N,\sigma_{t+1}}(\bm{\mu}_t))' (\hat{\nabla} F_{N,\sigma_{t+1}}(\bm{\mu}_t))\\
&-L\|\bm{\mu}_{t+1}-\bm{\mu}_t \|^2,\\
= & F_{N,\sigma_{t+1}}(\bm{\mu}_t) + \alpha_t (\nabla F_{N,\sigma_{t+1}}(\bm{\mu}_t))' (\hat{\nabla} F_{N,\sigma_{t+1}}(\bm{\mu}_t))\\
&-\alpha_t^2L\|\hat{\nabla} F_{N,\sigma_{t+1}}(\bm{\mu}_t)\|^2,
\end{split}
\end{equation*}
where $'$ denotes vector transpose, the third equality is from (\ref{zo-slph}), and the first inequality is because of the Cauchy-Schwarz inequality and (2) in Lemma \ref{F-property}.
Taking the expectation of the left-end and right-end gives
\begin{equation*}
\begin{split}
\mathbb{E}[&F_{N,\sigma_{t+1}}(\bm{\mu}_{t+1})] \geq \mathbb{E}[F_{N,\sigma_{t+1}}(\bm{\mu}_t)]+\alpha_t\sigma_{t+1}^2\mathbb{E}[ \|\nabla F_{N,\sigma_{t+1}}(\bm{\mu}_t)) \|^2]\\
& - \alpha_t^2 L \mathbb{E}[\| \hat{\nabla} F_{N,\sigma_{t+1}}(\bm{\mu}_t) \|^2]\\
&\geq^{\text{Lemma \ref{F-property}(4)}} \mathbb{E}[F_{N,\sigma_{t+1}}(\bm{\mu}_t)]+\alpha_t\sigma_{t+1}^2\mathbb{E}[ \|\nabla F_{N,\sigma_{t+1}}(\bm{\mu}_t)) \|^2] - \alpha_t^2 L G,  \\
&=^{\text{Lemma } \ref{F-property}}\mathbb{E}[F_{N,\sigma_{t+1}}(\bm{\mu}_t)]+\alpha_t\sigma_{t+1}^2\mathbb{E}[ \|\nabla F_{N,\sigma_{t+1}}(\bm{\mu}_t)) \|^2] - \alpha_t^2 2d^2e^{3Nf(\bm{x}^*)},
\end{split}
\end{equation*}
where the first inequality holds because 
\begin{align*}
\mathbb{E}[(\nabla F_{N,\sigma_{t+1}}(\bm{\mu}_t))' (\hat{\nabla} F_{N,\sigma_{t+1}}(\bm{\mu}_t))] &=\mathbb{E}\left[ \mathbb{E}[(\nabla F_{N,\sigma_{t+1}}(\bm{\mu}_t))' (\hat{\nabla} F_{N,\sigma_{t+1}}(\bm{\mu}_t)) | \bm{\mu}_t ] \right] \\
&=^{\text{by (\ref{nablaF-est})}}\sigma^2_{t+1}\mathbb{E}[ \|\nabla F_{N,\sigma_{t+1}}(\bm{\mu}_t)) \|^2].
\end{align*}
In sum, we arrived at an inequality which will be used later:
\begin{equation}
\begin{split}
\label{change-mu}
&\mathbb{E}[F_{N,\sigma_{t+1}}(\bm{\mu}_{t+1})] \geq  \mathbb{E}[F_{N,\sigma_{t+1}}(\bm{\mu}_t)]+\alpha_t\sigma_{t+1}^2\mathbb{E}[ \|\nabla F_{N,\sigma_{t+1}}(\bm{\mu}_t)) \|^2] - \alpha_t^2 2d^2e^{3Nf(\bm{x}^*)}.
\end{split}
\end{equation}
Next, we derive the relation between $F_{N,\sigma_{t+1}}(\bm{\mu}_{t})$ and $F_{N,\sigma_{t}}(\bm{\mu}_t)$. Let $f_N(\bm{x})=e^{Nf(\bm{x})}$. For any $\bm{\mu}\in\mathbb{R}^d$,
\begin{equation*}
\begin{split}
&|F_{N,\sigma_{t+1}}(\bm{\mu}) - F_{N,\sigma_{t}}(\bm{\mu})| \\
=& \left|\mathbb{E}_{\bm{\epsilon}\sim\mathcal{N}(\bm{0},I_d)}[f_N(\bm{\mu}+\sigma_{t+1}\bm{\epsilon})]\right.\left.-\mathbb{E}_{\bm{\epsilon}\sim\mathcal{N}(\bm{0},I_d)}[f_N(\bm{\mu}+\sigma_{t}\bm{\epsilon})] \right|\\
=&\left|(2\pi)^{-d/2}\int_{\bm{\epsilon}\in\mathbb{R}^d} (f_N(\bm{\mu}+\sigma_{t+1}\bm{\epsilon})-f_N(\bm{\mu}+\sigma_{t}\bm{\epsilon}))e^{-\|\bm{\epsilon}\|^2/2}d\bm{\epsilon} \right|\\
\leq &(2\pi)^{-d/2}\int_{\bm{\epsilon}\in\mathbb{R}^d} \left|f_N(\bm{\mu}+\sigma_{t+1}\bm{\epsilon})-f_N(\bm{\mu}+\sigma_{t}\bm{\epsilon})\right| e^{-\|\bm{\epsilon}\|^2/2}d\bm{\epsilon} \\
= &(2\pi)^{-d/2}\int_{\bm{\epsilon}\in\mathbb{R}^d} \left| e^{Nf(\bm{\mu}+\sigma_{t+1}\bm{\epsilon})} -e^{Nf(\bm{\mu}+\sigma_{t}\bm{\epsilon})} \right| e^{-\|\bm{\epsilon}\|^2/2}d\bm{\epsilon} \\
\leq &(2\pi)^{-d/2}\int_{\bm{\epsilon}\in\mathbb{R}^d}Ne^{Nf(\bm{x}^*)}|f(\bm{\mu}+\sigma_{t+1}\bm{\epsilon})-f(\bm{\mu}+\sigma_t\bm{\epsilon})| e^{-\|\bm{\epsilon}\|^2/2}d\bm{\epsilon} \\
\leq &\frac{1}{(2\pi)^{d/2}}\int_{\bm{\epsilon}\in\mathbb{R}^d}Ne^{Nf(\bm{x}^*)} L_f|\sigma_{t+1}-\sigma_t|\|\bm{\epsilon}\|  e^{-\|\bm{\epsilon}\|^2/2}d\bm{\epsilon} \\
\leq &N e^{Nf(\bm{x}^*)} L_f |\sigma_{t+1}-\sigma_t| \mathbb{E}\|\bm{\epsilon}\| \\
\leq &N e^{Nf(\bm{x}^*)} L_f |\sigma_{t+1}-\sigma_t| \sqrt{\mathbb{E}[\|\bm{\epsilon}\|^2]} \\
\leq &N e^{Nf(\bm{x}^*)} L_f |\sigma_{t+1}-\sigma_t| \sqrt{d} \\
= &N e^{Nf(\bm{x}^*)} L_f \sigma_0 (1-\beta)\beta^t \sqrt{d} \\
= & \sigma_0C_{N,d,f}(1-\beta) \beta^t, 
\end{split}
\end{equation*}
where $C_{N,d,f}=N e^{Nf(\bm{x}^*)} L_f \sqrt{d}$, the fifth relation is because of mean value theorem, the sixth relation is because of Assumption \ref{f-Lipschitz}, $S_d$ denotes the unit sphere surface in $\mathbb{R}^d$, and $\Gamma$ denotes the Gamma function. In sum, we have
\begin{equation*}
|F_{N,\sigma_{t+1}}(\bm{\mu}) - F_{N,\sigma_{t}}(\bm{\mu})| \leq \sigma_0C_{N,d,f}(1-\beta) \beta^t.
\end{equation*}
This result further implies
$$ F_{N,\sigma_{t+1}}(\bm{\mu}_{t}) \geq F_{N,\sigma_{t}}(\bm{\mu}_t)-\sigma_0C_{N,d,f}(1-\beta)\beta^t.$$
Plugging this result to (\ref{change-mu}) gives
\begin{align*}
\mathbb{E}[F_{N,\sigma_{t+1}}(\bm{\mu}_{t+1})] \geq & \mathbb{E}[F_{N,\sigma_{t}}(\bm{\mu}_t)]+\alpha_t\sigma_{t+1}^2\mathbb{E}[ \|\nabla F_{N,\sigma_{t+1}}(\bm{\mu}_t)) \|^2] \\
&- \alpha_t^2 2d^2e^{3Nf(\bm{x}^*)} - \sigma_0C_{N,d,f}(1-\beta)\beta^t. 
\end{align*}
Summing both sides over $t\in\{0,1,...,T-1\}$ gives
\begin{align*} 
\mathbb{E}[F_{N,\sigma_T}(\bm{\mu}_T)] \geq &\mathbb{E}[F_{N,\sigma_0}(\bm{\mu}_0)] + \sum_{t=0}^{T-1}\alpha_t\sigma_{t+1}^2\mathbb{E}[\|\nabla F_{N,\sigma_{t+1}}(\bm{\mu}_t) \|^2]\\
 &- 2d^2e^{3Nf(\bm{x}^*)}\sum_{t=0}^{T-1}\alpha_t^2 - \sigma_0C_{N,d,f}(1-\beta)\sum_{t=0}^{T-1}\beta^t. 
 \end{align*}
After re-organizing the terms, we arrive at
\begin{align*}  
\sum_{t=0}^{T-1}\alpha_t\sigma_{t+1}^2\mathbb{E}[&\|\nabla F_{N,\sigma_{t+1}}(\bm{\mu}_t) \|^2]
\leq  \mathbb{E}[F_{N,\sigma_T}(\bm{\mu}_T)] - \mathbb{E}[F_{N,\sigma_0}(\bm{\mu}_0)] \\
 +&  2d^2e^{3Nf(\bm{x}^*)}\sum_{t=0}^{T-1}\alpha_t^2 + \sigma_0 C_{N,d,f}(1-\beta)\sum_{t=0}^{T-1}\beta^t\\
\leq& f_N(\bm{x}^*) - F_{N,\sigma_0}(\bm{\mu}_0)+2d^2e^{3Nf(\bm{x}^*)}\sum_{t=0}^{T-1}\alpha_t^2 + \sigma_0 C_{N,d,f}(1-\beta)\sum_{t=0}^{T-1}\beta^t
\end{align*}
This finishes the proof for Theorem \ref{main-theorem}.
\end{proof}

\noindent\textbf{Corollary \ref{convergence-rate}}. \textit{Assume the conditions in Theorem \ref{main-theorem} and Assumption \ref{lr}. For any positive integer $T$, define 
$$\nu_T:= \min_{\tau\in\{0,1,...,T\}}\mathbb{E}[\|\nabla F_{N,\sigma_{\tau+1}}(\bm{\mu}_\tau)\|^2],$$ where $\{\bm{\mu}_{\tau}\}$ are generated from GS-PowerHP (\ref{zo-slph}) with $\beta\in [0,1)$. Then, we have $\nu_T=\mathcal{O}_N\left( (\sigma_0\beta^T+b)^{-2}(\sum_{t=0}^{T-1} \alpha_t )^{-1}\right)$, and the dependence on $N$ can be removed if $f(\bm{x}^*)<0$.} 

\begin{proof}
\begin{equation*}
\begin{split}
\label{C0}
&\sum_{t=0}^{T-1} \alpha_t \sigma^2_{T} \nu_{T} \leq \sum_{t=0}^{T-1} \alpha_t \sigma^2_{t+1} \nu_t \\
\leq&^{by (\ref{nu-def})} \sum_{t=0}^{T-1} \alpha_t  \sigma^2_{t+1}\mathbb{E}[\|\nabla F_{N,\sigma_{t+1}}(\bm{\mu}_t)\|^2] \\
 \leq&^{\text{Theorem \ref{main-theorem}}} f_N(\bm{x}^*) - F_{N,\sigma_0}(\bm{\mu}_{0}) + 2d^2f_N^3(\bm{x}^*) \sum_{t=0}^{\infty}\alpha_t^2 + \sigma_0 C_{N,d,f},\\
\end{split}
\end{equation*}
where $C_{N,d,f}$ is as defined in Theorem \ref{main-theorem}.
Therefore, we have
\begin{equation}
\begin{split}
\label{nu-T}
\nu_{T} &\leq  \frac{f_N(\bm{x}^*) - F_{N,\sigma_0}(\bm{\mu}_{0}) + 2d^2f_N^3(\bm{x}^*) \sum_{t=0}^{\infty}\alpha_t^2 + \sigma_0 C_{N,d,f}}{\sigma^2_T\sum_{t=0}^{T-1} \alpha_t} \\
&= \mathcal{O}_N\left( (\sigma_0\beta^T+b)^{-2}\left(\sum_{t=0}^{T-1} \alpha_t \right)^{-1}\right).
\end{split}
\end{equation}
If $f(\bm{x}^*)< 0$, then
\begin{itemize}
\item $f_N(\bm{x}^*)-F_{N,\sigma_0}(\bm{\mu}_0)=e^{Nf(\bm{x}^*)} - \mathbb{E}_{\bm{x}\sim\mathcal{N}(\bm{\mu}_0,\sigma^2_0 I_d)}[e^{Nf(\bm{x})}]\in [0, 1]$;
\item $f_N^3(\bm{x}^*)\leq 1$;
\item $\lim_{N\to \infty} C_{N,d,f}=\lim_{N\to \infty} Ne^{Nf(\bm{x}^*)L_f\sqrt{d}}=0$, which implies that $C_{N,d,f}$ has an upper bound that is independent from $N$.
\end{itemize}
Plugging these inequalities to the first inequality in (\ref{nu-T}) gives
$$ \nu_{T} = \mathcal{O}\left( (\sigma_0\beta^T+b)^{-2}\left(\sum_{t=0}^{T-1} \alpha_t \right)^{-1}\right). $$

\end{proof}

\noindent\textbf{Corollary \ref{iter-complexity}}. \textit{Assume the conditions in Theorem \ref{main-theorem} and Assumption \ref{lr}. Let $\alpha_t = (t+1)^{-(1/2+\gamma})$ for some $\gamma\in (0,1/2)$. For any $\varepsilon\in(0,1)$, after $T> (C_2C_1d^2\varepsilon^{-1})^{2/(1-2\gamma)}=O_N((d^2\varepsilon^{{-1}}b^{-2})^{2/(1-2\gamma)})$ times of parameter updating, we have
$$\min_{t\in\{0,1,2,...,T\}}\mathbb{E}[ \| \nabla F_{N,\sigma_{t+1}}(\bm{\mu}_t) \|^2 ]<\varepsilon,$$
where $C_0:=f_N(\bm{x}^*) - F_{N,\sigma_0}(\bm{\mu}_{0}) + 2f_N^3(\bm{x}^*) \sum_{t=1}^\infty t^{-(1+2\gamma)} + \sqrt{2}\sigma_0Nf_N(\bm{x}^*)L_f$, $C_1 = b^{-2}C_0(1-2\gamma)$, and $C_2=\max\{1,1/C_1\}$. Furthermore, the dependence of the $O_N$-factor on $N$ can be removed under the additional assumption of $f(\bm{x})< 0$ for all $\bm{x}\in\mathbb{R}^d$.}
\begin{proof}
For any non-negative integer $t$, define
\begin{equation} 
\label{nu-def}
\nu_t:= \min_{\tau\in\{0,1,...,t\}}\mathbb{E}[\|\nabla F_{N,\sigma_{\tau+1}}(\bm{\mu}_\tau)\|^2].
\end{equation}
Then, 
\begin{equation}
\begin{split}
\label{C0}
&\sum_{t=0}^{T-1} \alpha_t \sigma^2_{T} \nu_{T} \leq \sum_{t=0}^{T-1} \alpha_t \sigma^2_{t+1} \nu_t \\
\leq&^{by (\ref{nu-def})} \sum_{t=0}^{T-1} \alpha_t  \sigma^2_{t+1}\mathbb{E}[\|\nabla F_{N,\sigma_{t+1}}(\bm{\mu}_t)\|^2] \\
 \leq&^{\text{Theorem \ref{main-theorem}}} f_N(\bm{x}^*) - F_{N,\sigma_0}(\bm{\mu}_{0}) + 2d^2f_N^3(\bm{x}^*) \sum_{t=0}^{T-1}\alpha_t^2 \\
&+ \sigma_0 C_{N,d,f}(1-\beta)\sum_{t=0}^{T-1}\beta^t,\\
 \leq &f_N(\bm{x}^*) - F_{N,\sigma_0}(\bm{\mu}_{0}) + 2d^2f_N^3(\bm{x}^*) \sum_{t=1}^\infty t^{-(1+2\gamma)} +\sigma_0C_{N,d,f}\\
\leq &\left(f_N(\bm{x}^*) - F_{N,\sigma_0}(\bm{\mu}_{0}) + 2f_N^3(\bm{x}^*) \sum_{t=1}^\infty t^{-(1+2\gamma)}\right.\\
&\left. + \sqrt{2}\sigma_0Nf_N(\bm{x}^*)L_f \right)d^2,\\
 = &C_{0,N} d^2,
\end{split}
\end{equation}
where $C_{0,N}:=f_N(\bm{x}^*) - F(\bm{\mu}_{0}) + 2f_N^3(\bm{x}^*) \sum_{t=1}^\infty t^{-(1+2\gamma)} + \sqrt{2}\sigma_0Nf_N(\bm{x}^*)L_f  <\infty$, the second-to-last relation is because $d>1$ and $f_N(\bm{x}^*)\geq F_{N,\sigma_0}(\bm{\mu}_0)$, $C_{0,N}=\sqrt{2}\frac{\Gamma((d+1)/2)}{\Gamma(d/2)}Nf_N(\bm{x}^*) L_f$, and $\frac{\Gamma((d+1)/2)}{\Gamma(d/2)}\leq d$ for any $d\geq 1$.
Therefore, we have
$$\sum_{t=0}^{T-1} \alpha_t \sigma^2_T \nu_{T}\leq  C_{0,N} d^2. $$
Since $\alpha_t=(t+1)^{-(1/2+\gamma)}$, dividing both sides of the above inequality by $\sum_{t=0}^{T-1}\alpha_t\sigma^2_T$ gives
\begin{equation*}
\begin{split}
\nu_T&\leq\frac{C_{0,N}d^2}{\sigma^2_T\sum_{t=1}^{T} t^{-(1/2+\gamma)}}\\
&<\frac{C_{0,N}d^2}{b^2\int_{1}^Tt^{-(1/2+\gamma)}dt},\quad \text{recall }\sigma_T = \sigma_0\beta^T+b,\\
&<\frac{C_{0,N}d^2}{b^2\int_{1}^Tt^{-(1/2+\gamma)}dt}\\
&=\frac{C_{0,N}d^2}{b^2(T^{\frac{1}{2}-\gamma}-1)/(\frac{1}{2}-\gamma)}\\
&<\frac{C_{0,N}d^2}{b^2(T^{\frac{1}{2}-\gamma}/2)/(\frac{1}{2}-\gamma)},\quad\text{when }T>2^{2/(1-2\gamma)},\\
&= C_{1,N} \frac{d^2}{T^{\frac{1}{2}-\gamma}},\quad \text{where } C_{1,N}:=\frac{C_{0,N}}{b^2(1/2)/(\frac{1}{2}-\gamma)}.
\end{split}
\end{equation*}
In sum, we have
\begin{equation}
\label{nu-bound}
\nu_T \leq C_{1,N} \frac{d^2}{T^{\frac{1}{2}-\gamma}},\quad \text{ whenever }T>2^{2/(1-2\gamma)}.
\end{equation}
Define $C_{2,N}:=\max\{1, 2/C_{1,N}\}$. Given any $\varepsilon\in(0,1)$, whenever $$T> (C_{2,N}C_{1,N}d^2\varepsilon^{-1})^{2/(1-2\gamma)}=O((d^2b^{-2} \varepsilon^{-1})^{2/(1-2\gamma)}),$$ 
we have
\begin{align*} 
T&>(C_{2,N}C_{1,N}d^2\varepsilon^{-1})^{2/(1-2\gamma)}\\
&>(C_{2,N}C_{1,N})^{2/(1-2\gamma)}\\
&\geq 2^{2/(1-2\gamma)}, 
\end{align*}
and
\begin{align*} \nu_T &\leq^{\text{from }(\ref{nu-bound})}  C_{1,N} \frac{d^2}{T^{\frac{1}{2}-\gamma}} \\
&< C_{1,N} \frac{d^2}{(C_{2,N}C_{1,N}d^2\varepsilon^{-1})^{\frac{2}{1-2\gamma}(\frac{1}{2}-\gamma)}} \\
&=\frac{\varepsilon}{C_{2,N}}\leq \varepsilon.
\end{align*}
This finishes the proof for the first result.

Under the additional assumption that $f(\bm{x})< 0$ for all $\bm{x}\in\mathbb{R}^d$, we have $f_N(\bm{x}^*)=e^{Nf(\bm{x}^*)}< 1$ and $Nf_N(\bm{x}^*)<A^*$ for some $A^*>0$ that is independent from $N$. Hence, $C_{0,N}\leq C_0$, where 
$$C_{0}:=1 - F(\bm{\mu}_{0}) + 2(\bm{x}^*) \sum_{t=1}^\infty t^{-(1+2\gamma)} + \sqrt{2}\sigma_0A^*L_f.$$ 
The above proof still holds if we replace $C_{0,N}$ with $C_0$, $C_{1,N}$ with $C_1:=\frac{C_{0}}{b^2(1/2)/(\frac{1}{2}-\gamma)}$, and $C_{2,N}$ with $C_2:=\max\{1,1/C_{1}\}$. Then, the $O_N$ factor's dependence on $N$ is removed. This finishes the proof for Cororllary \ref{iter-complexity}.
\end{proof}

\noindent\textbf{Theorem} \ref{thm2} \textit{
Suppose Assumption \ref{coercivity} holds. Given any positive numbers $N$, $\sigma$, $\delta$, and $M$ such that $\delta<M$, there exists $N_{\delta,\sigma,M}>0$ such that whenever $N>N_{\delta,\sigma,M}$ the following statement holds true: if $\bm{\mu}\in\mathbb{R}^d$ and $\|\bm{\mu}\|<M$, then for any $i\in\{1,2,...,d\}$}
\begin{enumerate}
 \item $\mu_i>x^*+\delta \Rightarrow \frac{\partial F_{N,\sigma}(\bm{\mu})}{\partial \mu_i}<0$; \textit{and}
 \item $\mu_i<x^*-\delta \Rightarrow \frac{\partial F_{N,\sigma}(\bm{\mu})}{\partial \mu_i}>0$,
\end{enumerate}
\textit{where $\mu_i$ and $x^*_i$ denotes the $i^{th}$ entry of $\bm{\mu}$ and $\bm{x}^*$, respectively.}
\begin{proof}
Despite of some small difference between the assumptions for this theorem and those for \cite[Theorem 2.1]{GS-PowerOpt}, the two theorems share the same proof, which can be found in their appendix.
\end{proof}

\noindent\textbf{Corollary} \ref{summary-cor}
\textit{Suppose Assumption \ref{coercivity}, \ref{f-Lipschitz}, \ref{lr}, and \ref{mu-bounded} hold. For any $\delta>0$, $b>0$, and $\bm{\mu}_0\in\mathbb{R}^d$, let $N_{\delta,b,\bm{\mu}_0}$ and $M_{(\delta,b,\bm{\mu}_0)}$ be the threshold and $\|\bm{\mu}\|$-bound in Assumption \ref{mu-bounded}, respectively. Let $N_{\delta,b,M_{(\delta,b,\bm{\mu}_0)}}$ be the threshold stated in Theorem \ref{thm2}. Then, whenever $N>\max\{N_{\delta,b,\bm{\mu}_0},N_{\delta,b,M_{(\delta,b,\bm{\mu}_0)}}\}$, $\{\bm{\mu}_t\}$ generated by (\ref{zo-slph}) converges in expectation to $\mathcal{S}_{\bm{x}^*,\delta}$ with the iteration complexity of $O((d^2\varepsilon^{-1}b^{-2})^{2/(1-2\gamma)})$, given the learning rate of $\alpha_t = (t+1)^{-(1/2+\gamma)}$ for some $\gamma\in (0,1/2)$.}
\begin{proof}
Corollary \ref{iter-complexity} implies that $\bm{\mu}_t$ converges in expectation to a stationary point $\bm{\mu}^*$ of $F_{N,b}(\bm{\mu})$. Since $\|\bm{\mu}_t\|\leq M_{(\delta,b,\bm{\mu}_0)}$, $\|\bm{\mu}^*\|\leq M_{(\delta,b,\bm{\mu}_0)}$. Also, Theorem \ref{thm2} implies that, all the stationary points of $F_{N,b}$ within the region of $\|\bm{\mu}\|\leq M_{(\delta,b,\bm{\mu}_0)}$ lie in $\mathcal{S}_{\bm{x}^*;\delta}$. This completes the proof.
\end{proof}

\subsection{Performance of Square Attack as a Benchmark for Image Adversarial Attacks}
\label{appendix:square-attack}
This appendix reports the performance of Square Attack~\cite{Andriushchenko2020} as a specialized black-box adversarial-attack benchmark. Unlike the general-purpose zeroth-order optimization methods considered in the main experiments, Square Attack is specifically designed for image adversarial attacks and uses a different attack objective and stopping rule. Therefore, we include it as a strong reference baseline, while the results of GS-PowerHP are copied from Tables~\ref{tab:per-image-attacks-mnist-cifar} and~\ref{imagenet-attack} for comparison.

\begin{table}[h!]
\caption{Comparison between GS-PowerHP and Square Attack on the same 100 selected images from MNIST, CIFAR-10, and ImageNet used in Tables~\ref{tab:per-image-attacks-mnist-cifar} and~\ref{imagenet-attack}. The results of GS-PowerHP are copied from the corresponding main-text tables. Square Attack updates one candidate perturbation at each iteration, whereas the smoothing-based methods use a population size of $K=10$.}
\label{tab:square-attack-appendix}
\centering
\begin{tabular}{ >{\centering\arraybackslash} p{2.0cm} | >{\centering\arraybackslash}p{3cm} >{\centering\arraybackslash}p{1.0cm}>{\centering\arraybackslash}p{1.6cm} >{\centering\arraybackslash}p{2.0cm}>{\centering\arraybackslash}p{2.0cm} }  
\midrule
Imageset & Algorithm& SR & $\bar{R}^2$ & $\overline{\|\bm{\mu}^*\|_2}$ & $\bar{T}$ \\
\toprule
\multirow{2}{*}{MNIST} 
&GS-PowerHP & $100\%$ & $92\%$ & $4.68(0.97)$ & $1247(314)$ \\
&Square Attack & $100\%$  & $64\%$ & $9.76(0.31)$ & $887(1126)$ \\
\hline
\multirow{2}{*}{CIFAR-10} 
&GS-PowerHP & $100\%$ & $99\%$ & $1.67(0.36)$ & $597(212)$ \\
&Square Attack & $100\%$  & $78\%$ & $9.85(0.17)$ & $471(486)$ \\
\hline
\multirow{2}{*}{ImageNet} 
&GS-PowerHP & $78\%$      &  $95\%(7\%)$     & $35.8(6.5)$          & $1389(542)$\\ 
& Square Attack &$100\%$ &$90\%(9\%)$ &$49.94(0.11)$ &$2209(2001)$ \\
\bottomrule
\end{tabular}
\end{table}
Table~\ref{tab:square-attack-appendix} shows that Square Attack achieves a $100\%$ success rate on all three datasets, confirming its strength as a specialized image-adversarial-attack method. However, GS-PowerHP produces substantially smaller perturbations and higher image similarity on all three datasets; in particular, on MNIST and CIFAR-10, it matches the $100\%$ success rate of Square Attack while using much smaller perturbations, and on ImageNet it achieves a lower success rate but still produces smaller and more visually similar adversarial examples. This advantage is further illustrated by the adversarial examples reported in Tables~\ref{Adversarial-Image-MNIST},~\ref{Adversarial-Image-CIFAR}, and~\ref{Adversarial-Image-ImageNet}.

\subsection{Untarged Image Adversarial Attacks on An Adversarially Robust Classifier}
\label{sec-untargeted-IAA}
To test the robustness of our method across different tasks, under similar settings in Subsection \ref{Experiment-IAA}, we perform \textit{untargeted} image adversarial attacks on a ResNet18-based robust classifier for ImageNet (i.e., Salman2020Do\_R18 \cite{salman2020adversarially} implemented with the Python package of RobustBench \cite{croce2020robustbench}). We set the fitness function (maximization objective) for this task as
\begin{equation}
 f(\bm{x}):= \max( \max_{i\neq\mathcal{T}} \mathcal{G}(\bm{a+y})_i -  \mathcal{G}(\bm{a+y})_{\mathcal{T}}, \; \kappa) - \lambda \|\bm{y}\|,
\end{equation}
where $\bm{y}=\tanh{\bm{x}}$, $\lambda>0$ denotes the regularizing coefficient, and $\mathcal{T}$ denotes the label predicted by the classifier for image $\bm{a}$. For this task, we set $\kappa=0.001$. The selected hyper-parameter values and the corresponding candidate sets are reported in Appendix \ref{appendix-hyperparameters}.
\newpage
\begin{table}[h!]
\caption{Per-image Untargeted Attack on 100 Images from ImageNet. For each image attack, we set $\bm{\mu}_0=\bm{0}$ and $T_{total}=1,000$. The hyper-parameters for GS-PowerHP are selected as $N=8,\sigma_0=0.1$, $\alpha=2.0$, and $\beta=0.995$. The hyper-parameter values of other algorithms are reported in Table \ref{imagenet-hyper-parameters-untargeted}. For all algorithms except Square A (square attack), we perform $1,000$ iterations where 10 perturbations are sampled in each iteration. We omit the Square Attack results from the table as it is not a smoothing-based algorithm. However, it achieves a $100\%$ Success Rate (SR) with an average perturbation norm of 27.23. To ensure a fair comparison, the budget was set to 10,000 iterations; this is equivalent to 1,000 iterations of other algorithms, as Square Attack generates only a single perturbation per iteration. }
\label{attack-robust-imagenet}
\centering%
\begin{tabular}{ >{\centering\arraybackslash} p{3cm} | >{\centering\arraybackslash}p{0.9cm} >{\centering\arraybackslash}p{1.90cm}>{\centering\arraybackslash}p{1.9cm} >{\centering\arraybackslash}p{1.5cm} }  
\midrule
Algorithm & SR & $\bar{R}^2$ & $\overline{\|\bm{\mu}^*\|}$ & $\bar{T}$ \\
\toprule
\textbf{Our Algo.} & $\mathbf{100}\%$      &  $89\%(10\%)$     & $\mathbf{25.8}(11.6)$          & $208(188)$\\ 
GS-PowerOpt     & $82\%$  &$88\%(12\%)$     & $24.3(10.3)$    & $270(243)$\\
ZOSLGHr     & $94\%$  &$86\%(16\%)$     & $30.1(12.8)$    & $54(143)$\\
ZOAdaMM   & $97\%$  &$83\%(18\%)$     & $32.3(12.7)$    & $107(266)$\\
ZOSGD  & $53\%$  &$91\%(13)$     & $21.1(4.6)$    &  $985(100)$\\
ZOSLGHd     & $52\%$  &$95\%(5\%)$     & $17.0(2.6)$    & $243(276)$\\
STD-Htp  &$63\%$   &$90\%(12\%)$  & $22.1(9.2)$    &$274(266)$\\
\bottomrule
\end{tabular}
\end{table}

\subsection{Hyper-parameters Selection in Each Experiment}
\label{appendix-hyperparameters}
\subsubsection{Hyper-parameter Values for Optimizing Ackley and Rastrigin} For each algorithm and each test objective (either Ackley or Rastrigin), for each algorithm and each candidate set of hyperparameters, we perform 10 trials to optimize the objective function. In each trial, each algorithm is performed for $T$ iterations. Let $\bm{\mu}^*_n$ denote the candidate solution in the $n^{th}$ trial with the largest $f$ value. Then, for each algorithm, we select the hyper-parameter set that produces the largest average fitness: $\sum_{n=1}^{10}f(\bm{\mu}_n^*)/10$. The hyper-parameter candidates and selected values are reported in Table \ref{Ackley-hyper-parameters} for Ackley and Table \ref{Rastrigin-hyper-parameters} for Rastrigin.

\begin{table}[h!]
\caption{Hyper-parameters for Maximizing Ackley (Table \ref{tab:combined-ackley-rastrigin}). The candidate set for (initial) smoothing parameters is $\mathcal{S}:=\{1.0, 1.5, 2.0\}$. The hyper-parameter symbols for each algorithm are the same as their source publications. For example, $t_1$ in ZO-SLGHd and ZO-SLGHr denotes the initial scaling parameter, $\mu$ in ZO-AdaMM is the scaling parameter, and $\alpha$ denotes a constant learning rate. The set of candidate values that lead to the highest $f$ value (averaged over 10 trials) are selected.}
\label{Ackley-hyper-parameters}
\centering%
\begin{tabular}{ p{2.5cm}  p{4.3cm} | p{8cm} }
\midrule
Algorithm & Selected Values & Candidates \\
\toprule
GS-PowerHP & $\alpha =.1$, $N=5$, $\sigma_0=1.0$, $b=0$, $\beta=.999$.   &  $N\in\{.05, 1, 3, 5\}$, $\alpha \in\{1.0, .1, .05\}$, $\sigma_0\in\mathcal{S}$.       \\
GS-PowerGS & $\alpha =.1$, $N=5$, $\sigma=1.0$.   &  $N\in\{.05, 1, 3, 5\}$, $\alpha \in\{1.0, .1, .05\}$, $\sigma\in\mathcal{S}$.       \\

ZO-SLGHd   & $\beta=1.0$, $\eta=.001$, $t_1=2.0$, $\gamma=.999$      &   $\beta\in\{1.0, .1, .05, .01, .001\}$, $\eta\in\{.1,.01,.001\}$, $t_1\in\mathcal{S}$, $\gamma\in\{.999,.995\}$.    \\
ZO-SLGHr   & $\beta=1.0$, $t_1=1.0$,  $\gamma=.995$.     &   $\beta\in\{1.0, .1, .05, .01, .001\}$, $\gamma\in\{.999,.995\}$, $t_1\in\mathcal{S}$. \\
ZO-AdaMM   &$\beta_{1}=.9$, $\beta_2=.1$, $\alpha=1.0$, $\mu=1.0$.        &   $\alpha\in\{1.0, .1, .05, .01, .001\}$, $\beta_1\in\{.5,.7,.9\}$, $\beta_2\in\{.1,.3,.5\}$, $\mu\in\mathcal{S}$   \\
ZO-SGD   & $\alpha=1.0$, $\mu=2.0$.      &   $\alpha\in\{1.0, .1, .05, .01, .001\}$, $\mu\in\mathcal{S}$.       \\
STD-Homotopy  & $\alpha=1.0$, $\gamma=.8$, $\sigma=2.0$, $T_{\bm{\mu}}=500$, $\tau=100$, $N_\sigma=10$.   &  $\gamma\in\{.2,.5,.8\}$, $\alpha\in\{1.0, .1, .05, .01, .001\}$, $\sigma\in\mathcal{S}$.    \\
CMA-ES & $\sigma_0=.5$, $c=1.0$. & $\sigma_0\in\{.1,.5,1.0\}$, $c\in\{.5,1.0,1.5\}$.\\
\bottomrule
\end{tabular}
\footnotetext{Total number of generations: 1000. $\bm{m}_1=[-.5,-.5]$. $\bm{m}_2=[.5,.5]$.}
\footnotetext{\textbf{Notation.} $\mathcal{P}_0$ denotes the initial population. For any real random vector $\bm{z}$ and any two real scalar $a<b$, $\bm{z} \sim$ uniform$[a, b]$ denotes that each entry of $\bm{z}$ is sampled uniformly from the interval $[a, b]$.}
\end{table}

\begin{table}[h!]
\caption{Hyper-parameters for Maximizing Rastrigin (Table \ref{tab:combined-ackley-rastrigin}). The candidate set for (initial) smoothing parameters is $\mathcal{S}:=\{1.0, 1.5, 2.0\}$. The hyper-parameter symbols for each algorithm are the same as their source publications. For example, $t_1$ in ZO-SLGHd and ZO-SLGHr denotes the initial scaling parameter, $\mu$ in ZO-AdaMM is the scaling parameter, and $\alpha$ denotes a constant learning rate. The set of candidate values that lead to the highest $f$ value (averaged over 10 trials) are selected.}
\label{Rastrigin-hyper-parameters}
\centering%
\begin{tabular}{ p{2.5cm}  p{4.3cm} | p{8cm} }
\midrule
Algorithm & Selected Values & Candidates \\
\toprule
GS-PowerHP & $\alpha =.05$, $N=.02$, $\sigma_0=1.5$, $b=0$, $\beta=.999$.   &  $N\in\{.02, .2, 2.0\}$, $\alpha \in\{1.0, .1, .05\}$, $\sigma_0\in\mathcal{S}$.       \\
GS-PowerGS & $\alpha =.05$, $N=.02$, $\sigma=2.0$.   &  $N\in\{.02, .2, 2.0\}$, $\alpha \in\{1.0, .1, .05\}$, $\sigma\in\mathcal{S}$.       \\

ZO-SLGHd   & $\beta=.001$, $\eta=.01$, $t_1=2.0$, $\gamma=.999$      &   $\beta\in\{1.0, .1, .05, .01, .001,.0001\}$, $\eta\in\{.1,.01,.001\}$, $t_1\in\mathcal{S}$, $\gamma\in\{.999,.995\}$.    \\
ZO-SLGHr   & $\beta=.001$, $t_1=1.0$,  $\gamma=.995$.     &   $\beta\in\{1.0, .1, .05, .01, .001,.0001\}$, $\gamma\in\{.999,.995\}$, $t_1\in\mathcal{S}$. \\
ZO-AdaMM   &$\beta_{1}=.5$, $\beta_2=.1$, $\alpha=.01$, $\mu=1.0$.        &   $\alpha\in\{1.0, .1, .05, .01, .001,.0001\}$, $\beta_1\in\{.5,.9\}$, $\beta_2\in\{.1,.5\}$, $\mu\in\mathcal{S}$   \\
ZO-SGD   & $\alpha=.0001$, $\mu=1.0$.      &   $\alpha\in\{1.0, .1, .05, .01, .001,.0001\}$, $\mu\in\mathcal{S}$.       \\
STD-Homotopy  & $\alpha=.05$, $\gamma=.5$, $\sigma=1.0$, $T_{\bm{\mu}}=500$, $\tau=100$, $N_\sigma=10$.   &  $\gamma\in\{.2,.5,.8\}$, $\alpha\in\{1.0, .1, .05, .01, .001,.0001\}$, $\sigma\in\mathcal{S}$.    \\
CMA-ES & $\sigma_0=1.0$, $c=1.0$. & $\sigma_0\in\{.1,.5,1.0\}$, $c\in\{.5,1.0,1.5\}$.\\
\bottomrule
\end{tabular}
\footnotetext{Total number of generations: 1000. $\bm{m}_1=[-.5,-.5]$. $\bm{m}_2=[.5,.5]$.}
\footnotetext{\textbf{Notation.} $\mathcal{P}_0$ denotes the initial population. For any real random vector $\bm{z}$ and any two real scalar $a<b$, $\bm{z} \sim$ uniform$[a, b]$ denotes that each entry of $\bm{z}$ is sampled uniformly from the interval $[a, b]$.}
\end{table}

\newpage
\subsubsection{Hyper-parameters for the Targeted Image Adversarial Attacks.}
For the image adversarial attack experiment on MNIST and CIFAR-10, we report the hyper-parameter candidates and selections for GS-PowerHP, GS-PowerOpt, and Square Attack in Table \ref{mnist-hyper-parameters} and \ref{cifar-hyper-parameters}, respectively. The hyper-parameters for other algorithms are chosen the same as those used in \cite[Table 8-9]{GS-PowerOpt} since the same experiments are performed there, except for the hyper-parameters of GS-PowerOpt, for which we found better values: $(N,\sigma,\alpha)=(0.1,0.05,0.1)$ for MNIST and $(N,\sigma,\alpha)=(0.15,0.03,0.06)$ for CIFAR-10. 

\begin{table}[h]
\caption{Hyper-parameters for targeted attacks against MNIST (Table \ref{tab:per-image-attacks-mnist-cifar} and \ref{tab:square-attack-appendix}). }
\label{mnist-hyper-parameters}
\centering%
\begin{tabular}{ p{3cm}  p{4cm} | p{7cm} }
\midrule
 & Selected Values & Candidates $(\bm{\mu}^*)$ \\
\toprule
GS-PowerHP & $\alpha =.065$, $N=.5$, $\sigma=.05$.   &  $N\in\{.05,.1,.5\}$, $\alpha\in\{.05, .055, .06, .065, .07\}$, $\sigma\in \{.05, .1 \} $\\
Square Attack & $\varepsilon = 10$, $p=0.5$ & $\varepsilon\in \{ 1,5,10,20,100\}$, $p\in\{0.1,0.5,0.9\}$\\
\bottomrule
\end{tabular}
\end{table}

\begin{table}[h]
\caption{Hyper-parameters for targeted attacks against CIFAR-10 (Table \ref{tab:per-image-attacks-mnist-cifar} and \ref{tab:square-attack-appendix}). }
\label{cifar-hyper-parameters}
\centering%
\begin{tabular}{ p{3cm}  p{4cm} | p{7cm} }
\midrule
 & Selected Values & Candidates $(\bm{\mu}^*)$ \\
\toprule
GS-PowerHP & $\alpha =.06$, $N=.15$, $\sigma_0=.03$, $\beta=.999$.   &  $N\in\{.1,.15,.2\}$, $\alpha\in\{.05, .06, .07\}$, $\sigma\in \{.03, .05, .1\}$, $\beta\in\{.995,.999\}$.\\
Square Attack & $\varepsilon = 10$, $p=0.5$ & $\varepsilon\in \{1,10,15,20,40\}$, $p\in\{0.05,0.1,0.5,0.9\}$.\\
\bottomrule
\end{tabular}
\end{table}
For the least-likely targeted image adversarial attacks on ImageNet, the hyper-parameters are reported in Table \ref{imagenet-hyper-parameters}. Note that $b$ for GS-PowerHP is set as 0 for all experiments.

\begin{table}[h]
\caption{Hyper-parameters for Targeted Attack on ImageNet (Table \ref{imagenet-attack} and \ref{tab:square-attack-appendix}). The candidate set for (initial) smoothing parameters is $\mathcal{S}:=\{0.1, .05\}$. The hyper-parameter symbols for each algorithm are the same as their source publications. For example, $t_1$ in ZO-SLGHd and ZO-SLGHr denotes the initial scaling parameter, $\mu$ in ZO-AdaMM is the scaling parameter, and $\alpha$ denotes a constant learning rate. The set of candidate values that lead to the highest success rate and minimal perturbation size (averaged over the 10 image attackes) are selected.}
\label{imagenet-hyper-parameters}
\centering%
\begin{tabular}{ p{3cm}  p{5cm} | p{6.9cm} }
\midrule
Algorithm & Selected Values & Candidates $(\bm{\mu}^*)$ \\
\toprule
GS-PowerHP & $\alpha =1.0$, $N=8$, $\sigma=.1$.   &  $N\in\{5, 8\}$, $\alpha_t\in\{1.0, 1.5\}$, $\sigma_0\in\mathcal{S}$.       \\
GS-PowerGS & $\alpha =1.5$, $N=5$, $\sigma=.1$.   &  $N\in\{5, 8\}$, $\alpha_t\in\{1.0, 1.5\}$, $\sigma\in\mathcal{S}$.       \\

ZO-SLGHd   & $\beta=1/3072$, $\eta=10^{-5}$, $t_1=.1$, $\gamma=.995$      &   $\beta\in\{.1, 1/3072, .01/3072\}$, $\eta\in\{.0001/3072,10^{-5}\}$, $t_1\in\mathcal{S}$, $\gamma\in\{.999,.995\}$.    \\
ZO-SLGHr   & $\beta=1/3072$, $t_1=.1$,  $\gamma=.999$.     &   $\beta\in\{.1, 1/3072, .01/3072\}$, $\gamma\in\{.999,.995\}$, $t_1\in\mathcal{S}$. \\
ZO-AdaMM   &$\beta_{1}=.9$, $\beta_2=.1$, $\alpha=.1$, $\mu=.05$.        &   $\alpha\in\{0.5/3072, .1, 1.0\}$, $\beta_1\in\{.5,.9\}$, $\beta_2\in\{.1,.3\}$, $\mu\in\mathcal{S}$   \\
ZO-SGD   & $\alpha=0.01$, $\mu=.1$.      &   $\alpha\in\{.01, 1/3072, .01/3072\}$, $\mu\in\mathcal{S}$.       \\
Square Attack & $\varepsilon = 250$, $p=0.001$ & $\varepsilon\in \{ 20,60,100,150,200,250,300\}$, $p\in\{0.0001, 0.001, 0.005,0.01 \}$\\
STD-Homotopy  & $\alpha=.5$, $\gamma=.8$, $\sigma=.1$, $T_{\bm{\mu}}=300$, $\tau=100$, $N_\sigma=10$.   &  $\gamma\in\{.5,.8\}$, $\alpha\in\{.5,.1\}$.     \\
\bottomrule
\end{tabular}
\footnotetext{Total number of generations: 1000. $\bm{m}_1=[-.5,-.5]$. $\bm{m}_2=[.5,.5]$.}
\footnotetext{\textbf{Notation.} $\mathcal{P}_0$ denotes the initial population. For any real random vector $\bm{z}$ and any two real scalar $a<b$, $\bm{z} \sim$ uniform$[a, b]$ denotes that each entry of $\bm{z}$ is sampled uniformly from the interval $[a, b]$.}
\end{table}

\newpage
\subsubsection{Hyper-parameters for the Untrageted Image Adversarial Attacks.}
The hyper-parameters' value for each algorithm in the experiment of untargeted image adversarial attacks are selected in the same way as the targeted-image-attack experiments. The candidate and selected values are reported in Table \ref{imagenet-hyper-parameters-untargeted}.

\begin{table}[h]
\caption{Hyper-parameters for the Untargeted Attack on ImageNet (Table \ref{attack-robust-imagenet}). The hyper-parameter symbols for each algorithm are the same as their source publications. For example, $t_1$ in ZO-SLGHd and ZO-SLGHr denotes the initial scaling parameter, $\mu$ in ZO-AdaMM is the scaling parameter, and $\alpha$ denotes a constant learning rate. The set of candidate values that lead to the highest success rate and minimal perturbation size (averaged over the 10 image attackes) are be selected.}
\label{imagenet-hyper-parameters-untargeted}
\centering%
\begin{tabular}{ p{3cm}  p{6cm} | p{5cm} }
\midrule
 & Selected Values & Candidates $(\bm{\mu}^*)$ \\
\toprule
GS-PowerHP & $\alpha =2.0$, $N=8$, $\sigma_0=.1$, $b=0$.   &  $N\in\{5, 8\}$, $\alpha\in\{0.5, 1.0, 2.0\}$.       \\
GS-PowerGS & $\alpha =2.0$, $N=8$, $\sigma=.1$.   &  $N\in\{5, 8\}$, $\alpha\in\{0.5, 1.0, 2.0\}$.       \\

STD-Homotopy  & $\alpha=2.0$, $\gamma=.5$, $\sigma=0.1$, $T_{\bm{\mu}}=300$, $\tau=100$, $N_\sigma=10$.   &  $\gamma\in\{.5,.8\}$, $\alpha\in\{1.0, 2.0\}$.     \\

ZO-SLGHd   & $\beta=1/3072$, $\eta=.0001/3072$, $t_1=.1$, $\gamma=.995$      &   $\beta\in\{1.0, .01, 1/3072\}$.    \\
ZO-SLGHr   & $\beta=1/3072$, $t_1=.1$,  $\gamma=.995$.     &   $\beta\in\{.1, 1/3072\}$. \\
ZO-AdaMM   &$\beta_{1}=.9$, $\beta_2=.1$, $\alpha=.1$, $\mu=.1$.        &   $\alpha\in\{.1, 1.0\}$, $\beta_1\in\{.5,.9\}$, $\beta_2\in\{.1,.3\}$.   \\
ZO-SGD   & $\alpha=0.001$, $\mu=.1$.      &   $\alpha\in\{ .01/3072, .001, .01, .1, 2.0\}$.       \\
Square Attack & $\varepsilon = 30$, $p=0.3$. & $p\in\{0.1, 0.3\}$.\\
\bottomrule
\end{tabular}
\footnotetext{Total number of generations: 1000. $\bm{m}_1=[-.5,-.5]$. $\bm{m}_2=[.5,.5]$.}
\footnotetext{\textbf{Notation.} $\mathcal{P}_0$ denotes the initial population. For any real random vector $\bm{z}$ and any two real scalar $a<b$, $\bm{z} \sim$ uniform$[a, b]$ denotes that each entry of $\bm{z}$ is sampled uniformly from the interval $[a, b]$.}
\end{table}

\newpage
\subsection{Example Adversarial Images (Targeted)}
\label{appendix-adversarial}
In this section we produce the adversarial images for 4 randomly picked images from each image set, using the top three smoothing-based algorithms in our MNIST experiment, CMA-ES, and Square Attack, respectively. For ImageNet, we do not involve CMA-ES because it requires a large RAM which is beyond our computer's capacity. The results can be found in Table \ref{Adversarial-Image-MNIST}, \ref{Adversarial-Image-CIFAR}, and \ref{Adversarial-Image-ImageNet}.

\onecolumn
\begin{table}[h]
\centering
\caption{Four adversarial images for the least-likely targeted attacks on MNIST, produced by the top three smoothing-based algorithms in our MNIST experiment, CMA-ES, and Square Attack, respectively.}
\label{Adversarial-Image-MNIST}
\begin{tabular}{|c|c|c|c|c|}
\hline
\textbf{Test Image ID} & 9953 & 3850 & 4962 & 3886 \\ \hline
\textbf{Original Image} &
\includegraphics[scale=0.18]{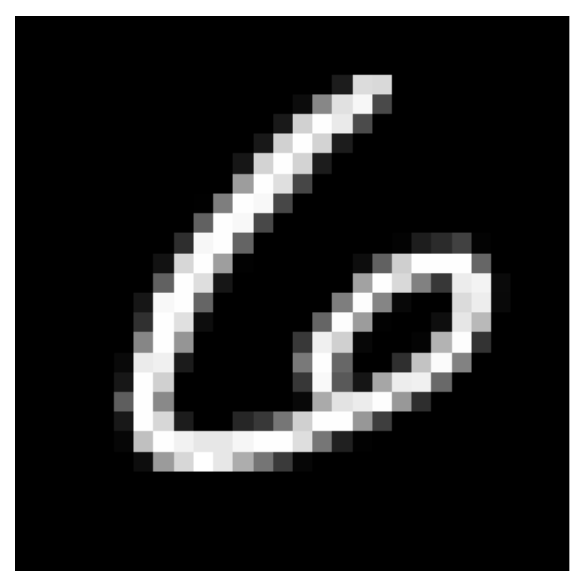} &
\includegraphics[scale=0.18]{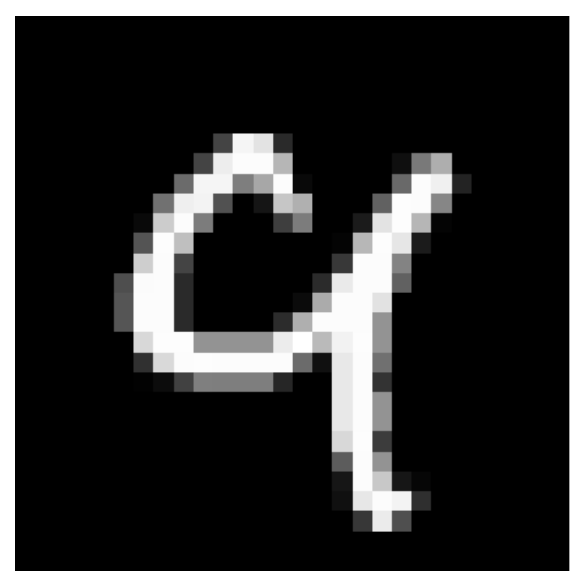} &
\includegraphics[scale=0.18]{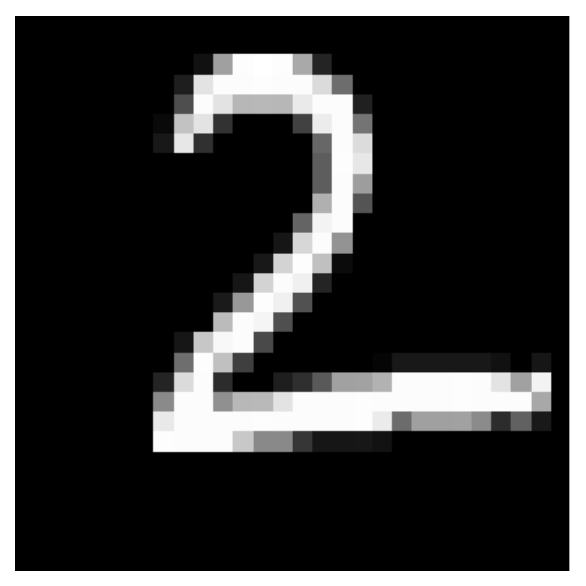} &
\includegraphics[scale=0.18]{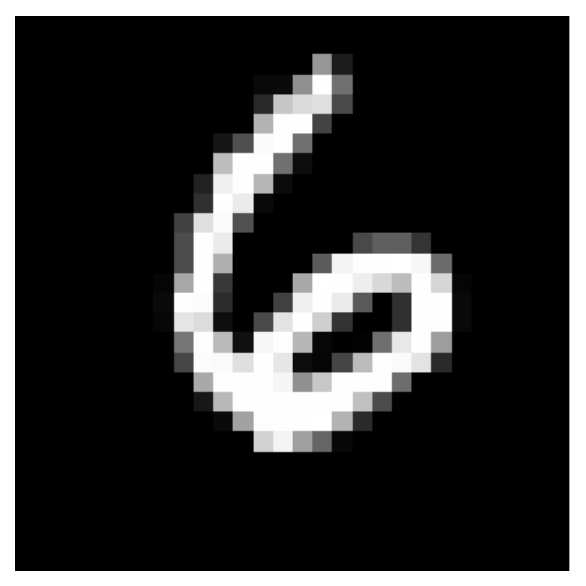} \\ 
\textbf{Target Label $\mathcal{T}$} & 9 & 1 & 9 & 1 \\ 
\hline
\textbf{Our Algo.} &
\includegraphics[scale=0.18]{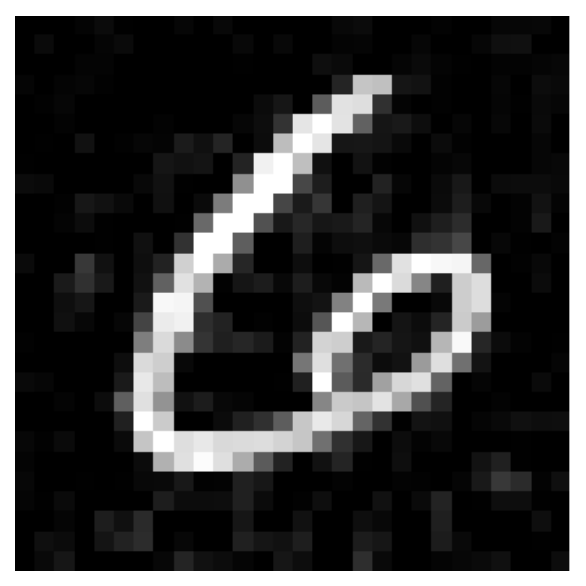} &
\includegraphics[scale=0.18]{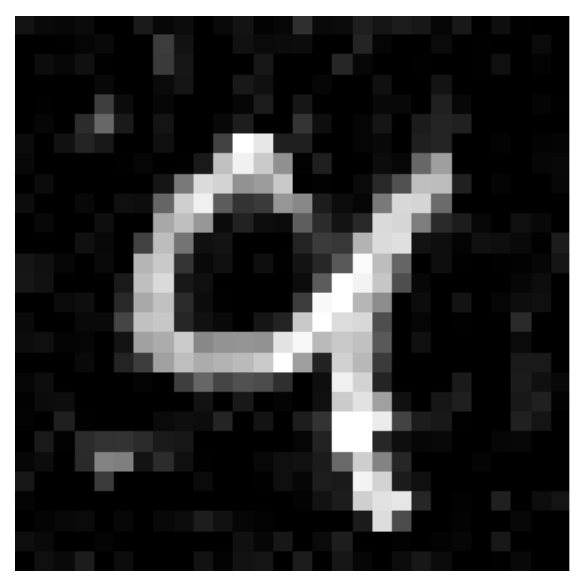} &
\includegraphics[scale=0.18]{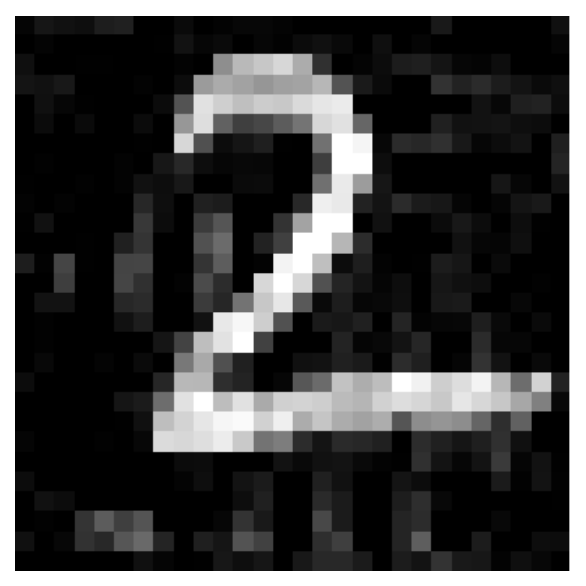} &
\includegraphics[scale=0.18]{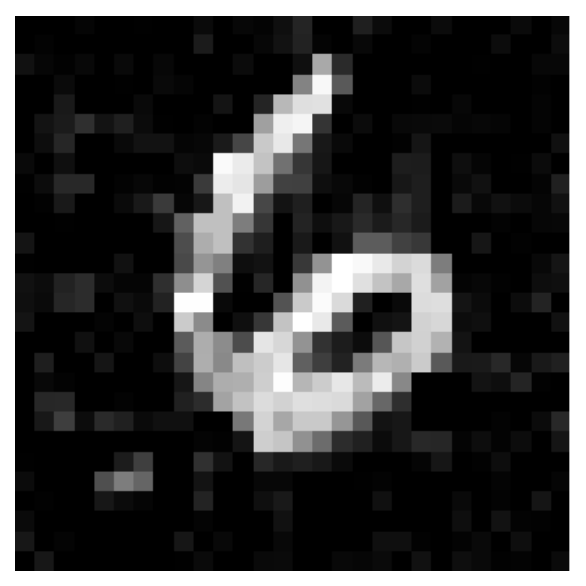} \\ 
$R^2(\bm{a},\bm{a+\mu^*})$ & $94.5\%$ & $88.3\%$ & $85.1\%$ & $89.3\%$ \\ 
\hline
\textbf{CMA-ES} &
\includegraphics[scale=0.18]{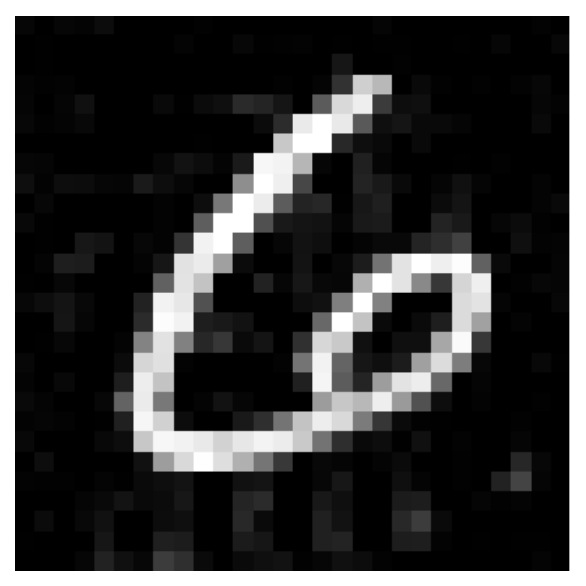} &
\includegraphics[scale=0.18]{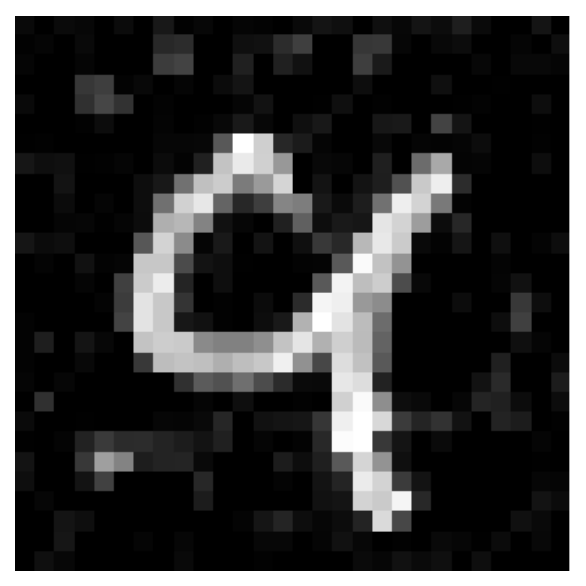} &
\includegraphics[scale=0.18]{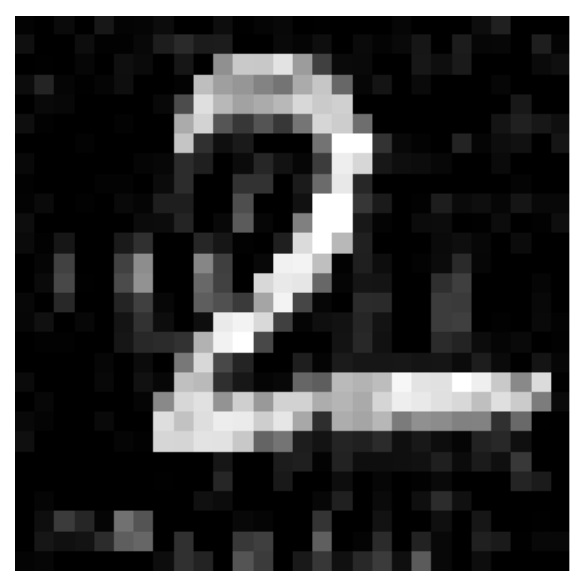} &
\includegraphics[scale=0.18]{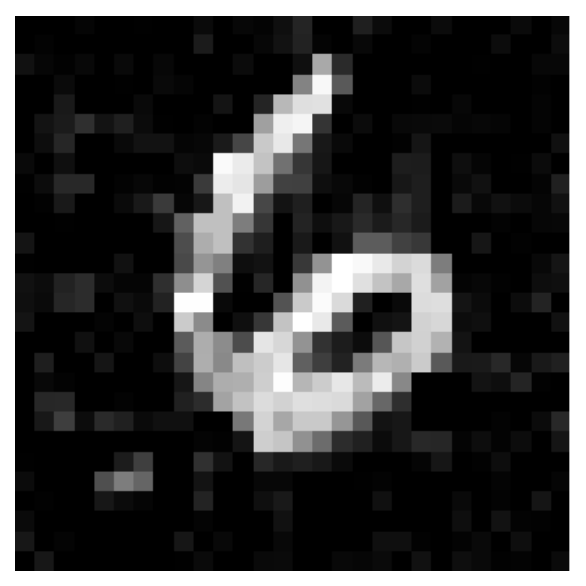} \\ 
$R^2(\bm{a},\bm{a+\mu^*})$ & $95.6\%$ & $88.2\%$ & $82.1\%$ & $89.3\%$ \\ 
\hline
\textbf{GS-PowerOpt} &
\includegraphics[scale=0.18]{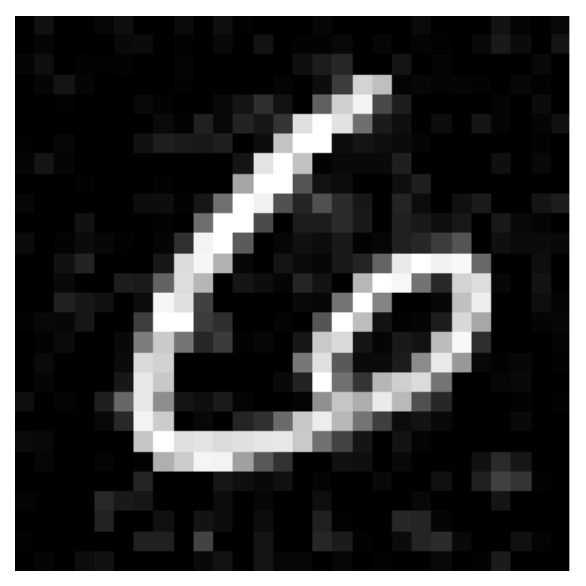} &
\includegraphics[scale=0.18]{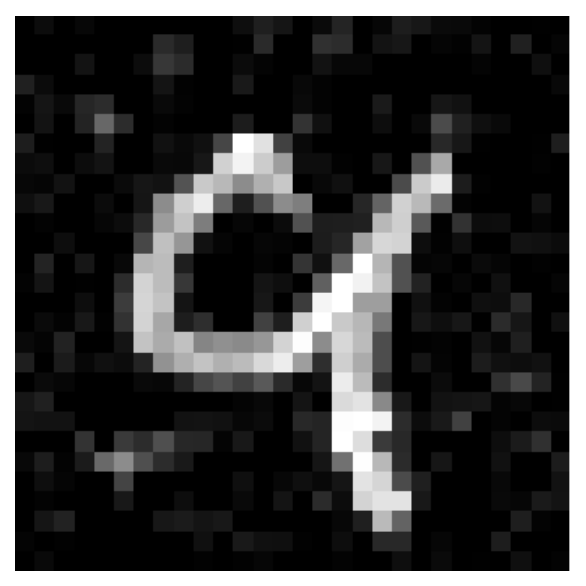} &
\includegraphics[scale=0.18]{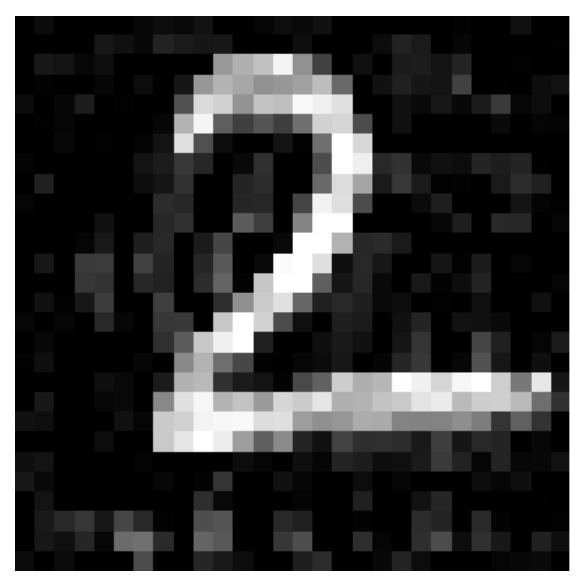} &
\includegraphics[scale=0.18]{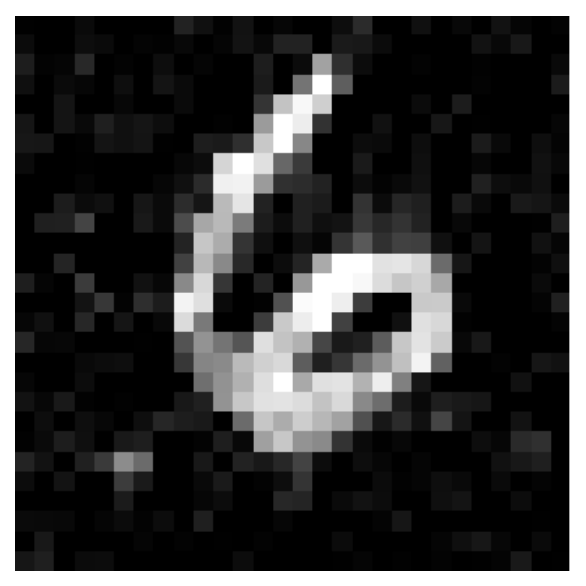} \\ 
$R^2(\bm{a},\bm{a+\mu^*})$ & $93.7\%$ & $86.0\%$ & $83.3\%$ & $86.7\%$ \\ 
\hline
\textbf{ZOSGD} &
\includegraphics[scale=0.18]{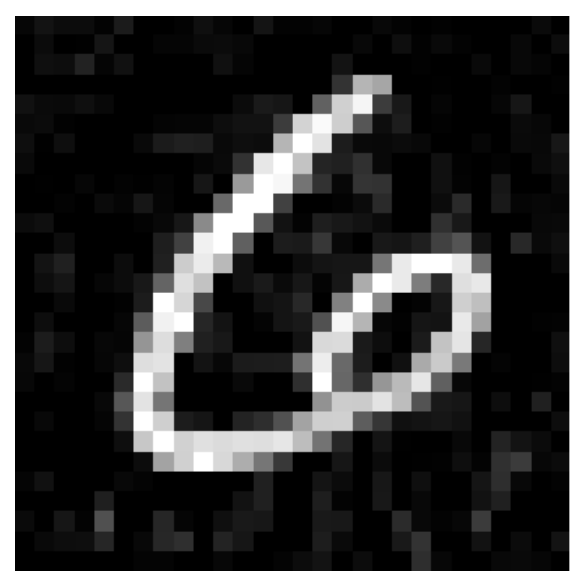} &
\includegraphics[scale=0.18]{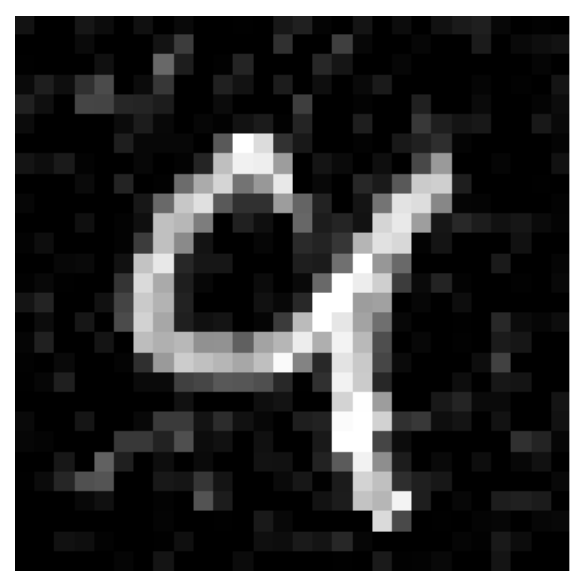} &
\includegraphics[scale=0.18]{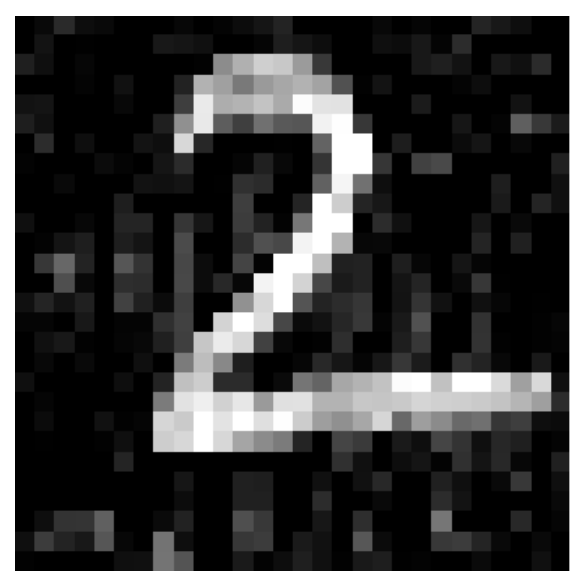} &
\includegraphics[scale=0.18]{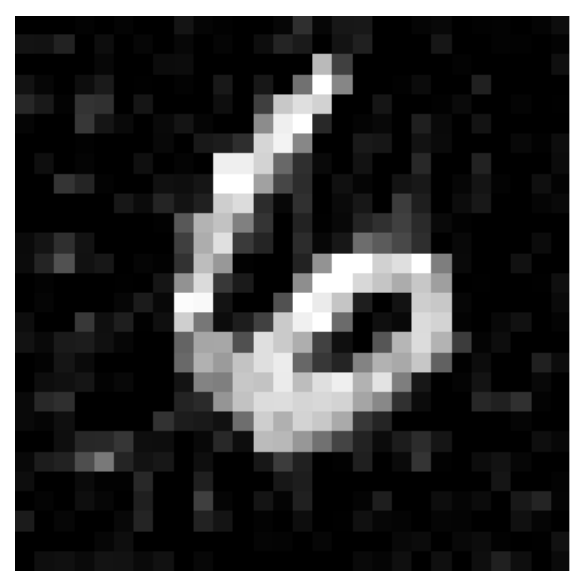} \\ 
$R^2(\bm{a},\bm{a+\mu^*})$ & $92.7\%$ & $85.8\%$ & $81.1\%$ & $87.5\%$ \\ 
\hline
\textbf{Square Attack} &
\includegraphics[scale=0.18]{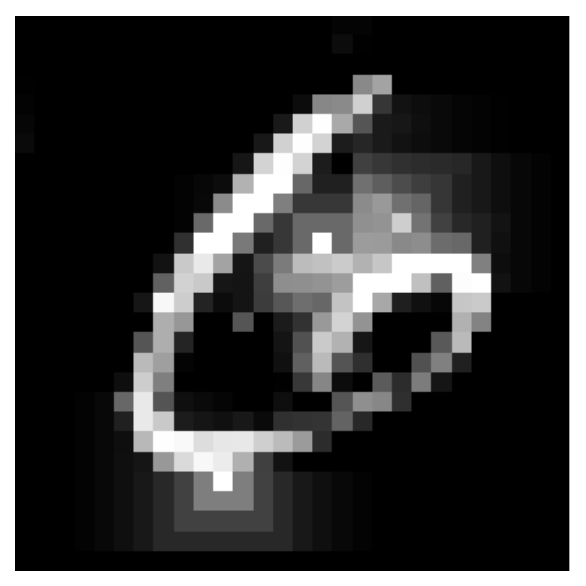} &
\includegraphics[scale=0.18]{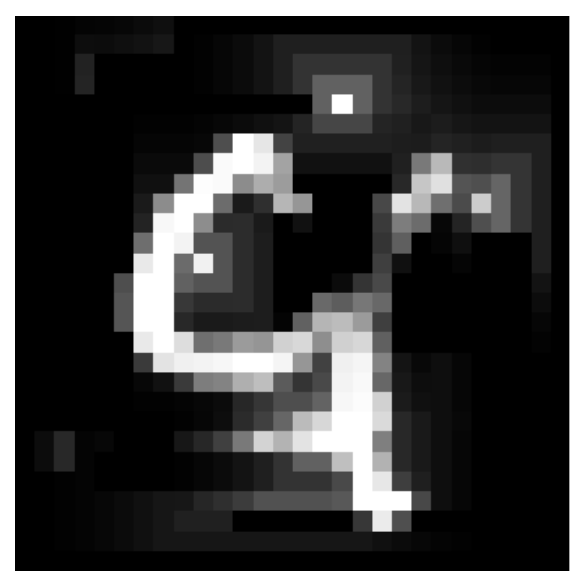} &
\includegraphics[scale=0.18]{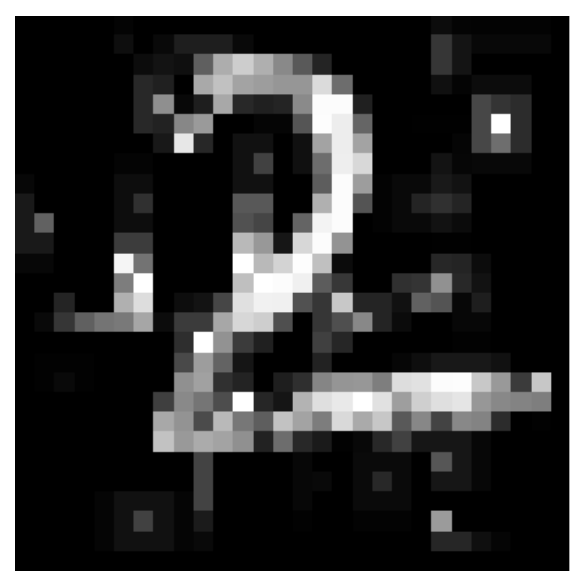} &
\includegraphics[scale=0.18]{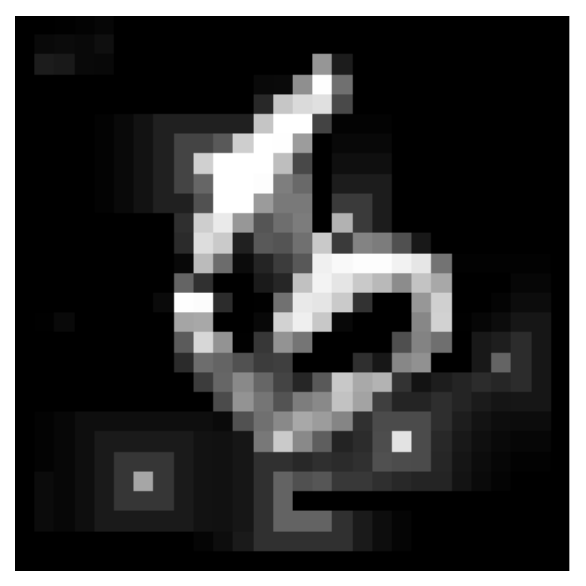} \\ 
$R^2(\bm{a},\bm{a+\mu^*})$ & $65.4\%$ & $66.0\%$ & $67.4\%$ & $67.5\%$ \\ 
\hline
\end{tabular}
\end{table}

\begin{table}[h!]
\centering
\caption{Four adversarial images for the least-likely targeted attacks on CIFAR-10, produced by the top three smoothing-based algorithms in our CIFAR-10 experiment, CMA-ES, and Square Attack, respectively. Unsuccessful attacks (i.e., predicted label is different from the adversarial target) are marked with `Unsuccessful'.}
\label{Adversarial-Image-CIFAR}
\begin{tabular}{|c|c|c|c|c|}
\hline
\textbf{Test Image ID} & 9953 & 3850 & 4962 & 3886 \\ \hline
\textbf{Original Image} &
\includegraphics[scale=0.15]{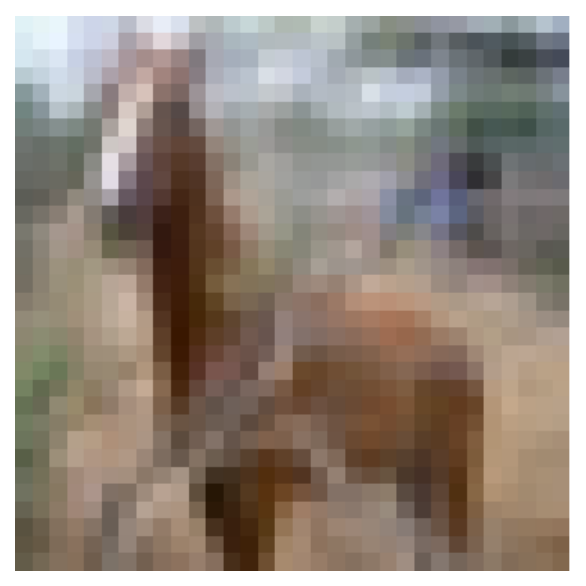} &
\includegraphics[scale=0.15]{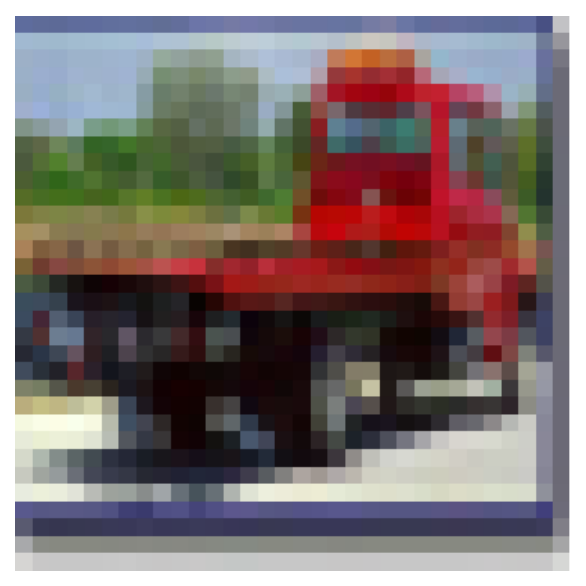} &
\includegraphics[scale=0.15]{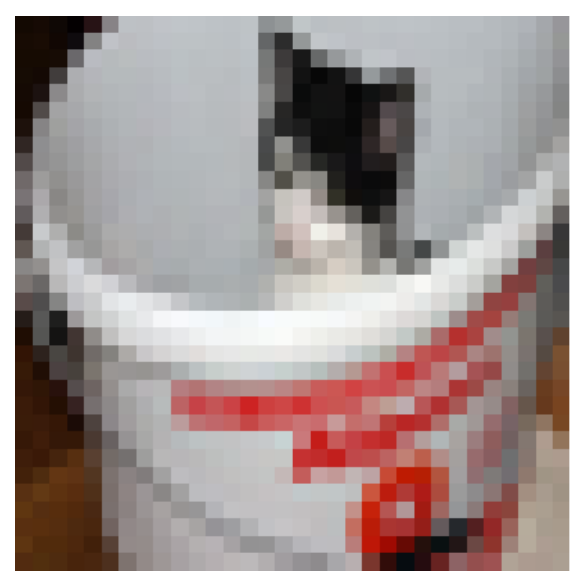} &
\includegraphics[scale=0.15]{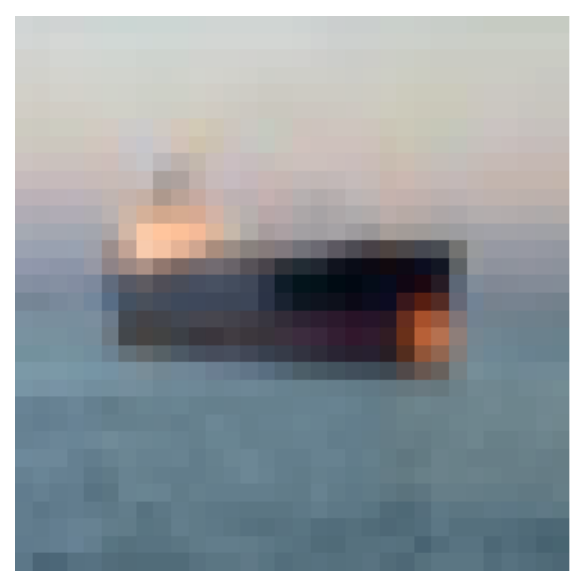} \\ 
\textbf{True Label} & 7 (i.e., Horse) & 9 (i.e., Truck) & 3 (i.e., Cat) & 8 (i.e., Ship) \\ 
\textbf{Target Label $\mathcal{T}$} & 1(i.e., Automobile) & 4 (i.e., Deer) & 7 (i.e., Horse) & 7 (i.e., Horse) \\ 
\hline
\textbf{Our Algo.} &
\includegraphics[scale=0.15]{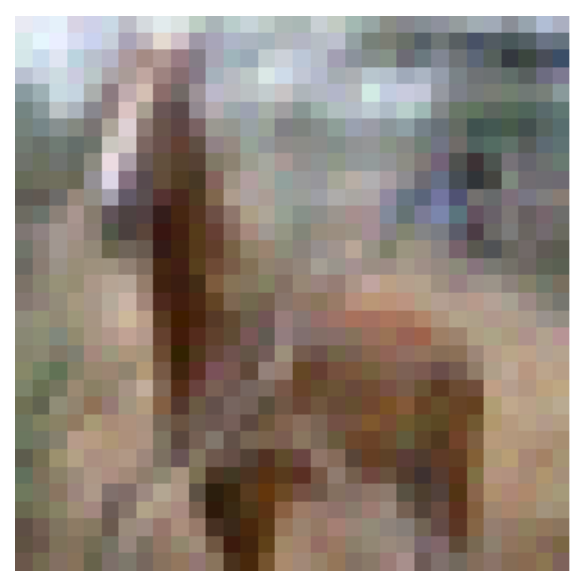} &
\includegraphics[scale=0.15]{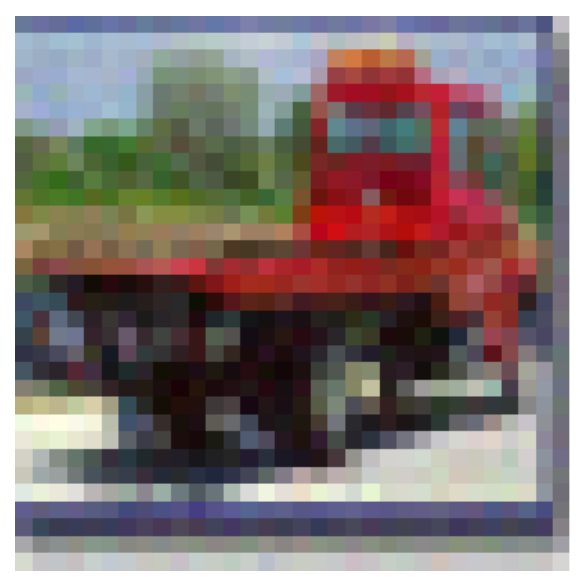} &
\includegraphics[scale=0.15]{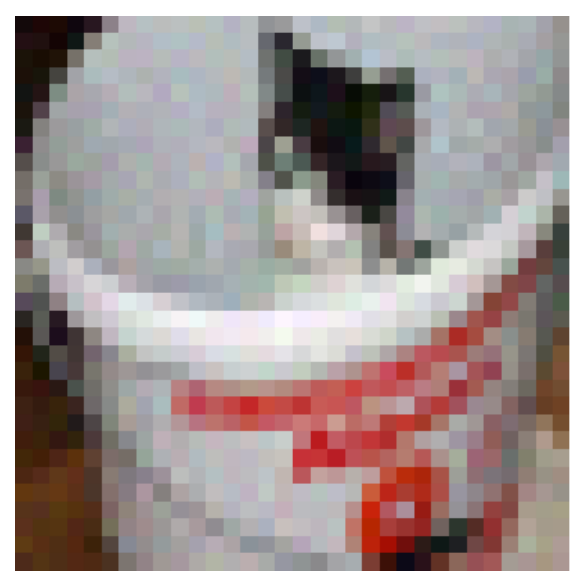} &
\includegraphics[scale=0.15]{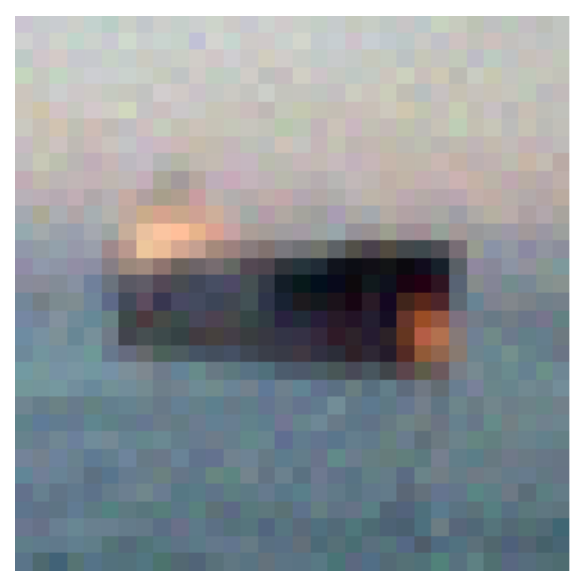} \\ 
$R^2(\bm{a},\bm{a+\mu^*})$ & $98.9\%$ & $99.2\%$ & $99.1\%$ & $98.7\%$ \\ 
\hline
\textbf{ZO-SLGHd} &
\includegraphics[scale=0.15]{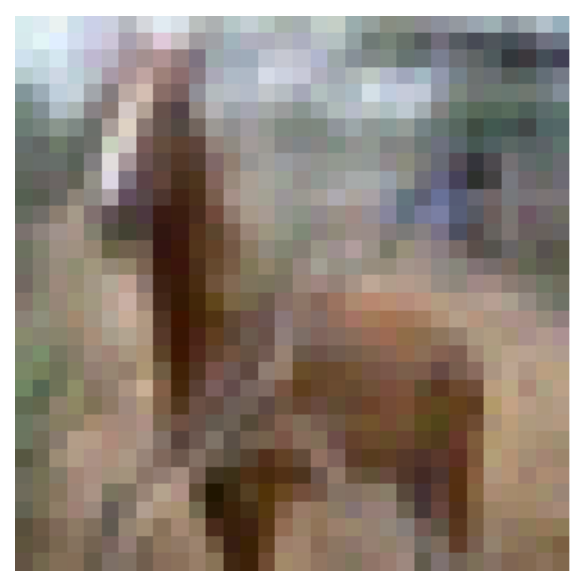} &
\includegraphics[scale=0.15]{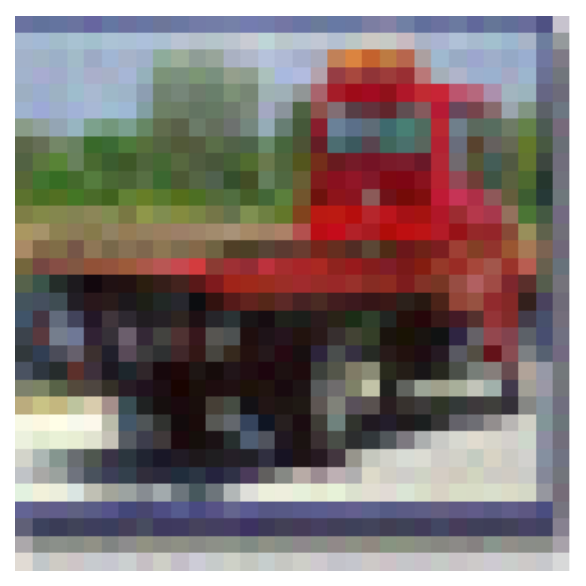} &
\includegraphics[scale=0.15]{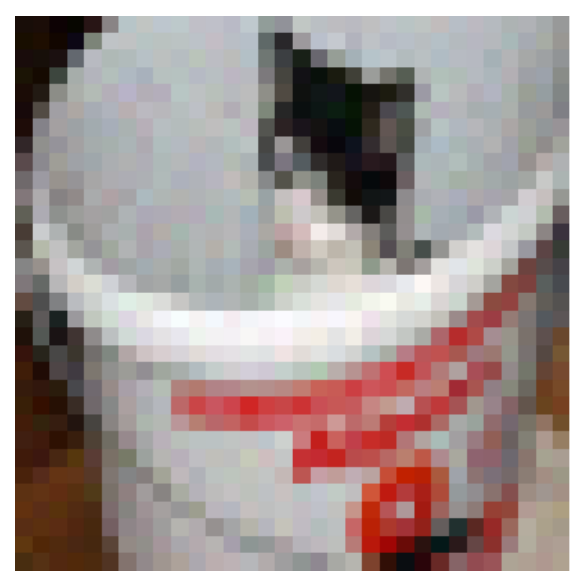} &
\includegraphics[scale=0.15]{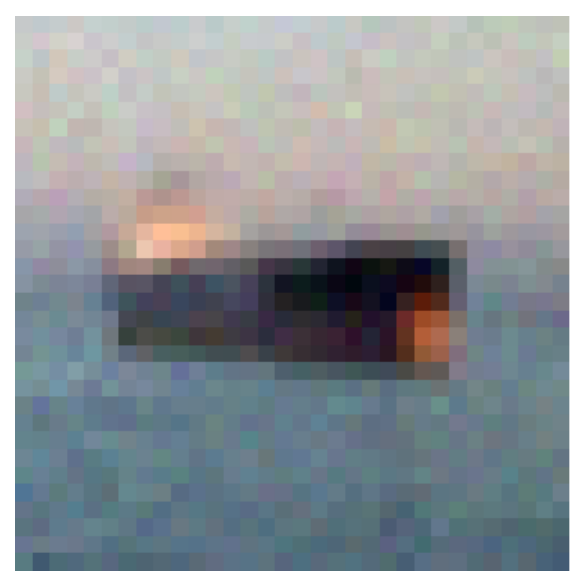} \\ 
$R^2(\bm{a},\bm{a+\mu^*})$ & $99.3\%$ & $99.4\%$ & $99.5\%$ & $98.7\%$ \\ 
\hline
\textbf{ZO-SLGHr} &
\includegraphics[scale=0.15]{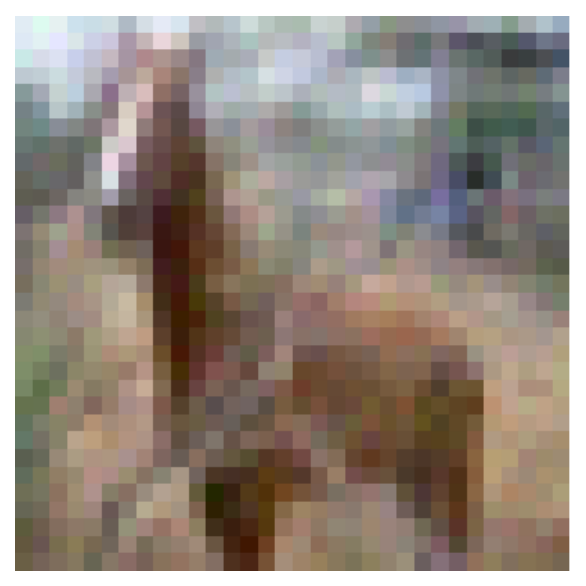} &
\includegraphics[scale=0.15]{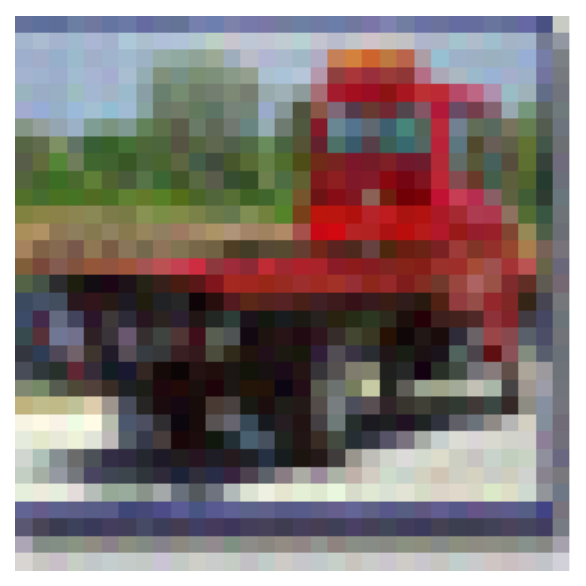} &
\includegraphics[scale=0.15]{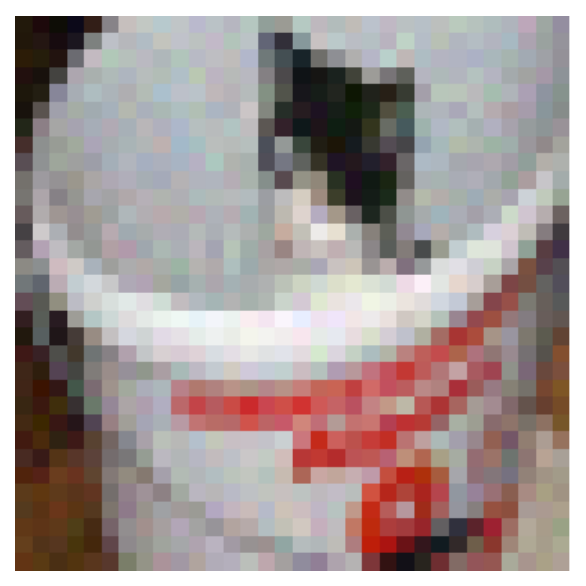} &
\includegraphics[scale=0.15]{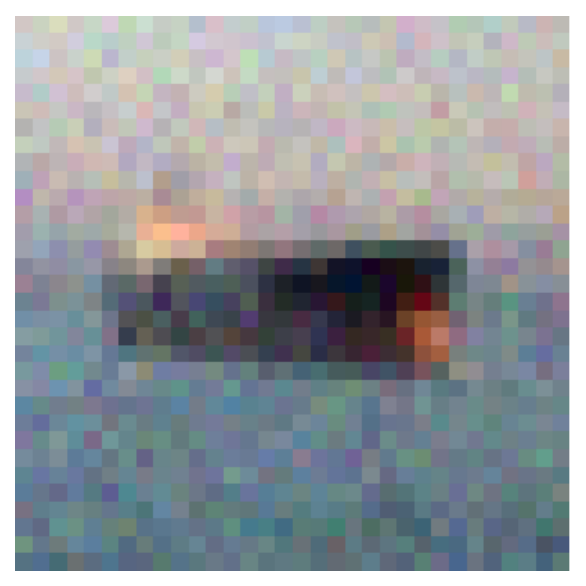} \\ 
$R^2(\bm{a},\bm{a+\mu^*})$ & $98.4\%$ & $99.2\%$ & $98.8\%$ & $95.1\%$ \\ 
\hline
\textbf{CMA-ES} &
\includegraphics[scale=0.15]{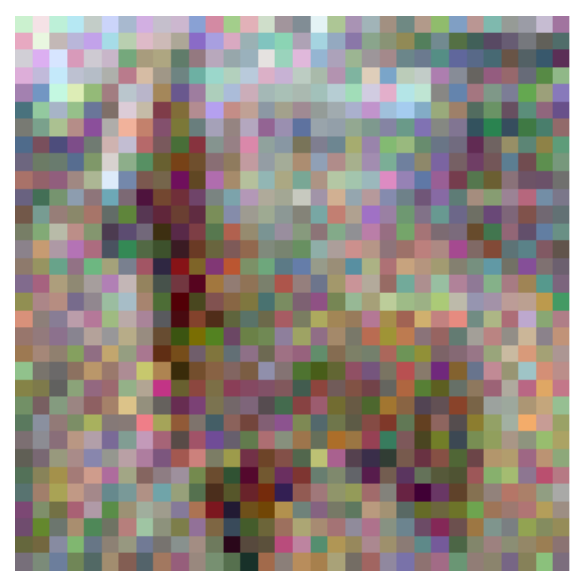} &
\includegraphics[scale=0.15]{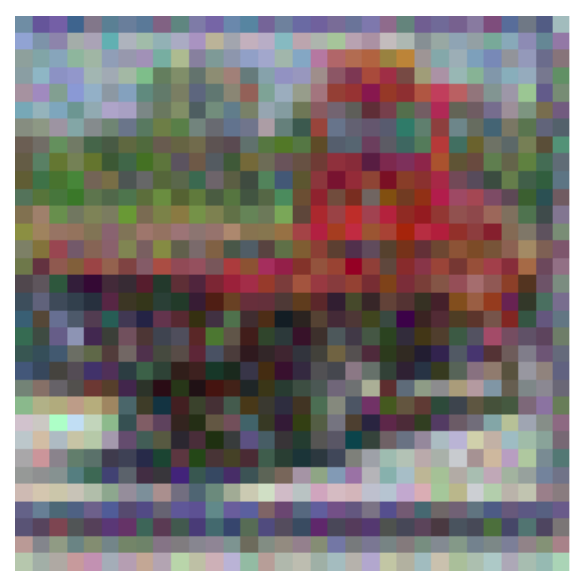} &
\includegraphics[scale=0.15]{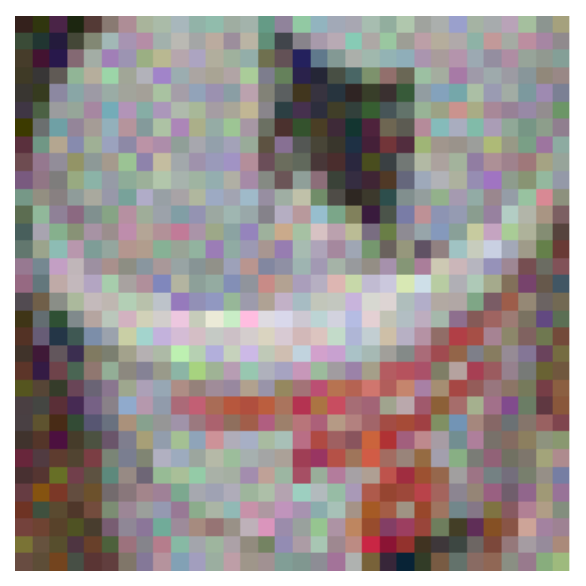} &
\includegraphics[scale=0.15]{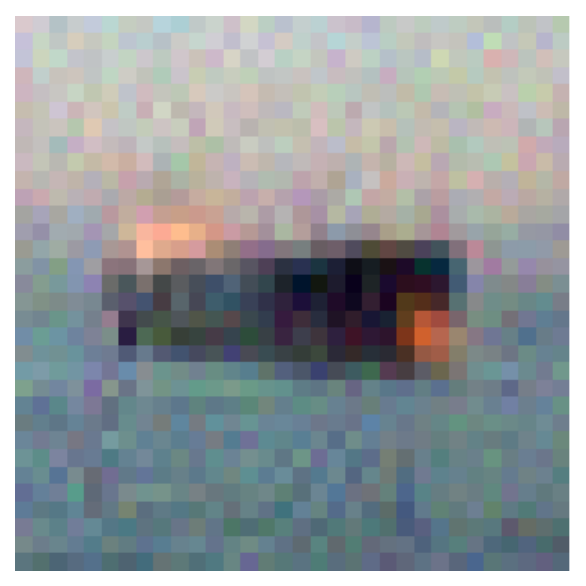} \\ 
$R^2(\bm{a},\bm{a+\mu^*})$ & $56.2\%$ & $87.9\%$ & $86.1\%$ & $95.5\%$ \\ 
\hline
\textbf{Square Attack} &
\includegraphics[scale=0.15]{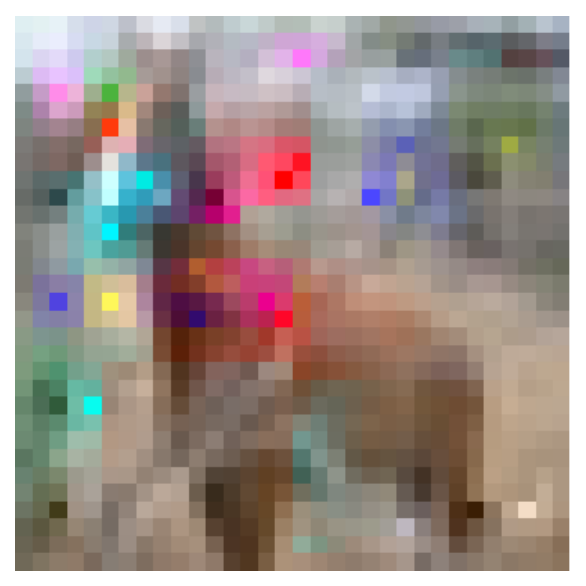} &
\includegraphics[scale=0.15]{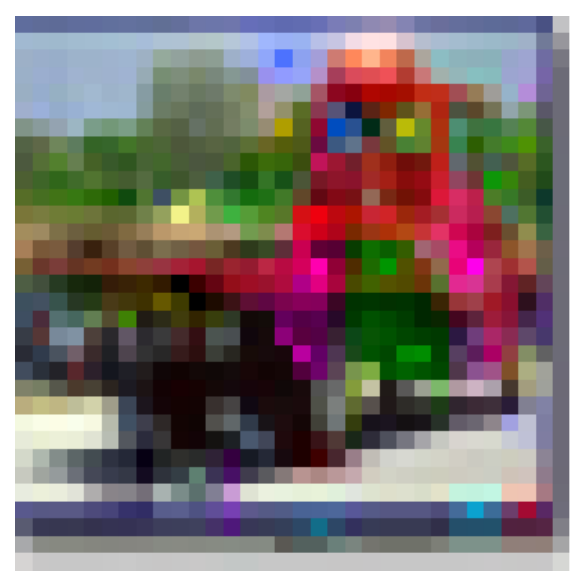} &
\includegraphics[scale=0.15]{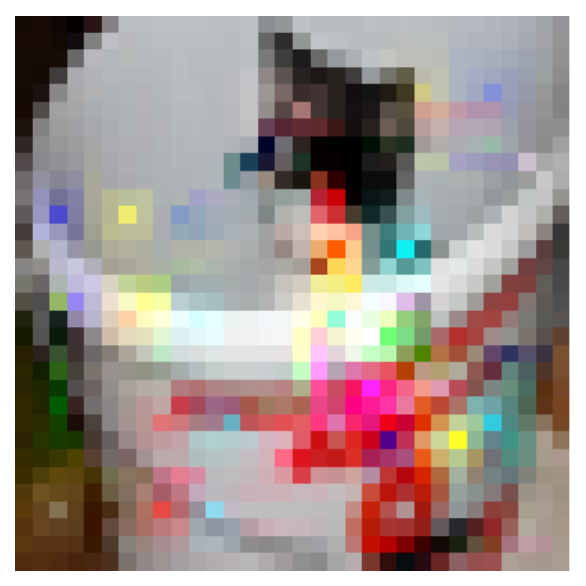} &
\includegraphics[scale=0.15]{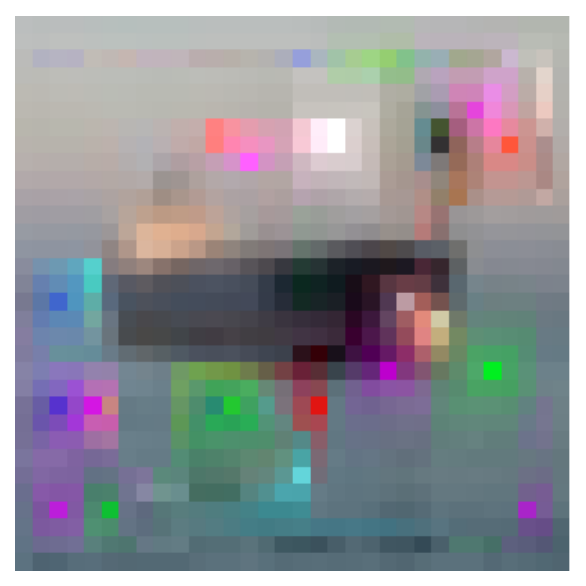} \\ 
$R^2(\bm{a},\bm{a+\mu^*})$ & $65.8\%$ &86.6\% & $85.9\%$ & $68.6\%$ \\ 
\hline
\end{tabular}
\end{table}

\begin{table}[h!]
\centering
\caption{Four adversarial images for the least-likely targeted attacks on ImageNet, produced by the top four smoothing-based algorithms in our CIFAR-10 experiment and Square Attack, respectively. Unsuccessful attacks (i.e., predicted label is different from the adversarial target) are marked with `Unsuccessful'.}
\label{Adversarial-Image-ImageNet}
\begin{tabular}{|c|c|c|c|c|}
\hline
\textbf{Original Image} &
\includegraphics[scale=0.18]{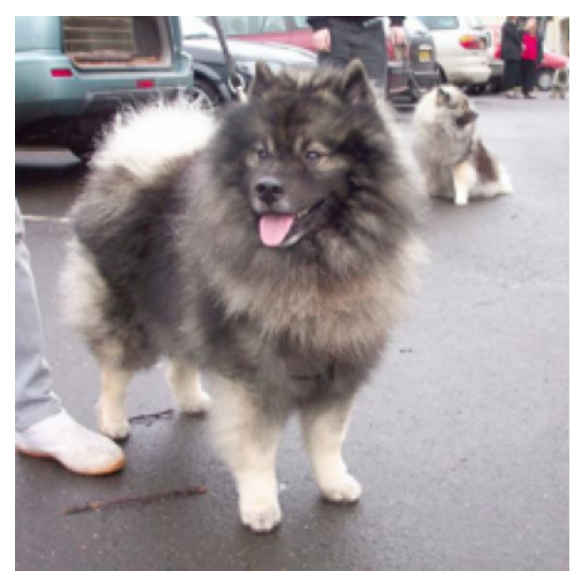} &
\includegraphics[scale=0.18]{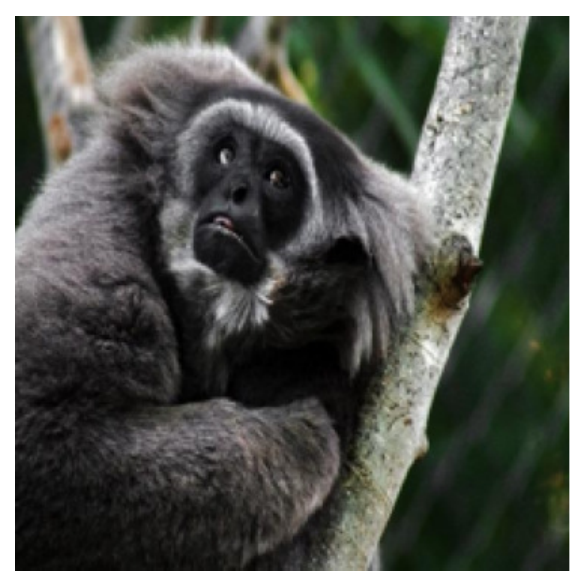}&
\includegraphics[scale=0.18]{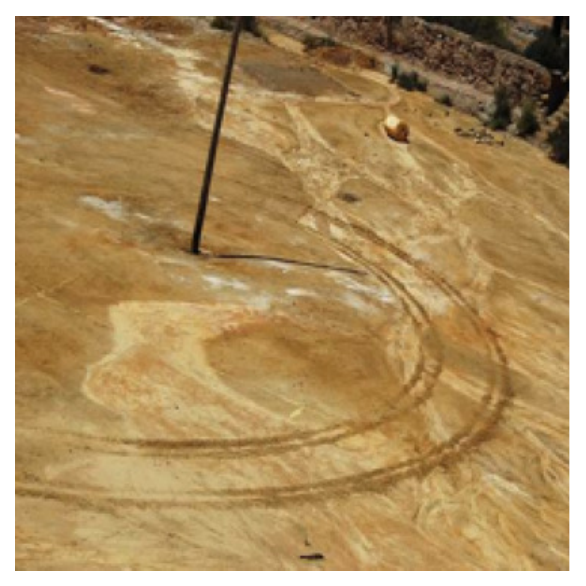}  &
\includegraphics[scale=0.18]{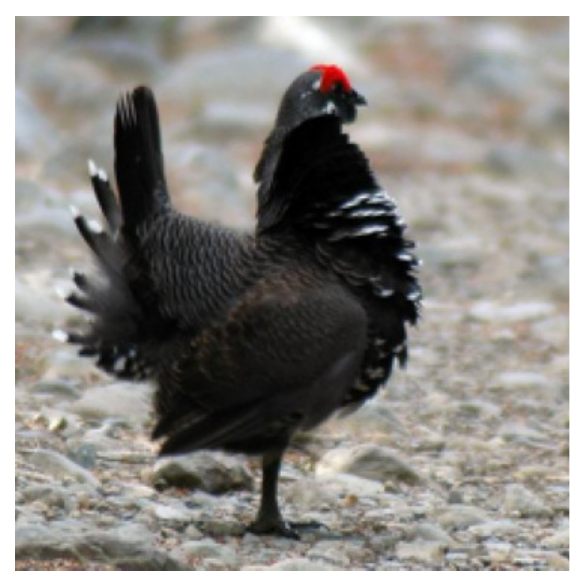} \\ 
\textbf{True Label} & 261 & 368  & 835 & 80 \\ 
\textbf{Target Label $\mathcal{T}$} & 344 & 4 & 487  & 95  \\ 
\hline
\textbf{Our Algo.} &
\includegraphics[scale=0.18]{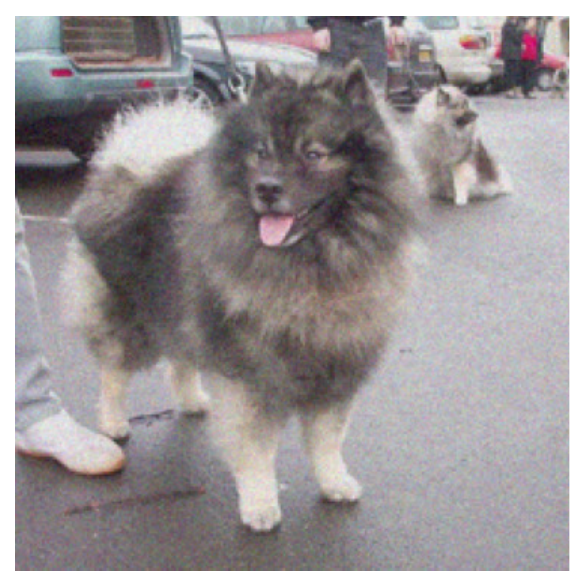} &
\includegraphics[scale=0.18]{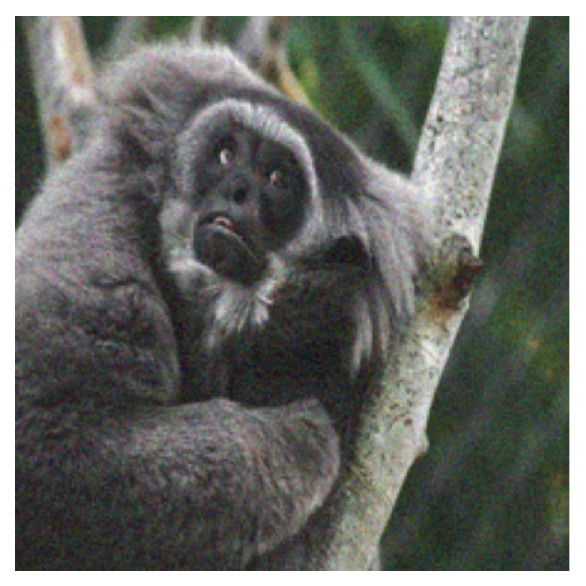} &
\includegraphics[scale=0.18]{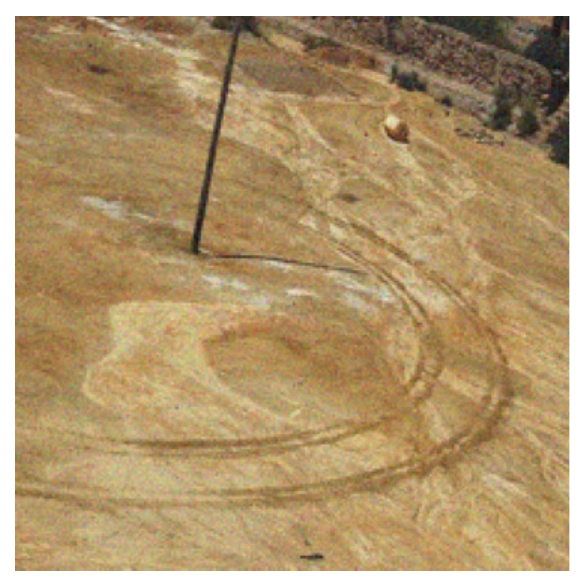} &
\includegraphics[scale=0.18]{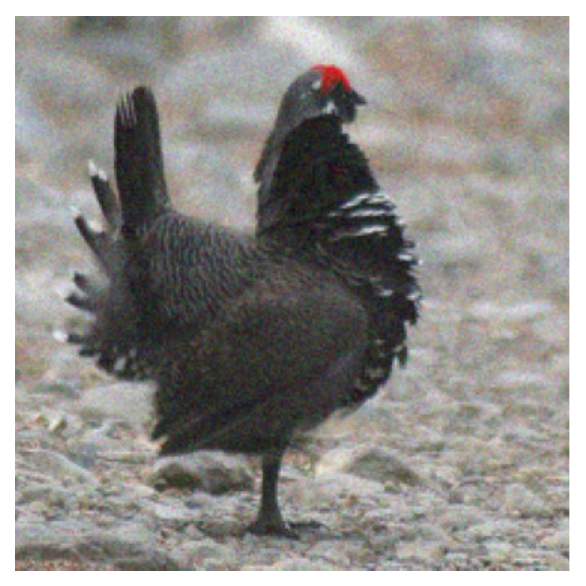} \\ 
$R^2(\bm{a},\bm{a+\mu^*})$ & $96.6\%$ & $95.3\%$ & $94.0\%$ & $95.8\%$ \\ 
\hline
\textbf{GS-PowerOpt.} &
\includegraphics[scale=0.18]{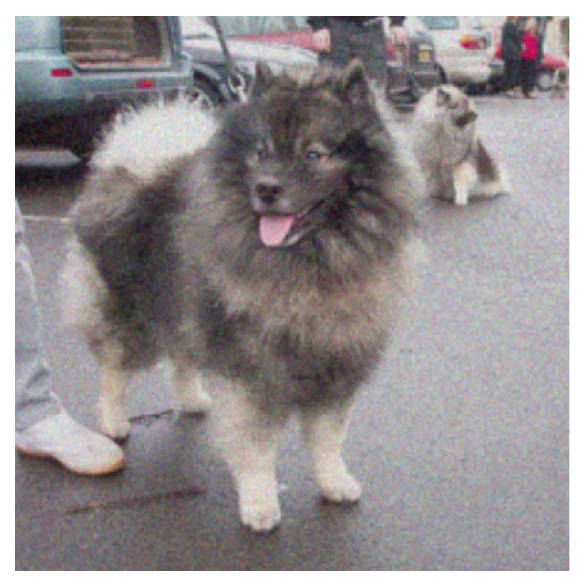} &
\includegraphics[scale=0.18]{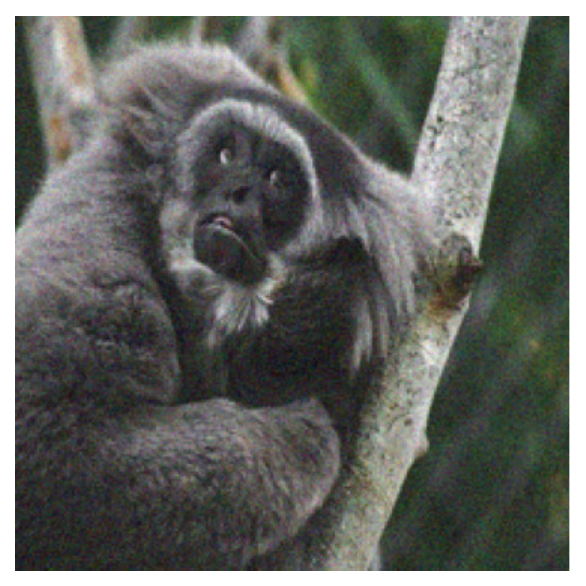}&
\includegraphics[scale=0.18]{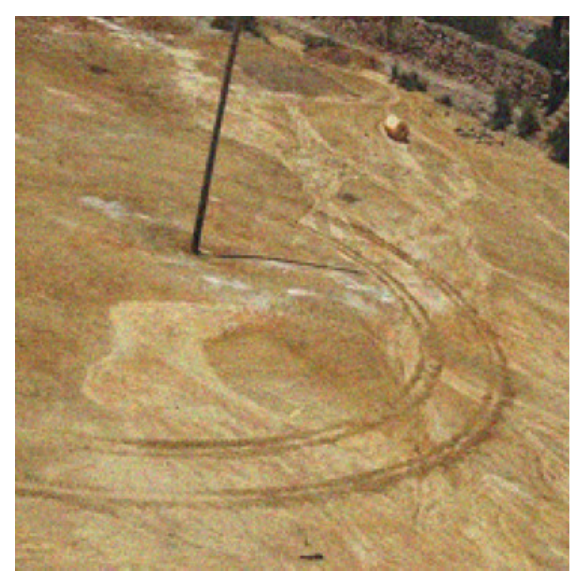} &
Unsuccessful. \\ 
$R^2(\bm{a},\bm{a+\mu^*})$ & $95.2\%$ & $94.8\%$  & $91.8\%$& NA. \\ 
\hline
\textbf{ZOSLGHd} &
\includegraphics[scale=0.18]{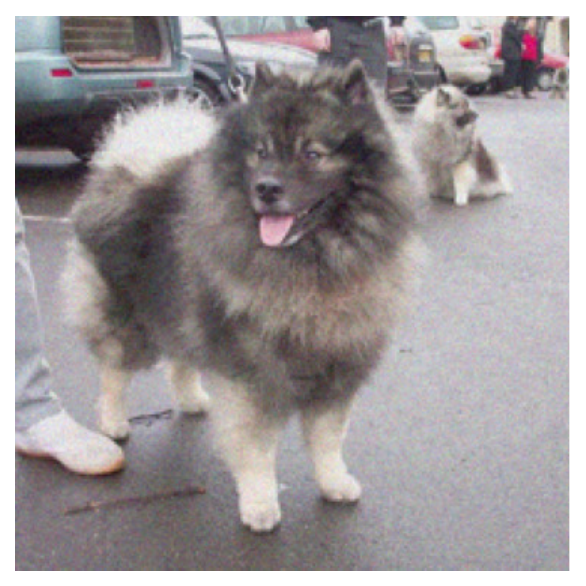} &
\includegraphics[scale=0.18]{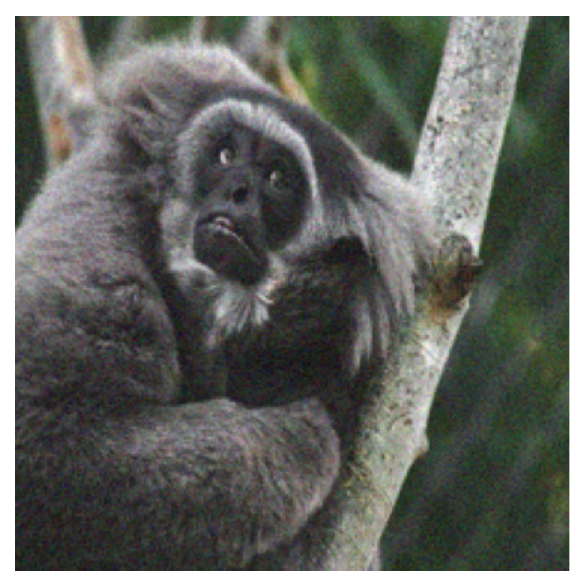}&
\includegraphics[scale=0.18]{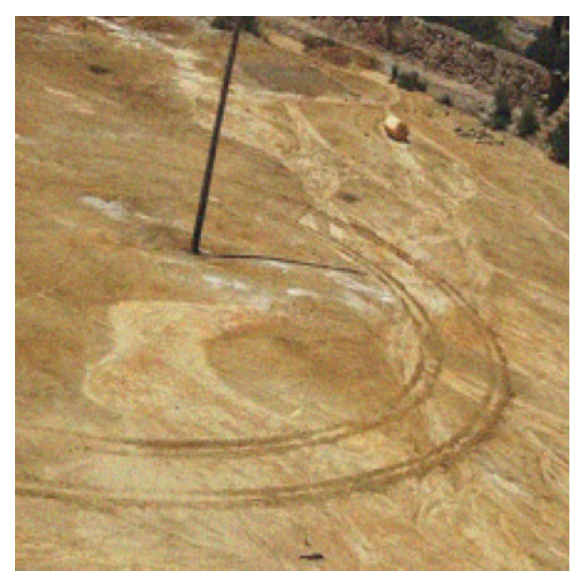} &
Unsuccessful. \\
$R^2(\bm{a},\bm{a+\mu^*})$ & $97.7\%$ & $95.9\%$ & $96.1\%$ & NA. \\ 
\hline
\textbf{ZOSLGHr} &
\includegraphics[scale=0.18]{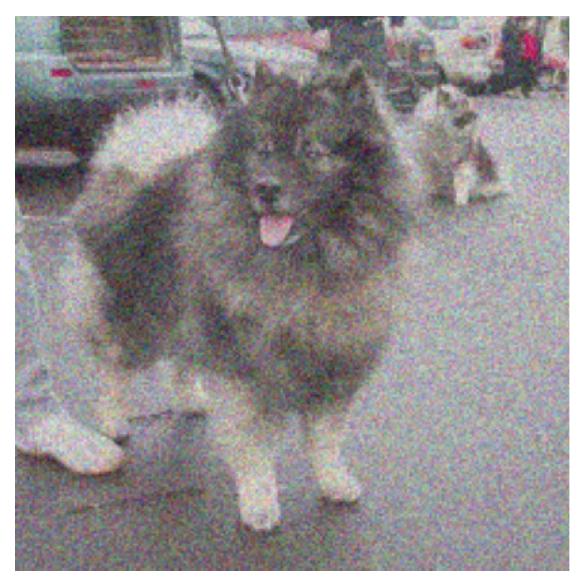} &
Unsuccessful.&
\includegraphics[scale=0.18]{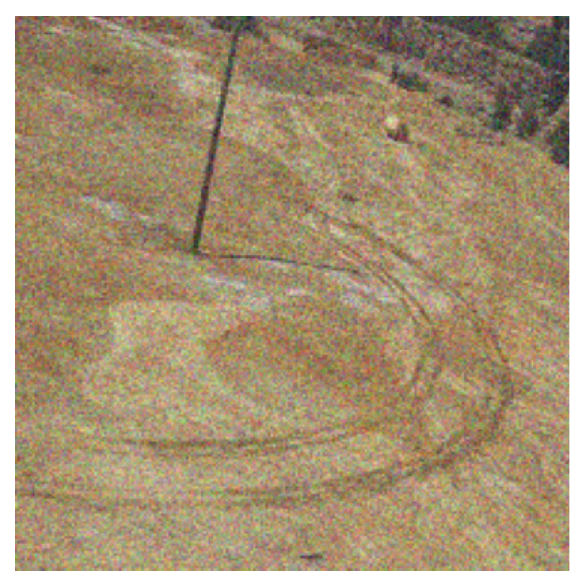} &
Unsuccessful. \\
$R^2(\bm{a},\bm{a+\mu^*})$ & $69.5\%$ & NA. & $48.4\%$ & NA. \\ 
\hline
\textbf{Square Attack} &
\includegraphics[scale=0.18]{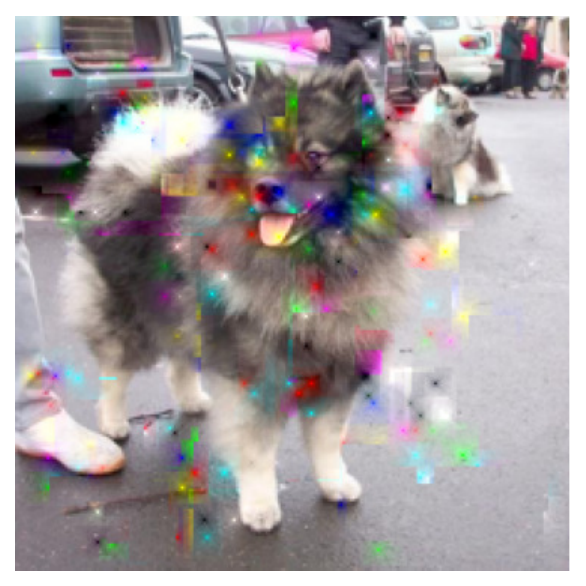} &
\includegraphics[scale=0.18]{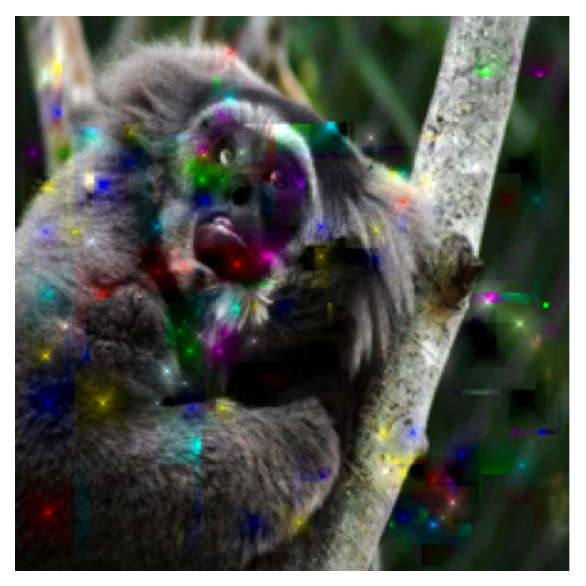}&
\includegraphics[scale=0.18]{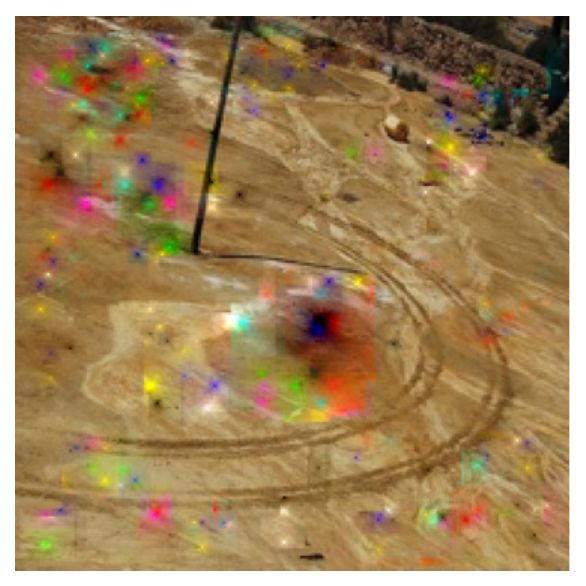} &
\includegraphics[scale=0.18]{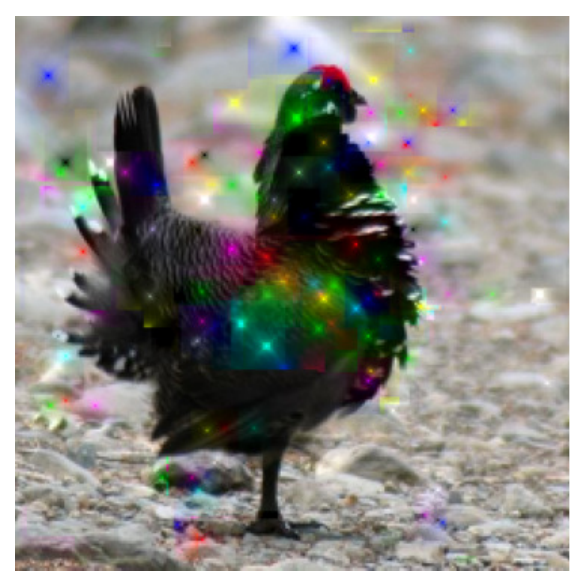} \\
$R^2(\bm{a},\bm{a+\mu^*})$ & $91.6\%$ & $91.2\%$ & $85.9\%$ & $93.8\%$ \\ 
\hline
\end{tabular}
\end{table}
\end{document}